\documentclass[12pt,a4paper]{book}
\input amssym
\catcode`\@=11

\font\decp=cmr6
\font\titl=cmr17 at 20pt
\font\capit=cmbx12 at 24pt
\font\nom=cmr17
\font\dedi=cmsl12 at24pt

\font\tenmsa=msam10 \font\sevenmsa=msam7 \font\fivemsa=msam5
\font\tenmsb=msbm10 at 12pt\font\sevenmsb=msbm7 at 9pt\font\fivemsb=msbm5 at 7pt
\newfam\msafam \newfam\msbfam
\textfont\msafam=\tenmsa  \scriptfont\msafam=\sevenmsa \scriptscriptfont\msafam=
\fivemsa
\textfont\msbfam=\tenmsb  \scriptfont\msbfam=\sevenmsb \scriptscriptfont\msbfam=
\fivemsb

\font\teneufm=eufm10 \font\seveneufm=eufm7 \font\fiveeufm=eufm5
\newfam\eufmfam 
\textfont\eufmfam=\teneufm \scriptfont\eufmfam=\seveneufm\scriptscriptfont
\eufmfam=\fiveeufm

\def\hexnumber@#1{\ifcase#1 0\or1\or2\or3\or4\or5\or6\or7\or8\or9\or
 A\or B\or C\or D\or E\or F\fi }

\edef\msa@{\hexnumber@\msafam} \edef\msb@{\hexnumber@\msbfam}

\mathchardef\qed="0\msa@03

\def\Bbb#1{{\fam\msbfam\relax#1}}

\def\R{{\Bbb R}} \def\C{{\Bbb C}}  \def\N{{\Bbb N}} 

\catcode`\@=12

\newtheorem{teo}{Theorem}[section]
\newtheorem{pro}[teo]{Proposition}
\newtheorem{cor}[teo]{Corollary} 
\newtheorem{eje}[teo]{Example} 
\newtheorem{lem}[teo]{Lemma} 
\newtheorem{obs}[teo]{Remark} 
\newtheorem{defin}[teo]{Definition}

\newcommand{\bteo}{\begin{teo}}
\newcommand{\eteo}{\end{teo}}
\newcommand{\bpro}{\begin{pro}}
\newcommand{\epro}{\end{pro}}
\newcommand{\bcor}{\begin{cor}}
\newcommand{\ecor}{\end{cor}}
\newcommand{\blem}{\begin{lem}}
\newcommand{\elem}{\end{lem}}
\newcommand{\bdefi}{\begin{defin}\rm}
\newcommand{\edefi}{\end{defin}}
\newcommand{\beje}{\begin{eje}\rm}
\newcommand{\eeje}{\end{eje}}
\newcommand{\bobs}{\begin{obs}\rm}
\newcommand{\eobs}{\end{obs}}
\newcommand{\bequ}{\begin{equation}}
\newcommand{\eequ}{\end{equation}}

\def\bdem{{\sl Proof. }}
\def\edem{\par\bigskip}

\def\S{{\cal S}}
\def\SS{{\cal S}_0}
\def\LL{L^\infty_0}
\def\M{{\cal M}}
\def\NN{\hbox{N}}

\def\loX{{\Lambda^p_X(w)}}
\def\loinX{{\Lambda^{p,\infty}_X(w)}}
\def\lox#1#2#3{{\Lambda^{#1}_{#2}(#3)}}
\def\loinx#1#2#3{{\Lambda^{#1,\infty}_{#2}(#3)}}

\def\lo{{\Lambda^p_u(w)}}
\def\loz{{\Lambda^{p_0}_{u_0}(w_0)}}
\def\lou{{\Lambda^{p_1}_{u_1}(w_1)}}
\def\loin{{\Lambda^{p,\infty}_u(w)}}

\def\loinu{{\Lambda^{p_1,\infty}_{u_1}(w_1)}}
\def\llo{{\Lambda^p(w)}}
\def\lloz{{\Lambda^{p_0}(w_0)}}

\def\lloin{{\Lambda^{p,\infty}(w)}}

\def\lloinu{{\Lambda^{p_1,\infty}(w_1)}}
\def\lpw{{L^p(w)}}
\def\lpwz{{L^{p_0}(w_0)}}
\def\lpwzdec{{L^{p_0}_{\hbox{\decp dec}}(w_0)}}

\def\lpwdec{{L^{p}_{\hbox{\decp dec}}(w)}}
\def\lpinc{{L^{p}_{\hbox{\decp inc}}}}
\def\lpdec{{L^{p}_{\hbox{\decp dec}}}}
\def\lpwu{{L^{p_1}(w_1)}}
\def\lpwinu{{L^{p_1,\infty}(w_1)}}

\def\BBinnux#1{B_{#1,\infty}}

\def\fiq{\phi_Q}
\def\Fiu{\Phi_u}
\def\limto#1{\smash{\mathop{\longrightarrow}\limits_{#1}}}
\def\Fd{\M_{\hbox{\decp dec}}(\R^+)}

\def\d{\downarrow}
\def\u{\uparrow}
\def\dip{\hbox{\decp d}}
\def\di{\hbox{\rm d}}
\def\Rm{\R^+}
\def\Rn{\R^n}
\def\lesssim{\mathrel{\mathop<\limits_{\sim}}}

\makeindex

\begin{document}

\title{Recent Developments in the Theory of Lorentz Spaces and
Weighted Inequalities}
\author{Mar\'{\i}a J. Carro \and Jos\'e A. Raposo \and Javier Soria}
\maketitle
\newpage
\pagenumbering{roman}\pagestyle{empty}

\centerline{\titl Recent Developments in the Theory of Lorentz}

\medskip

\centerline{\titl Spaces and Weighted Inequalities}

\vskip 3cm

\centerline{\nom Mar\'{\i}a J. Carro \hfill Jos\'e A. Raposo \hfill Javier Soria}

\vskip 2cm

\centerline{ADDRESS}

\medskip

\centerline{Dept. Appl. Math. and Analysis}
\centerline{Gran Via 585}
\centerline{University of Barcelona}
\centerline{E-08071 Barcelona, SPAIN}

\vskip 1cm

\centerline{E-MAIL: {\tt carro@mat.ub.es\hskip 1cm soria@mat.ub.es}}

\vfill


\medskip\noindent
This work has been partly supported by DGICYT PB97--0986 and CIRIT 1999SGR00061.
\newpage
\
\newpage
\ \vfill
\rightline{\dedi Per a la Beatriu}
\vfill
\newpage\pagestyle{empty}
\
\newpage
\pagestyle{plain}
\tableofcontents
\newpage\pagestyle{empty}
\
\newpage
\pagestyle{plain}
\addcontentsline{toc}{chapter}{Foreword}
\ \vskip 2cm
\noindent{\capit Foreword} 
\vskip2cm
\noindent 
On September 4, 1998, Jos\'e Antonio Raposo, the second named author, died when he was only 39. He had
just received his Ph.D. degree a few months earlier, under our supervision. He was really happy for all the
new projects he had for the future, and so were we, since he was an extraordinary mathematician, and a very
valuable friend. 

This work is an updated version of his thesis, written as a self-contained text, with most of
the motivations, examples and applications available in the literature. 

We want to thank all the people who have encouraged us to write this book, and specially Jos\'e Antonio's
family. We also thank  Joan Cerd\`a who has read the whole manuscript and has given us many good advises,
improving the final version of these notes.

\medskip
\rightline{Mar\'{\i}a J. Carro and Javier Soria}
\rightline{Barcelona, September 2000}

\newpage\pagestyle{empty}
\
\newpage

\pagenumbering{arabic}
\addcontentsline{toc}{chapter}{Introduction}
\ \vskip 2cm
\noindent{\capit Introduction} 
\vskip2cm
\pagestyle{plain}

\noindent 
The main objective of this work is to bring together two well known and, a priori, unrelated theories dealing
with weighted inequalities for the Hardy-Littlewood maximal operator
$$
Mf(x)\index{$M$}=\sup_{x\in Q}{1\over|Q|}\int_Q |f(y)|\,dy,\quad x\in\Rn.
$$
For this, we consider the boundedness of $M$ in the weighted Lorentz space $\lo$, introduced by G.G. Lorentz
(\cite{Lor}) in 1951, 
\bequ\label{intro}
M:\lo\longrightarrow\lo.
\eequ
Two examples are historically relevant as a motivation:
If $w=1$,  (\ref{intro}) corresponds to the study of the boundedness
$$
M:L^p(u)\longrightarrow L^p(u),
$$
which was characterized by B. Muckenhoupt (\cite{Muck1}) in 1972, and the solution is given by the so called
$A_p$ weights. The second case is when we take in (\ref{intro}) $u=1$. This is a more recent theory, and
was completely solved by  M.A. Ari\~no and  B. Muckenhoupt (see \cite{AM2}) in 1991. It turns out that the
boundedness
$$
M:\llo\longrightarrow\llo,
$$
can be seen to be equivalent to the boundedness of the Hardy operator $A$ restricted to decreasing functions
of
$L^p(w)$, since the nonincreasing rearrangement of $Mf$ is pointwise equivalent to $Af^*$. The class of
weights satisfying this boundedness is known as  $B_p$. Another results related to this problem can be
found in \cite{Saw2}, \cite{Ne1}, \cite{CS1}, \cite{CS2}, \cite{Sp1}, \cite{Sp2}, and \cite{HM}.

Even though the $A_p$ and $B_p$ classes enjoy some similar features, they come from very different
theories, and so are the techniques used on each case: Calder\'on--Zygmund decompositions and covering
lemmas for $A_p$,  rearrangement invariant properties and positive integral operators for $B_p$.

\medskip

It is our aim to give a unified version of these two theories. Contrary to what one could expect, the
solution is not given in terms of the limiting cases above considered (i.e., $u=1$ and $w=1$), but in a
rather more complicated condition, which reflects the difficulty  of estimating the distribution function of
the Hardy-Littlewood maximal operator with respect to general measures (some previous results on this
direction can be found in \cite{CHK} and \cite{HK}, for the particular case of the $L^{p,q}(u)$ Lorentz
spaces with a power weight $w$).

\pagestyle{myheadings}
\markboth{INTRODUCTION}{M.J.\ Carro, J.A.\ Raposo, and J.\ Soria}

\medskip

In order to carry out this program as a self-contained monograph, we study in the first two chapters the
main results needed for the rest of the book. Hence, in Chapter~\ref{ch:characteristic} we consider the  
boundedness of integral operators on monotone functions. The most important  result is Theorem~\ref{te:
I.2.13}, where we prove that under very general conditions, the boundedness is determined by the action of
the operator on some characteristic functions (see also \cite{BPS} and \cite{CS2}). We then consider the main
example in this context, namely the Hardy operator, and introduce the $B_p$ class.

\medskip

In Chapter~\ref{ch:lor} we make an exhaustive study of the functional properties of the Lorentz spaces over
general measure spaces, and arbitrary weights. In particular we consider the discrete case (sequence Lorentz
spaces) where we answer several open questions about normability and  duality (see Theorem~\ref{te:
II.4.18}). Some previous results were already proved (in the case $w$ decreasing) in \cite{Al} and
\cite{AEP}.

\medskip

Finally, in  Chapter~\ref{ch:maximal} we consider the solution to (\ref{intro}), which is proved in
Theorem~\ref{te:principal} (some sufficient conditions had been already obtained in \cite{CS3} and
\cite{Ne2}). We see that the characterization that we obtain is based upon the so called
$p\,-\,\epsilon$ condition, which is natural since both $A_p$ and $B_p$ enjoy this property. However is worth
noticing that  there exist non-doubling weights $u$ for which $M$ satisfies (\ref{intro}) for suitable
weights $w$ (see Theorem~\ref{te: pnd}). We also study the weak-type and the restricted weak-type
version of (\ref{intro}).  

\medskip

As far as possible, we have always tried to give precise bibliographic information about the results which
were previously known, and also about the techniques used in the proofs, although due to the big number of
people   working in these theories, this a difficult task. We have also followed the standard notation found
in the main reference books (e.g., \cite{BS}, \cite{BL}, \cite{GR}, \cite{Ru2}, \cite{St}, \cite{SW}). 

\chapter{Boundedness of operators on characteristic functions and the Hardy operator} 
\label{ch:characteristic}
\pagestyle{myheadings}
\markboth{CHAPTER 1. BOUNDEDNESS OF OPERATORS}{M.J.\ Carro, J.A.\ Raposo, and J.\ Soria}

\section{Introduction}
\label{se I: Introduction} 

This chapter includes general results on boundedness of operators on  
$L^p\/$ such as the 
Hardy operator\index{Hardy operator}
\bequ
\label{hardy}
Af(t)\index{$A$}={1\over t}\int_0^t f(s)\,ds,\quad t>0,
\eequ 
and the   Hardy-Littlewood maximal operator\index{Hardy-Littlewood maximal operator}
\bequ
\label{hardylit}
Mf(x)\index{$M$}=\sup_{x\in Q}{1\over|Q|}\int_Q |f(y)|\,dy,\quad x\in\Rn.
\eequ

The main result is  Theorem~\ref{te: I.2.13} and its Corollaries~\ref{cor: I.2.14}  and 
\ref{cor: I.2.19} where it is proved that a great variety of operators satisfy that their boundedness
is characterized by the restriction to characteristic functions. 

In the third section, we study the  Hardy operator, recalling some known results that we shall use
in the following chapters and including a detailed study of the discrete  Hardy operator\index{Hardy
operator!discrete},
\bequ
\label{dischardy}
A_d f(n)\index{$A_d$}={1\over n+1}\sum_{k=0}^n f(k),\qquad n\in\N.
\eequ

The following standard notations will be used very often: Letters such as $X$ or $Y$ are always
$\sigma$-finite measure spaces and ${\cal M}(X)$\index{${\cal M}(X)$} (${\cal M}^+(X)$)\index{${\cal
M}^+(X)$} denotes the space of measurable (resp. measurable and nonnegative) functions on
$X$.  The distribution function\index{distribution function} of $f$ is 
$$
\lambda_f(t)\index{$\lambda_f$}=\mu(\{x:|f(x)|>t\}),
$$
the nonincreasing rearrangement\index{nonincreasing rearrangement} is
$$
f^*(t)\index{$f^*$}=\inf\{s>0:\lambda_f(s)\le t\},
$$
and
$$
f^{**}(t)\index{$f^{**}$}={1\over t}\int_0^t f^*(s)\,ds.
$$
If $f$ is a positive nonincreasing (nondecreasing) function we will write
$f\downarrow$ (resp. $f\uparrow$). The space $\lpdec$\index{$\lpdec$} ($\lpinc$)\index{$\lpinc$} is the cone
of all nonincreasing (resp. nondecreasing) functions in $L^p$. Two positive quantities $A$ and $B$, are said
to be equivalent ($A\approx B$) if there exists a constant $C>1$ (independent of the essential parameters
defining
$A$ and
$B$) such that $C^{-1}A\le B\le CA$. If only $ B\le CA$, we write $B\lesssim A$. The
undetermined cases $0\cdot\infty,\,{\infty\over\infty},\,{0\over 0}\/$, will always be taken equal to 
$0\/$.

\section{Boundedness on characteristic functions}
\label{se I: characteristic functions}

In this section we consider the  problem of characterizing the boundedness of 
\bequ
\label{I.2.1}
T\,:\,\big(L^{p_0}(X)\cap L\,,\,\|\cdot\|_{p_0}\big)\longrightarrow L^{p_1}(Y),
\eequ 
where
$L\/$ is a subclass of ${\cal M}(X)\/$ and  $T\/$ is an operator (usually linear or sublinear). That
is, we are interested in studying the inequality
\bequ
\label{I.2.2}
\|Tf\|_{L^{p_1}(Y)}\le C\|f\|_{L^{p_0}(X)},\quad f\in L.
\eequ  
In some cases inequality  (\ref{I.2.2}) holds for every  $f\in L\/$ if and only if it holds on
characteristic functions
$f=\chi_A\in L\/$. This happens (see \cite{CS2} and \cite{Sp1}) in the case 
$L=\lpwzdec$, $X=(\Rm,w_0(t)\,dt)$ and 
$Y=(\Rm,w_1(t)\,dt)\/$ in the range of indices  $0<p_0\le1,\;p_0\le p_1<\infty\/$ and operators of
the type 
\bequ
\label{I.2.3}
Tf(r)=\int_0^\infty k(r,t)f(t)\,dt\,,\quad r>0,
\eequ  
for positive kernels $k$. Using this, one can easily obtain a useful characterization of   (\ref{I.2.1}). For
example, in the previous case the condition on the weights 
$w_0,w_1\/$ (nonnegative locally integrable functions
in $\R^+$)  for  which $T:\lpwzdec\to\lpwu\/$ (defined by (\ref{I.2.3})) is:
$$
\bigg(\int_0^\infty\bigg(\int_0^r k(s,t)\,dt\bigg)^{p_1}w_1(s)\,ds\bigg)^{1/p_1}\le
C\,\bigg(\int_0^r w_0(s)\,ds\bigg)^{1/p_0}\,,
\quad r>0. 
$$ 
The range $0<p_0\le1,\;p_0\le p_1<\infty\/$ is fundamental and, in fact, the previous principle
also works in other situations, as in the case of  increasing functions (see \cite{Sp1}) 
$$
T:L^{p_0}_{\hbox{\decp inc}}(w_0)\to\lpwu
$$ 
where $T\/$ is as in  (\ref{I.2.3}).

As we shall see (Theorem~\ref{te: I.2.13}) all this is a consequence of a general principle that 
can be applied to a more general class of operators than those of integral type  and for a
bigger class of functions that includes the monotone functions.  

We first need some definitions.

\bdefi
\label{def: I.2.4} 
We say that $\emptyset\neq L\subset \M(X)\/$ is a  regular class\index{regular class} in 
$X\/$ if, for every  $f\in L\/$,   
\begin{enumerate}
\item[(i)] $|\alpha f|\in L\/$, for every  $\alpha\in\R\/$,

\item[(ii)] $\chi_{\{|f|>t\}}\in L\/$, for every  $t>0\/$, and 

\item[(iii)] there exists a sequence of simple functions $(f_n)_n\subset L\/$ such that $0\le f_n(x)\le
f_{n+1}(x)\to |f(x)|\hbox{\ a.e.\ }x\in  X\/$.
\end{enumerate}
\edefi

\beje
\label{eje: I.2.5}

(i) If $X\/$ is an arbitrary measure space,  every functional lattice\index{functional lattice} in $X\/$
(i.e., a vector space
$L\subset\M(X)\/$ with 
$g\in L\/$ if $|g|\le|f|,\;f\in L\/$) with the  Fatou property\index{Fatou property} (see \cite{BS}), is a
regular class. In particular the  Lebesgue space
$L^p(X)\/$ and the  Lorentz space $L^{p,q}(X),\;0<p,q\le\infty\/$ (see (\ref{normlor})) are regular classes.

(ii) If $C\/$ is a class formed by measurable sets in $X\/$  containing the empty set,
$$
L_C=\big\{f\in\M(X)\,:\,\{|f|>t\}\in C,\;\forall t\ge0\big\}
$$ 
is a  regular class in $X\/$. To see this we observe that the conditions (i) and (ii) of 
Definition~\ref{def: I.2.4} are immediate to be checked while to prove (iii) we observe that if 
$0\le f\in L_C\/$, the simple functions
$$
f_n=\sum_{k=0}^{n2^{n-1}}k2^{-n}\chi_{\{k2^{-n}<f\le(k+1)2^{-n}\}}+n\chi_{\{f>n\}},\quad
n=1,2,\dots,
$$ 
form an increasing sequence that converges pointwise to    $f\/$, whose  level sets are also level sets
of   $f\/$  and hence they belong to  $C\/$.

(iii) If $L\/$ is a  regular class in  $X\/$,
$$
L^\ast=\big\{f^\ast\,:\,f\in L\big\}
$$ 
is a  regular class in $\Rm\/$.

(iv) In $\Rm\/$ the positive decreasing functions are a regular class and the same holds for the
increasing functions.

(v) In $\Rn\/$   radial functions, positive and decreasing  radial functions, or  positive and
increasing  radial functions are also regular classes.
\eeje

\bdefi
\label{def: I.2.6}
Let  $Y\/$ be a   measure space, $L\/$ a regular class and $T:L\to \M(Y)\/$
an operator.

(i) $T\/$ is sublinear if $|T(\alpha_1 f_1+\dots+\alpha_n f_n)(y)|\le|\alpha_1
Tf_1(y)|+\dots+|\alpha_n Tf_n(y)|\/$ a.e. 
$y\in Y\/$, for every $\alpha_1,\dots,\alpha_n\in\R$, and every $f_1,\dots,\,f_n\in L\/$ such that $\alpha_1
f_1+\dots+\alpha_n f_n\in L\/$.

(ii) $T\/$ is monotone if $|Tf(y)|\le|Tg(y)|\hbox{\ a.e.\ }y\in Y\/$, if  $|f(x)|\le|g(x)|\hbox{\
a.e.\ }x\in X,\;f,g\in L.\/$

(iii) We say that $T\/$ is order continuous\index{order continuous} if it is  monotone and if for every
sequence
$(f_n)_n\subset L\/$  with
$0\le f_n(x)\le f_{n+1}(x)\to f(x)\in L\hbox{\ a.e.\ }x\in X\/$,  we have that $\lim_n
|Tf_{n}(y)|=|Tf(y)|\hbox{\ a.e.\ }y\in Y\/$.
\edefi

\bobs
\label{obs: I.2.7}
It is immediate that a sublinear and  monotone  operator satisfies the following properties:

(i) $T(0)(y)=0\hbox{\ a.e.\ }y\in Y\/$,

(ii) $\big|T|f|(y)\big|=|Tf(y)|\hbox{\ a.e.\ }y\in Y,\;f\in L\/$.
\eobs

\bobs
\label{obsI.2.77}
(i) Every maximal operator of the form
$$
T^\ast f(x)=\sup_{T\in B}|Tf(x)|\,,\quad x\in X,\quad f\in L,
$$ 
is order continuous, where, for every  $x\in X\/$,
$B\/$ is a set of order continuous operators $T$.

(ii) In particular every  integral operator $Tf(r)=\displaystyle\int_0^{\infty} k(r,t)f(t)\,dt,\,r>0,\;f\in
L\subset
\M^+(\Rm)\/$ (with $k:\Rm\times\Rm\to[0,\infty)\/$), is order continuous. For example the Hardy operator
(see (\ref{hardy})) and its conjugate 
$$
Qf(r)=\int_r^\infty f(t)\,{dt\over t}.
$$

(iii)  Also the identity operator 
$f\mapsto f\/$ is obviously order continuous.
\eobs

The following  property  of the sublinear and monotone operators will be  fundamental.

\bteo
\label{te: I.2.8}
 Let  $L\subset \M(X)\/$ be a regular class and  $T:L\to \M(Y)\/$ a
sublinear and monotone operator. Then, for every simple function     $f\in L\/$,
$$
|Tf(y)|\le\int_0^\infty\big|T\chi_{\{|f|>t\}}(y)\big|\,dt\,,\quad\hbox{a.e.\ } y\in Y. 
$$
\eteo

\bdem      
By Remark~\ref{obs: I.2.7} (ii) we can assume $f\ge 0\/$. Then,
$$
f=\sum_{n=1}^N a_n\,\chi_{B_n},
$$ 
with $a_1,a_2,\dots,a_N>0\/$ and  $(B_j)_j\/$ is an increasing sequence of measurable sets in
$X\/$:
$B_1\subset B_2\subset\dots\subset B_N\/$. Set 
$A_n=\sum_{j=n}^N a_j\,,\;n=1,\dots,N,$\break $A_{N+1}=0\/$ and note that
$\{f>t\}=\emptyset\/$ if $t\ge A_1\/$ and $\{f>t\}=B_n\/$ for $A_{n+1}\le t<A_n\,,\;n=1,\dots,N\/$.
In particular   $\chi_{B_n}\/$ is in  $L\/$ (by Definition~\ref{def: I.2.4} (ii)) and since  $T\/$
is sublinear, it follows that 
\bequ
\label{I.2.9}
|Tf(y)|\le\sum_{n=1}^N a_n|T\chi_{B_n}(y)|\,,\quad\hbox{a.e.\ }y\in Y. 
\eequ
On the other hand
$T\chi_\emptyset(y)=0\/$ a.e. $y\in Y\/$ (by Remark~\ref{obs: I.2.7} (i)) and hence
$$
\int_0^\infty\big|T\chi_{\{f>t\}}(y)\big|\,dt=\sum_{n=1}^N\int_{A_{n+1}}^{A_n}|T
\chi_{B_n}(y)|\,dt=\sum_{n=1}^N
a_n|T\chi_{B_n}(y)|\,,
$$ 
since $A_n-A_{n+1}=a_n\/$. The last expression coincides with the right hand side of  (\ref{I.2.9}) and the 
theorem is proved.$\qquad\qed$
\edem

\bobs
\label{obs: I.2.10}
If the  operator   $T\/$ is linear and  positive ($Tf\ge0\/$ if $f\ge0\/$), the inequality
in the previous theorem is, in fact, an equality when
$f\ge0\/$.  
\eobs

The following equality will be very much useful later on and it has been proved by several authors.
See, for example,  \cite{HM}.

\blem
\label{lem: I.2.11} 
 If $0<p\le 1\/$,
$$
\sup_{f\downarrow}\,{\bigg(\displaystyle\int_0^\infty f(t)\,dt\bigg)^p\over\displaystyle\int_0^\infty
f^p(t)t^{p-1}\,dt}=p\,. 
$$
\elem

We introduce now some new spaces which play an important role in this theory, and make some important
comments about its definition.

\bdefi\label{def: Lor}
If $0<p,q<\infty$, we define the Lorentz space\index{Lorentz space}
\bequ\label{normlor}
L^{p,q}(X)\index{$L^{p,q}(X)$}=\bigg\{f:\Vert
f\Vert_{L^{p,q}(X)}=\bigg(\int_0^{\infty}(t^{1/p}f^*(t)\big)^q\,{dt\over t}\bigg)^{1/q}<\infty\bigg\}.
\eequ
If $q=\infty$ the space $L^{p,\infty}(X)$\index{$L^{p,\infty}(X)$} is defined with the usual modification.
When the space
$X$ is clearly understood on the context, we will simply write $\Vert f\Vert_{p,q}.$
\edefi

\bobs
\label{obs: I.2.12} 

(i) If $1\le q\le p\/$, the functional $\|\cdot\|_{p,q}\/$ is a norm. In general, it is only a quasi-norm
(see
\cite{SW}).

(ii) In the case $1<p<q\le\infty\/$ we can define
$$
\|f\|_{(p,q)}=\bigg(\int_0^\infty \big(t^{1/p}f^{\ast\ast}(t)\big)^q\,{dt\over t}\bigg)^{1/q}
$$
(with the usual modification  if $q=\infty\/$) which is a norm (see \cite{BS}) and it is equivalent
to the original quasi-norm:
\bequ\label{eqlo}
\|f\|_{p,q}\le\|f\|_{(p,q)}\le{p\over p-1}\,\|f\|_{p,q},\quad f\in\M(X).
\eequ 

(iii) It is easy to show that $\Vert
f\Vert_{p,q}=\displaystyle\bigg(\int_0^{\infty}(t\lambda_f^{1/p}(t))^q\,{dt\over t}\bigg)^{1/q}$ (see
Proposition \ref{pro: II.2.5}).
\eobs

Let us state the main result of this section.  

\bteo
\label{te: I.2.13}
Let $L\subset \M(X)\/$ be a regular class, $T:L\to \M(Y)\/$ an order continuous  sublinear
operator and  $0<q_0\le1,\,0<p_0<\infty\/$. Then
\begin{enumerate}
\item[(a)] If $q_0\le q_1\le p_1<\infty$,
$$
\sup_{f\in L}{\|Tf\|_{L^{p_1,q_1}(Y)}\over \|f\|_{L^{p_0,q_0}(X)}}=\sup_{\chi_B\in
L}{\|T\chi_B\|_{L^{p_1,q_1}(Y)}\over \|\chi_B\|_{L^{p_0,q_0}(X)}}. 
$$

\item[(b)] If $q_0<p_1<q_1\le\infty\/$,  
$$
\sup_{f\in L}{\|Tf\|_{L^{p_1,q_1}(Y)}\over \|f\|_{L^{p_0,q_0}(X)}}\le\bigg({p_1\over
p_1-q_0}\bigg)^{1/q_0}\,\sup_{\chi_B\in L}{\|T\chi_B\|_{L^{p_1,q_1}(Y)}\over
\|\chi_B\|_{L^{p_0,q_0}(X)}}. 
$$
\end{enumerate}
\eteo

\bdem       
Let
$$
C=\sup_{\chi_B\in L}{\|T\chi_B\|_{L^{p_1,q_1}(Y)}\over \|\chi_B\|_{L^{p_0,q_0}(X)}}. 
$$
We have to prove that  $\|Tf\|_{L^{p_1,q_1}(Y)}\le K\,C\,\|f\|_{L^{p_0,q_0}(X)}\/$ for every  $f\in
L\/$ with $K=1\/$ (in the case (a)) or $K=\big(p_1/(p_1-q_0)\big)^{1/q_0}\/$ (case (b)). By
Remark~\ref{obs: I.2.7} (ii) we can assume $f\ge0\/$ and since there exists an increasing sequence
of positive simple functions  
$(f_n)_n\subset L\/$  converging to  $f\/$ a.e., with
$$
|Tf_n(y)|\le|Tf_{n+1}(y)|\rightarrow|Tf(y)|\qquad\hbox{a.e.\ }y,
$$ 
(Definitions~\ref{def: I.2.4} and \ref{def: I.2.6}), by the monotone convergence theorem, it is
enough to prove the previous inequality for 
$(f_n)_n\/$. That is, we can assume that $0\le f\in L\/$ is a simple function. 

Let $p=p_1/q_0,\,q=q_1/q_0\/$. By Theorem~\ref{te: I.2.8} and Lemma~\ref{lem: I.2.11},
$$
\|Tf\|_{p_1,q_1}^{q_0}=\big\||Tf|^{q_0}\big\|_{p,q}\le\bigg\|\int_0^\infty q_0
t^{q_0-1}\big|T\chi_{\{f>t\}}(\cdot)\big|^{q_0}\,dt\,\bigg\|_{p,q}. 
$$ 
In the case (a),
$\|\cdot\|_{p,q}\/$ is a norm and we obtain that  
$$
\|Tf\|_{p_1,q_1}^{q_0}\le\int_0^\infty q_0
t^{q_0-1}\big\|\big|T\chi_{\{f>t\}}\big|^{q_0}\big\|_{p,q}\,dt. 
$$ 
In the case  (b),
$\|\cdot\|_{(p,q)}\/$ is a norm satisfying (\ref{eqlo})
and it follows that
\begin{eqnarray*}
\|Tf\|_{p_1,q_1}^{q_0}&\le&{p\over p-1}\int_0^\infty q_0
t^{q_0-1}\big\|\big|T\chi_{\{f>t\}}\big|^{q_0}\big\|_{p,q}\,dt\\ &=&{p_1\over p_1-q_0}\int_0^\infty
q_0 t^{q_0-1}\big\|\big|T\chi_{\{f>t\}}\big|^{q_0}\big\|_{p,q}\,dt.
\end{eqnarray*}
That is, in any case
\begin{eqnarray*}
\|Tf\|_{p_1,q_1}^{q_0}&\le& K^{q_0}\int_0^\infty q_0
t^{q_0-1}\big\|\big|T\chi_{\{f>t\}}\big|^{q_0}\big\|_{p,q}\,dt\\ &=&K^{q_0}\int_0^\infty q_0
t^{q_0-1}\big\|T\chi_{\{f>t\}}\big\|_{p_1,q_1}^{q_0}\,dt\\ &\le& (KC)^{q_0}\int_0^\infty q_0
t^{q_0-1}\big\|\chi_{\{f>t\}}\big\|_{p_0,q_0}^{q_0}\,dt\\ &=&(KC)^{q_0}\int_0^\infty p_0
t^{q_0-1}\big(\lambda_f(t)\big)^{q_0/p_0}\,dt\\ &=&(KC)^{q_0}\|f\|_{p_0,q_0}^{q_0}.\qquad\qed
\end{eqnarray*}
\edem

Applying the previous result to the strong (i.e., diagonal) case $T:L^{p_0}\rightarrow L^{p_1}\/$ we obtain:

\bcor
\label{cor: I.2.14}
 Let $L\subset \M(X)\/$ be a regular class and  $T:L\to \M(Y)\/$ an order continuous
sublinear operator. If
$0<p_0\le 1,\,p_0\le p_1<\infty\/$ we have 
$$
\sup_{f\in L}{\|Tf\|_{L^{p_1}(Y)}\over \|f\|_{L^{p_0}(X)}}=\sup_{\chi_B\in
L}{\|T\chi_B\|_{L^{p_1}(Y)}\over \|\chi_B\|_{L^{p_0}(X)}}. 
$$
\ecor

\bobs
\label{obs: I.2.15}
 The range of exponents $p_0,p_1\/$ in Corollary~\ref{cor: I.2.14} is optimal. To see this,
observe:
\begin{enumerate}
\item[(i)] The result is not true if $p_0>1\/$. A counterexample is the 
 Hardy operator $A\/$. Given $1<p_0\le p_1\/$, a necessary and sufficient condition  to have the
boundedness 
$A:L^{p_0}_{\hbox{\decp dec}}(w_0)\to L^{p_1}(w_1)\/$ is (see \cite{Saw2}),
\begin{eqnarray}
\label{I.2.16}
 W_1^{1/p_1}(r)&\le& C\,W_0^{1/p_0}(r),\nonumber\\ 
\bigg(\int_r^\infty{w_1(t)\over
t^{p_1}}\,dt\bigg)^{1/p_1}\bigg(\int_0^r\bigg({W_0(t)\over
t}\bigg)^{-p_0^\prime}w_0(t)\,dt\bigg)^{1/p_0^\prime}&\le& C.
\end{eqnarray}
It is easy to see that the condition  
$$\big\|A\chi_{(0,r)}\big\|_{L^{p_1}(w_1)}\le
C\,\big\|\chi_{(0,r)}\big\|_{L^{p_0}(w_0)},\;r>0\/$$ 
obtained from Corollary~\ref{cor: I.2.14}  is equivalent to
the inequalities
\begin{eqnarray*}
 W_1^{1/p_1}(r)&\le& C_1\,W_0^{1/p_0}(r),\quad r>0,\\ 
\bigg(\int_r^\infty\bigg({r\over t}\bigg)^{p_1}w_1(t)\,dt\bigg)^{1/p_1}&\le&
C_1\,W_0^{1/p_0}(r),\quad r>0.
\end{eqnarray*}
Observe that these two last inequalities are true for the weights 
$w_0(t)=t^{p_0-1}\,,w_1(t)=t^{p_1-1}\chi_{(1,2)}(t)\/$  while these weights do not satisfy
(\ref{I.2.16}).

\item[(ii)] The statement does not hold if  $p_1<p_0\/$. To see this, it is enough to consider 
$T=\hbox{\rm Id}:\lpwzdec\to\lpwu\/$. Then,
\bequ
\label{I.2.17}
\sup_{f\downarrow}{\displaystyle\int_0^\infty
f^{p_1}(t)w_1(t)\,dt\over\bigg(\displaystyle\int_0^\infty
f^{p_0}(t)w_0(t)\,dt\bigg)^{p_1/p_0}}=\sup_{g\downarrow}{\displaystyle\int_0^\infty
g(t)w_1(t)\,dt\over\bigg(\displaystyle\int_0^\infty g^{p_0/p_1}(t)w_0(t)\,dt\bigg)^{p_1/p_0}}.
\eequ 
The last supremum is equivalent (c.f.
\cite{Saw2}), up to multiplicative constants depending only on  $p_0/p_1\/$, to 
\bequ
\label{I.2.18}
\bigg(\int_0^\infty\bigg({W_1(t)\over
W_0(t)}\bigg)^{p_1/(p_0-p_1)}w_1(t)\,dt\bigg)^{(p_0-p_1)/p_0}.
\eequ
On the other hand, if  Corollary~\ref{cor: I.2.14}  were true in this case, (\ref{I.2.17}) would be
equal to 
$$
\sup_{r>0}{\|\chi_{(0,r)}\|_{\lpwu}\over\|\chi_{(0,r)}\|_{\lpwz}}=\sup_{r>0}{W_1^{1/p_1}(r)\over
W_0^{1/p_0}(r)}\,.
 $$ 
This last supremum is finite for the weights 
$w_i(t)=t^{\alpha_i}\,,\;i=0,1\/$, if
$(1+\alpha_1)/p_1=(1+\alpha_0)/p_0\,,\;\alpha_0,\alpha_1>-1\/$, while   (\ref{I.2.18}) is always
infinite in this case whatever 
$\alpha_1,\alpha_2\/$ are.
\end{enumerate}
\eobs

Considering the weak-type case $T:L^{p_0}\rightarrow L^{p_1,\infty}\/$, we  have as a
consequence of Theorem~\ref{te: I.2.13}  the following statement.

\bcor
\label{cor: I.2.19}
 Let $L\subset \M(X)\/$ be a  regular class and let $T:L\to \M(Y)\/$ be an order continuous sublinear
operator.   If $0<p_0\le 1,\,p_0<p_1<\infty\/$,
$$
\sup_{f\in L}{\|Tf\|_{L^{p_1,\infty}(Y)}\over \|f\|_{L^{p_0}(X)}}\le\bigg({p_1\over
p_1-p_0}\bigg)^{1/p_0}\sup_{\chi_B\in L}{\|T\chi_B\|_{L^{p_1,\infty}(Y)}\over
\|\chi_B\|_{L^{p_0}(X)}}.
 $$
\ecor

\bobs
\label{obs: I.2.20} The previous corollary is not true if  $p_0>1\/$, as  can be seen in the case
\bequ
\label{I.2.21}
A\,:\,\lpwzdec\longrightarrow\lpwinu.
\eequ
As  was proved in \cite{CS2} (see also \cite{Ne1}), a necessary and sufficient condition to have 
(\ref{I.2.21}) is
\begin{eqnarray} 
\label{I.2.22}
\bigg(\int_0^r\bigg({W_0(t)\over
t}\bigg)^{-p_0^\prime}w_0(t)\,dt\bigg)^{1/p_0^\prime}W_1^{1/p_1}(r)&\le& Cr,\quad r>0,\\ 
W_1^{1/p_1}(r)&\le& CW_0^{1/p_0}(r),\quad r>0,\nonumber
\end{eqnarray} 
that   does not coincide with the condition
$$
{W_1^{1/p_1}(t)\over t}\le C\,{W_0^{1/p_0}(r)\over r}\,,\quad 0<t<r<\infty\,,
$$
which can be obtained applying Corollary~\ref{cor: I.2.19}.  For example, the two weights
in Remark~\ref{obs: I.2.15} (iii.1) satisfy this last condition but  not  
(\ref{I.2.22}).
\eobs

A consequence of  Theorem~\ref{te: I.2.13}
and the   Marcinkiewicz interpolation theorem is the following result  on interpolation of
operators of restricted weak-type, which is a generalization of the well known 
Stein-Weiss theorem (Theorem 3.15 in \cite{SW} or Theorem 5.5 in \cite{BS}).

\bteo
\label{te: I.2.21}
 Let $0<p_0,p_1<\infty,\;0<q_0,q_1\le\infty,\;p_0\neq p_1,\;q_0\neq q_1\/$ and let us assume that 
$T:\big(L^{p_0,1}+L^{p_1,1}\big)(X)\to\M(Y)\/$ is a sublinear order continuous  operator satisfying
\begin{eqnarray*}
\big\|T\chi_B\big\|_{q_0,\infty}&\le& C_0\|\chi_B\|_{p_0,1},\quad B\subset X,\\ 
\big\|T\chi_B\big\|_{q_1,\infty}&\le& C_1\|\chi_B\|_{p_1,1},\quad B\subset X.
\end{eqnarray*}
 Then
$$
T\,:\,L^{p,r}(X)\longrightarrow L^{q,r}(Y),\qquad 0<r\le\infty,
$$ 
if
$$
{1\over p}={1-\theta\over p_0}+{\theta\over p_1},\qquad {1\over q}={1-\theta\over
q_0}+{\theta\over q_1},\qquad0<\theta<1. 
$$
\eteo

\bdem      
Since $L^{s,t_0}\subset L^{s,t_1}\/$ if  $t_0<t_1\/$, there exist indices
$r_0,r_1\in(0,1)\/$ with $r_i<q_i,\;i=0,1\/$ and such that 
$$
\big\|T\chi_B\big\|_{q_i,\infty}\le C_i\|\chi_B\|_{p_i,r_i},\quad B\subset X,\qquad i=0,1. 
$$
Theorem~\ref{te: I.2.13}
tells us that the previous inequality holds for every function   $f\in
L^{p_i,r_i}\/$  and the result follows from the general Marcinkiewicz interpolation theorem 
 (Theorem 5.3.2 in \cite{BL}).$\qquad\qed$
\edem

\bobs
\label{obs: I.2.22}

(i) The norm of a characteristic  function    in $L^{p,q}\/$ does not depend (up to constants)
on $q\/$. Therefore, the previous  theorem remains true if we substitute the original spaces 
$L^{p_i,1}\/$ by  $L^{p_i,r_i}\/$ with $0<r_i\le\infty,\;i=0,1\/$.

(ii) In the classical result of  Stein-Weiss mentioned above, the previous result is proved for a
more restricted set of indices:
$$
1\le p_0,p_1<\infty,\qquad 1\le q_0,q_1\le\infty. 
$$

(iii) Corollary~\ref{cor: I.2.14} can be also generalized to consider the case of two operators. The
following result gives such extension (the proof is similar to the one given in Corollary~\ref{cor: I.2.14},
and it is also based upon some ideas found in \cite{BPS}).
\eobs

\bteo\label{twoop}
Let $L\subset{\cal M}(X)$ be a regular class, and let $T_i:L\longrightarrow {\cal M}(Y_i)$, $i=0,1$ be two
order continuous operators. Assume that $T_0$ is positive and linear and $T_1$ is sublinear. Then, for each
of the following cases:
\begin{enumerate}
\item[(a)] $0<p_0\le1\le p_1<\infty$,
\item[(b)] $T_1=\hbox{\rm Id}$, $1\le p_1<\infty$, $0<p_0\le p_1$,
\item[(c)] $T_0=\hbox{\rm Id}$, $0< p_0\le1$, $p_0\le p_1<\infty$,
\item[(d)] $T_0=T_1$,  $0<p_0\le p_1<\infty$,
\end{enumerate}
we have that
$$
\sup_{f\in L}{\Vert T_1 f\Vert_{L^{p_1}(Y_1)}\over \Vert T_0 f\Vert_{L^{p_0}(Y_0)}}=
\sup_{\chi_B\in L}{\Vert T_1 \chi_B\Vert_{L^{p_1}(Y_1)}\over \Vert T_0 \chi_B\Vert_{L^{p_0}(Y_0)}}.
$$

\eteo

\section{The Hardy operator and the classes $B_p$}
\label{se I: Hardy operator}

The Hardy operator $A$ 
defined in (\ref{hardy}) 
will play a fundamental role in the following chapters. In particular we shall be interested in the
boundedness
\bequ
\label{I.6.1}
A\,:\,\lpwzdec\longrightarrow\lpwu 
\eequ 
and, also,
\bequ
\label{I.6.2} 
A\,:\,\lpwzdec\longrightarrow\lpwinu.
\eequ
The  diagonal case $A:L^p_{\hbox{\decp
dec}}(w)\to\lpw\,,p>1\/$ was solved by   Ari\~no and Muckenhoupt in \cite{AM1}. 
The condition on the weight  $w\/$ is known as $B_p\/$. This motivates the following
definition.

\bdefi
\label{def: I.6.3}
We write  $w\in B_p\/$\index{$B_p$} if
$$
A\,:\,L^p_{\hbox{\decp dec}}(w)\longrightarrow\lpw,
$$ 
and  $w\in B_{p,\infty}\/$\index{$B_{p,\infty}$} if
$$
A\,:\,L^p_{\hbox{\decp dec}}(w)\longrightarrow L^{p,\infty}(w). 
$$ 
Analogously we write
$(w_0,w_1)\in B_{p_0,p_1}\/$\index{$B_{p_0,p_1}$} if the boundedness  (\ref{I.6.1}) holds and we say that 
$(w_0,w_1)\in B_{p_0,p_1,\infty}\/$\index{$B_{p_0,p_1,\infty}$} if  (\ref{I.6.2}) holds.
\edefi

\bobs
\label{obs: I.6.4} 
Since $\lpw\subset L^{p,\infty}(w)\/$ we have that
$$
B_{p_0,p_1}\subset B_{p_0,p_1,\infty},\quad0<p_0,p_1<\infty. 
$$ 
In particular $B_p\subset
B_{p,\infty},\,0<p<\infty\/$.
\eobs

The characterization of the weights satisfying (\ref{I.6.2}) is obtained applying directly  Theorem 
3.3 in \cite{CS2}:

\bteo 
\label{te: I.6.5} 
Let $0<p_1<\infty\/$. Then,
\begin{enumerate}
\item[(a)]\ \ If $p_0>1\/$ the following conditions are equivalent:
\begin{enumerate}
\item[(i)] $(w_0,w_1)\in B_{p_0,p_1,\infty},$
\item[(ii)] $\bigg(\displaystyle\int_0^r\bigg({W_0(t)\over
t}\bigg)^{1-p_0^\prime}\,dt\bigg)^{1/p_0^\prime}W_1^{1/p_1}(r)\le Cr,\quad r>0,$
\item[(iii)] 
\begin{eqnarray*}
\bigg(\int_0^r\bigg({W_0(t)\over
t}\bigg)^{-p_0^\prime}w_0(t)\,dt\bigg)^{1/p_0^\prime}W_1^{1/p_1}(r)&\le& Cr,\\ 
W_1^{1/p_1}(r)&\le&
CW_0^{1/p_0}(r),\quad r>0.\\
\end{eqnarray*}
\end{enumerate}

\item[(b)]\ \ If $p_0\le1\/$ the following conditions are equivalent:
\begin{enumerate}
\item[(i)] $(w_0,w_1)\in B_{p_0,p_1,\infty},$
\item[(ii)] ${W_1^{1/p_1}(r)\over r}\le
C{W_0^{1/p_0}(t)\over t},\quad0<t<r.$
\end{enumerate}
\end{enumerate}
\eteo

The strong boundedness $A:\lpwzdec\to\lpwu\/$ is not so easy and, in fact, the case
 $0<p_1<p_0<1\/$ is still open. The result that follows characterizes the classes $B_p\/$. The
proof of 
$(i) \Leftrightarrow (ii)\/$ can be found in  \cite{AM1}  (case $p\ge1\/$) and it is a consequence of 
Corollary~\ref{cor: I.2.14}  in the case $p<1\/$. The equivalences with  (iii) and (iv) are proved
in  \cite{Sor}.

\bteo{(Ari\~no-Muckenhoupt, J. Soria)}
\label{te: I.6.6} 
For $0<p<\infty\/$ the following statements are  equivalent:
\begin{enumerate}
\item[(i)] $ w\in B_p.$
\item[(ii)] $\displaystyle\int_r^\infty\bigg({r\over t}\bigg)^p w(t)\,dt\le
C\int_0^r w,\quad r>0.$
\item[(iii)] $\displaystyle\int_0^r{t^{p-1}\over W(t)}\,dt\le C{r^p\over W(r)},\quad
r>0.$
\item[(iv)] $\displaystyle\int_0^r {1\over W^{1/p}(t)}\,dt\le C{r\over W^{1/p}(r)},\quad r>0.$
\end{enumerate}
\eteo

The classes $B_{p_0,p_1}\/$ with $p_0\le1,\,p_0\le p_1\/$ can be characterized  using
Corollary~\ref{cor: I.2.14} (see also \cite{CS2}). The case $p_1,p_0>1\/$ was solved by  Sawyer
(\cite{Saw2}). Stepanov (\cite{Sp2}) solved the case
$0<p_1\le 1<p_0\/$ and Sinnamon and Stepanov (\cite{SS}) solved the case $0<p_1<1=p_0\/$.

We shall also be interested in the  boundedness  of the discrete  Hardy operator $A_d\/$\index{Hardy
operator!discrete} acting on sequences  $f=\big(f(n)\big)_{n\ge0}\/$ in the form
$$
A_df(n)={1\over n+1}\sum_{k=0}^n f(k),\qquad n=0,1,2,\dots
$$
In the next chapter, it will be very much useful to know the boundedness of 
\bequ
\label{I.6.7}
A_d\,:\,\ell^p_{\hbox{\decp dec}}(w)\longrightarrow \ell^{p,\infty}(w), 
\eequ
and also,
\bequ
\label{I.6.8}
A_d\,:\,\ell^{p,\infty}_{\hbox{\decp dec}}(w)\longrightarrow \ell^{p,\infty}(w), 
\eequ
where
$w=\big(w(n)\big)_n\/$ is a weight in  $\N^\ast\/$, that is, a sequence of positive numbers, and
$\ell^p(w)\/$\index{$\ell^p(w)\/$}  is the  Lebesgue space $L^p\/$ in $\N^\ast\/$, with measure $\sum_n
w(n)\delta_{\{n\}}\/$, that is,
$$
\ell^p(w)=\bigg\{f=\big(f(n)\big)_n\subset\C\,:\,\|f\|_{\ell^p(w)}=\bigg(\sum_{n=0}^\infty|f(n)|^p
w(n)\bigg)^{1/p}<\infty\,\bigg\}. 
$$
$\ell^p_{\hbox{\decp dec}}(w)\/$\index{$\ell^p_{\hbox{\decp dec}}(w)$} is the class of positive and
decreasing sequences $f\/$  (we shall use the notation 
$f\d\/$) in $\ell^p(w)\/$ while  $\ell^{p,\infty}(w)\/$\index{$\ell^{p,\infty}(w)\/$} is the weak version
of  $\ell^p(w)\/$ and it is defined by the semi-norm
$$
\|f\|_{\ell^{p,\infty}(w)}=\sup_{t>0}t^{1/p}f^\ast_w(t),
$$ 
where $f^\ast_w\/$ is the decreasing rearrangement of  $f=\big(f(n)\big)_n\/$ 
in the measure space
$\big(\N^\ast,\sum_n w(n)\delta_{\{n\}}\big)\/$. Similarly as  was done in  $\Rm\/$, for each
weight 
$w\/$ (resp. $w_0,u,\dots\/$) in $\N^\ast\/$ we denote by  $W\/$ (resp. $W_0,U,\dots\/$)   the
sequence,
\bequ
\label{I.6.9}
W(n)=\sum_{k=0}^n w(k),\qquad n=0,1,2,\dots 
\eequ
With this notation, it can be easily seen that
\bequ
\label{I.6.10}
\|f\|_{\ell^{p,\infty}(w)}=\sup_{n\ge0}W^{1/p}(n)f(n),\qquad f\d. 
\eequ

It is now the moment to recall an important result due to    E. Sawyer
which  will be used several times in what follows. The original proof can be seen in \cite{Saw2}
and other proofs in  \cite{CS2} and \cite{Sp2}. 

\bteo
\label{te: I.5.7} 
If $1<p<\infty\/$ we have that
\begin{eqnarray*}
\sup_{f\downarrow}{\|f\|_{L^1(w_1)}\over\|f\|_{L^p(w_0)}}&\approx&\bigg(\int_0^\infty\bigg({W_1(t)\over
W_0(t)}\bigg)^{p^\prime-1}w_1(t)\,dt\bigg)^{1/p^\prime}\\ &\approx&\bigg(\int_0^\infty\bigg({W_1(t)\over
W_0(t)}\bigg)^{p^\prime}w_0(t)\,dt\bigg)^{1/p^\prime}+{W_1(\infty)\over W_0^{1/p}(\infty)}.
\end{eqnarray*}
\eteo
 
Let us now formulate a discrete version of the above result.

\bteo
\label{te: I.6.11}  Let $w=\big(w(n)\big)_n,\;v=\big(v(n)\big)_n\/$ be weights in $\N^\ast\/$ and let  
$$
S=\sup_{f\d}{\sum_{n=0}^\infty f(n)v(n)\over\big(\sum_{n=0}^\infty f(n)^p w(n)\big)^{1/p}}. 
$$
Then,
\begin{enumerate}
\item[(i)] If $0<p\le1\/$,
$$
S=\sup_{n\ge0}{V(n)\over W^{1/p}(n)},
$$ 
with $W\/$ defined by (\ref{I.6.9}) and  $V\/$ analogously.

\item[(ii)] If $1<p<\infty\/$,
\begin{eqnarray*}
S&\approx&\bigg(\int_0^\infty\bigg({\widetilde{V}(t)\over\widetilde{W}(t)}\bigg)^{p^\prime-1}
\tilde{v}(t)\,dt\bigg)^{1/p^\prime}\\ 
&\approx&\bigg(\int_0^\infty\bigg({\widetilde{V}(t)\over\widetilde{W}(t)}\bigg)^{p^\prime}
\tilde{w}(t)\,dt\bigg)^{1/p^\prime}+{\widetilde{V}(\infty)\over\widetilde{W}^{1/p}(\infty)},
\end{eqnarray*}
where $\tilde{v}\/$ is the weight in $\Rm\/$ defined by 
$$
\tilde{v}=\sum_{n=0}^\infty v(n)\chi_{[n,n+1)}
$$ 
and $\widetilde{V}(t)=\displaystyle\int_0^t\tilde{v}(s)\,ds\/$ and analogously for 
$\tilde{w}$ and $\widetilde{W}\/$.

Moreover, the constants implicit in the symbol $\approx\/$ only depend on $p\/$.
\end{enumerate}
\eteo

\bdem       (i) is obtained applying  Corollary~\ref{cor: I.2.14}  with
$p_1=1,\,p_0=p,\,X=Y=\N^\ast,\,T=\hbox{\rm Id}\/$ to the regular class $L\/$ of decreasing
sequences in 
$\N^\ast\/$.

(ii) can be deduced from Theorem~\ref{te: I.5.7}  observing that
\bequ
\label{I.6.12}
S=\sup_{\tilde{f}\downarrow}{\displaystyle\int_0^\infty\tilde{f}(t)\tilde{v}(t)\,dt\over\bigg(\displaystyle\int_0^\infty
\tilde{f}^p(t)\tilde{w}(t)\,dt\bigg)^{1/p}}.
\eequ
To see this, note that if  $f=\big(f(n)\big)_n\/$ is a decreasing sequence in  $\N^\ast\/$ and we
define 
$\tilde{f}=\sum_{n=0}^\infty f(n)\chi_{[n,n+1)}\in\Fd\/$, it is obvious that
$$
{\displaystyle\sum_{n=0}^\infty f(n)v(n)\over\bigg(\displaystyle\sum_{n=0}^\infty f(n)^p
w(n)\bigg)^{1/p}}={\displaystyle\int_0^\infty\tilde{f}(t)\tilde{v}(t)\,dt\over\bigg(\displaystyle\int_0^\infty
\tilde{f}^p(t)\tilde{w}(t)\,dt\bigg)^{1/p}}. 
$$ 
Therefore $S\/$ is less than or equal to the second member in 
 (\ref{I.6.12}). On the other hand, if $g\ge0\/$ is a  decreasing function     in $\Rm\/$ 
and we define the  decreasing sequence $f(n)=\big(\int_n^{n+1}g^p(s)\,ds\big)^{1/p},\;n=0,1,\dots\/$, we
obtain
$$
\int_0^\infty g^p(t)\tilde{w}(t)\,dt=\sum_n f(n)^p w(n),
$$ 
while by  H\"older's inequality
\begin{eqnarray*}
\int_0^\infty g(t)\tilde{v}(t)\,dt&=&\sum_n v(n)\int_n^{n+1}g(t)\,dt\\
&\le&\sum_n
v(n)\bigg(\int_n^{n+1}g^p(t)\,dt\bigg)^{1/p}=\sum_n v(n)f(n). 
\end{eqnarray*} 
Hence,
$$
{\displaystyle\int_0^\infty g(t)\tilde{v}(t)\,dt\over\bigg(\displaystyle\int_0^\infty
g^p(t)\tilde{w}(t)\,dt\bigg)^{1/p}}\le{\displaystyle\sum_{n=0}^\infty
f(n)v(n)\over\bigg(\displaystyle\sum_{n=0}^\infty f(n)^p w(n)\bigg)^{1/p}}\le S. 
$$ 
Thus
(\ref{I.6.12}) is proved  and  (ii) follows by applying  Theorem~\ref{te: I.5.7}.$\qquad\qed$
\edem

The following result characterizes the boundedness   (\ref{I.6.7}).

\bteo
\label{te: I.6.13}
\ \par
\begin{enumerate}
\item[(a)] If  $0<p\le1\/$ we have that $A_d\,:\,\ell^p_{\hbox{\decp dec}}(w)\longrightarrow
\ell^{p,\infty}(w)\/$ if and only if
\bequ
\label{I.6.14}
{W^{1/p}(n)\over n+1}\le C\,{W^{1/p}(m)\over m+1},\qquad 0\le m\le n.
\eequ

\item[(b)] If $1<p<\infty\/$ the following statements are equivalent:
\begin{enumerate}
\item[(i)] $ A_d\,:\,\ell^p_{\hbox{\decp dec}}(w)\longrightarrow \ell^{p,\infty}(w),$
\item[(ii)] $\tilde{w}=\sum_{n=0}^\infty w(n)\chi_{[n,n+1)}\in B_p,$
\item[(iii)] $\sum_{k=0}^n{1\over W^{1/p}(k)}\le C\,{n+1\over W^{1/p}(n)},\quad
n=0,1,2,\dots$
\end{enumerate}
\end{enumerate}
\eteo

\bdem      The boundedness  of $A_d\,:\,\ell^p_{\hbox{\decp dec}}(w)\longrightarrow
\ell^{p,\infty}(w)\/$ is equivalent (c.f. (\ref{I.6.10})) to the inequality 
$$
W^{1/p}(n)A_d f(n)\le C\|f\|_{\ell^p(w)},\qquad f\d,\,n=0,1,2,\dots
$$ 
Equivalently
\bequ
\label{I.6.15}
{W^{1/p}(n)\over n+1}\sup_{f\d}{\sum_{k=0}^n f(k)\over\big(\sum_{k=0}^\infty f(k)^p
w(k)\big)^{1/p}}\le C.
\eequ 
Hence, (a) is obtained applying  Theorem~\ref{te: I.6.11} (i). Also observe that the condition 
(\ref{I.6.14}) is obtained applying (\ref{I.6.15}) to sequences of the form
$f=(1,1,1,\dots,1,0,0,\dots)\/$ (characteristic functions in  $\N^\ast\/$)  and thus it is also
necessary in the case   $p>1\/$. 

To prove  (b), that is the case $p>1\/$, we note, as we already did in  Theorem~\ref{te: I.6.11},
that the discrete supremum in  (\ref{I.6.15}) can be substituted by a supremum on decreasing
functions in  
$\Rm\/$ with weights 
$\chi_{(0,n+1)}$ and $\tilde{w}\/$:
$$
\sup_{g\d}{\displaystyle\int_0^{n+1} g(t)\,dt\over\bigg(\displaystyle\int_0^\infty
g^p(t)\tilde{w}(t)\,dt\bigg)^{1/p}}. 
$$ 
Therefore
(\ref{I.6.15}) is equivalent to
$$
\widetilde{W}^{1/p}(n+1)Ag(n+1)\le C\|g\|_{L^p(\tilde{w})},\qquad g\d,\;n=0,1,2,\dots
$$ 
In fact
the previous inequality is true (with constant $2C\/$) substituting $n+1\/$ by any value 
$t>0\/$. To see this, if $0<t<1\/$,
\begin{eqnarray*}
\widetilde{W}^{1/p}(t)Ag(t)&=&\big(w(0)t\big)^{1/p}\,t^{-1}\int_0^t g(s)\,ds\\ 
&\le&\big(w(0)t\big)^{1/p}\,t^{-1}\bigg(\int_0^t g^p(s)\,ds\bigg)^{1/p}t^{1/p^\prime}\\ 
&\le&\bigg(\int_0^t
g^p(s)\tilde{w}(s)\,ds\bigg)^{1/p}\le\|g\|_{L^p(\tilde{w})},
\end{eqnarray*} 
and if  $t\in[n,n+1),\;n\ge1,\/$ we obtain, using  (\ref{I.6.14}),
$$
\widetilde{W}^{1/p}(t)Ag(t)\le\widetilde{W}^{1/p}(n+1)Ag(n)\le
2C\widetilde{W}^{1/p}(n)Ag(n)\le\|g\|_{L^p(\tilde{w})}. 
$$ 
Summarizing, in the case $p>1\/$ the
boundedness  of $A_d\,:\,\ell^p_{\hbox{\decp dec}}(w)\longrightarrow \ell^{p,\infty}(w)\/$ 
is equivalent to the inequality 
$$
\widetilde{W}^{1/p}(t)Ag(t)\le C\|g\|_{L^p(\tilde{w})},\qquad t>0,
$$ 
which holds for every function
$g\d\/$ in $\Rm\/$. But this means that we have the boundedness of the continuous Hardy operator 
$$
A\,:\,L^p_{\hbox{\decp dec}}(\tilde{w})\longrightarrow L^{p,\infty}(\tilde{w}),
$$ 
which is equivalent  (Definition~\ref{def: I.6.3}) to  $\tilde{w}\in B_{p,\infty}=B_p\/$ (the classes
$B_p\/$ and
$B_{p,\infty}\/$ coincide if  $p>1\/$, see \cite{Ne1}). We have then proved the equivalence (i)
and (ii). To see the equivalence with  (iii) let us first observe that this condition 
also implies (\ref{I.6.14}) since if  $m\le n\/$,
$$
{m+1\over W^{1/p}(m)}\le\sum_{k=0}^m{1\over W^{1/p}(k)}\le\sum_{k=0}^n{1\over W^{1/p}(k)}
$$
and by 
(iii), the last expression is bounded above by  $C(n+1)/W^{1/p}(n)\/$. Hence, in the proof of 
(ii)$\Leftrightarrow$(iii) we can assume that (\ref{I.6.14}) holds. With this (and taking into
account the definition of 
$\tilde{w}\/$) it is immediate to show the equivalence between the condition 
$\tilde{w}\in B_p\/$ (that is (ii)), that by Theorem~\ref{te: I.6.6}  (iv) is
$$
\int_0^r{1\over\widetilde{W}^{1/p}(t)}\,dt\le C\,{r\over\widetilde{W}^{1/p}(r)},\qquad r>0,
$$ 
and condition  (iii):
$$
\sum_{k=0}^n{1\over W^{1/p}(k)}\le C\,{n+1\over W^{1/p}(n)},\quad n=0,1,2,\dots
$$ 
It is enough to discretize the  integral and to note that 
$$
\int_n^{n+1}{1\over\widetilde{W}^{1/p}(t)}\,dt\approx{1\over W^{1/p}(n)},\quad n=0,1,2,\dots\qquad\qed
$$
\edem

The characterization of the  boundedness  (\ref{I.6.8}) can be seen very easily and it is  
equivalent to the condition (b)(iii) of the previous  theorem.

\bteo
\label{te: I.6.16} 
For $0<p<\infty\/$ we have that
$$
A_d\,:\,\ell^{p,\infty}_{\hbox{\decp dec}}(w)\longrightarrow \ell^{p,\infty}(w)
$$
if and only if
$$
\sum_{k=0}^n{1\over W^{1/p}(k)}\le C\,{n+1\over W^{1/p}(n)},\quad n=0,1,2,\dots
$$
\eteo

\bdem      The boundedness  in the statement is  equivalent to the inequality
\bequ
\label{I.6.17}
\|A_d f\|_{\ell^{p,\infty}(w)}\le C,\qquad f\d,\;\|f\|_{\ell^{p,\infty}(w)}\le1. 
\eequ
Now, if  $f\/$ is decreasing $\|f\|_{\ell^{p,\infty}(w)}=\sup_n W^{1/p}(n)f(n)\/$ (cf.\ 
(\ref{I.6.10})) and $\|f\|_{\ell^{p,\infty}(w)}\le1\/$ implies $f(n)\le W^{-1/p}(n),$ for every $ n\/$.
On the other hand the sequence  $W^{-1/p}\/$ is decreasing and with norm  1. Therefore 
$f=W^{-1/p}\/$ is the sequence that attains the maximum in the left hand side of  (\ref{I.6.17}) and the
characterization will be 
$$
\|A_d W^{-1/p}\|_{\ell^{p,\infty}(w)}=C<\infty. 
$$ 
Taking into account that  $A_d W^{-1/p}\/$ is decreasing  and  (\ref{I.6.10}) for the norm in 
$\ell^{p,\infty}(w)\/$ of such sequences,  we obtain the condition of the statement. $\qquad\qed$
\edem

\bcor
\label{cor: I.6.18} 
If $1<p<\infty\/$ the following conditions are equivalent:
\begin{enumerate}
\item[(i)] $ \tilde{w}=\sum_{n=0}^\infty w(n)\chi_{[n,n+1)}\in B_p,$
\item[(ii)] $ \sum_{k=0}^n{1\over W^{1/p}(k)}\le C\,{n+1\over W^{1/p}(n)},\quad n=0,1,2,\dots,$
\item[(iii)] $ A_d\,:\,\ell^p_{\hbox{\decp dec}}(w)\longrightarrow \ell^{p,\infty}(w),$ 
\item[(iv)] $ A_d\,:\,\ell^{p,\infty}_{\hbox{\decp dec}}(w)\longrightarrow \ell^{p,\infty}(w),$
\item[(v)] $ A_d\,:\,\ell^p_{\hbox{\decp dec}}(w)\longrightarrow \ell^p(w).$
\end{enumerate}
\ecor

\bdem      The first four equivalences are a consequence of the two previous theorems. On the other
hand, (v) implies (iii) and (i) implies
$$
A\,:\,L^p_{\hbox{\decp dec}}(\tilde{w})\longrightarrow L^p(\tilde{w}). 
$$ 
It is immediate to check that this  boundedness  of the  Hardy operator in  $\Rm\/$ implies the
corresponding boundedness for the discrete  Hardy operator 
$A_d\/$, that is (v). $\qquad\qed$
\edem

\chapter{Lorentz Spaces}
\label{ch:lor}
\pagestyle{myheadings}
\markboth{CHAPTER 2. LORENTZ SPACES}{M.J.\ Carro, J.A.\ Raposo, and J.\ Soria}

\section{Introduction}
\label{se:introlor}

If $(X,\mu)$ is a measure space, $w$ is a weight in $\R^+$ (see Definition \ref{defgen} below) and
$0<p<\infty$, the Lorentz space\index{Lorentz space} $\Lambda_X^p(w)$\index{$\Lambda_X^p(w)$} is
defined as the class of all measurable functions $f$ in $X$ for which
\bequ
\label{normx}
\Vert f \Vert_{\Lambda_X^p(w)}\buildrel{\hbox{\decp def.}}\over=\bigg(\int_0^{\infty}
(f^*(t))^p\,w(t)\,dt\bigg)^{1/p}<\infty,
\eequ
where $f^*$ is the nonincreasing rearrangement\index{nonincreasing rearrangement} of
$f$ with respect to $\mu$. These spaces were first introduced by G.G.\ Lorentz in
\cite{Lor} for the case $X=(0,1)$. By choosing $w$ properly, one obtains the
spaces $L^{p,q}$\index{$L^{p,q}$} defined in (\ref{normlor}). As we shall see, we also have a
weak-type version denoted by
$\Lambda^{p,\infty}_X(w)$\index{$\Lambda^{p,\infty}_X(w)$}.

In this chapter we are going to consider the study of analytical properties of these
function spaces, giving a complete description of some previously known partial
results. Already in the work of Lorentz (\cite{Lor}) there exists a characterization
of when (\ref{normx}) defines a norm. Later, A. Haaker (\cite{Haa})
extends this result to the case of $X=\R^+$ and $w\notin L^1$, and considers also
the existence of a dual space. E. Sawyer (\cite{Saw2}) gives
several equivalent conditions on $w$ so that $\Lambda_{\R^n}^p(w)$ is a Banach
space, when $p>1$, and M.J. Carro, A. Garc\'{\i}a del Amo, and J. Soria
(\cite{CGS}) study the normability in the case $p\ge1$ and $X$ is a nonatomic
space. For these same conditions, J. Soria (\cite{Sor}) characterizes the
normability of the weak-type Lorentz spaces\index{Lorentz space!weak-type}. We
present a throughout account of all these properties in the general setting of
resonant spaces, allowing us to also consider the discrete
case $X=\N$ and the sequence spaces $\ell^p$.  These discrete spaces (which are also
known as $d(w,p)$\index{$d(w,p)$}) had been previously studied for decreasing
weights $w$ (and hence $\Vert\cdot\Vert_{\Lambda^p(w)}$ is always a norm if
$p\ge 1$). Here we will consider general weights and also the weak-type spaces
$d^{\infty}(w,p)$\index{$d^{\infty}(w,p)$}.

This chapter is divided into several sections. After the definition and first
properties of section~\ref{se:def} we consider in section~\ref{se:cuasi}
the quasi-normed Lorentz spaces (for which only minimal assumptions on the
monotonicity of the primitive of $w$ are required). Density properties for simple
functions and absolute continuity of the quasi-norm are also considered.
Section~\ref{se:dual} is the largest and we study the duality results. In
particular we characterize  when $\Lambda^*=\Lambda'$ and when these are the
trivial space. In section~\ref{se:norma} we give necessary and sufficient
conditions to have that the Lorentz spaces are Banach spaces. Finally in
section~\ref{se:oper} we study interpolation properties and the boundedness of
certain operators.

We fix now the main definitions and notations used in what follows.
\bdefi
\label{defgen}
We denote by $\R^+=(0,\infty)$. The letters $w,\ \tilde w,\ w_0,...,$ are used for
weight\index{weight} functions in $\R^+$ (nonnegative locally integrable functions
in $\R^+$). For a given weight $w$ we write $W(r)=\displaystyle\int_0^rw(t)\,dt<\infty$, $0\le r<\infty$.

If $(X,\mu)$ is a measure space,  $f\in {\cal M}(X)$, we denote 
the distribution function\index{distribution function} and the nonincreasing
rearrangement\index{nonincreasing rearrangement} of
$f$ as
 $\lambda_f(t)$ and $f^*(t)$, and 
if the measure $d\mu(t)=w(t)\,dt$, we will  write $\lambda_f^w$ and $f^*_w$ to
show the dependence on the weight $w$.

 If
$(X,\mu)=(\R^m,w(t)dt)\/$ we write $L^{p,q}(w)$\index{$L^{p,q}(w)$} instead of $L^{p,q}(X)\/$.  
\edefi

\section{$\Lambda_X^p(w)$ spaces}
\label{se:def}

In this section $(X,\mu)\/$ denotes, except if otherwise mentioned, a general measure space.

\bdefi\label{def:norm}
Let $w\/$ be a weight in $\R^+$. For $0<p<\infty\/$ we define the functional
$\|\cdot\|_\loX:\M(X)\rightarrow[0,\infty]\/$ as
$$\|f\|_\loX=\bigg(\int_0^\infty (f^\ast(t))^p w(t)\,dt\bigg)^{1/p},\quad f\in \M(X).$$ The Lorentz
space\index{$\Lambda^p(w)$}
$\Lambda^p(w)=\loX\/$ is the class
$$\loX=\big\{f\in \M(X)\,:\,\|f\|_\loX<\infty\big\}.$$ Observe that $\|f\|_\loX=\|f^\ast\|_{L^p(w)}\/$. This
allows us to extend the previous definition. For
$0<p,q\le\infty\/$ set
$$\lox{p,q}{X}{w}\index{$\lox{p,q}{X}{w}$}=\big\{f\in
\M(X)\,:\,\|f\|_{\lox{p,q}{X}{w}}=\|f^\ast\|_{L^{p,q}(w)}<\infty\big\}.$$
\edefi

From now on, the notation $\loX\/$ or $\lox{p,q}{X}{w}\/$ without any reference to $w$, means that $w\/$ is a
weight in $\R^+$ not identically zero on $\big(0,\mu(X)\big)\/$. The symbol $\Lambda\/$ will denote any of the
spaces previously defined. Sometimes we will write $\Lambda^p\/$ as a shorthand for $\Lambda^p_X(w)\/$ if
$X\/$ and
$w\/$ are clearly determined. Similarly $\Lambda(w),\Lambda_X,\Lambda^p(w)\/$, etc.
Observe that
$\Lambda^{\infty,q}=\{0\}\/$ if $0<q<\infty\/$.

These spaces were introduced by Lorentz in \cite{Lor2} and \cite{Lor} for the case $X=(0,l)\subset\Rm\/$.
They are invariant under rearrangement\index{rearrangement invariant} and generalize the  $L^p\/$ 
Lebesgue spaces and 
$L^{p,q}\/$ (see Examples~\ref{Eje2.3}(i) and \ref{Eje2.3}(ii)). The spaces $\Lambda^p(w)\/$ we have defined
are usually called the \lq\lq classical\rq\rq\ Lorentz spaces  to distinguish them from $L^{p,q}\/$. 

\bobs\label{Obs2.2}

(i) $f^\ast(t)=0\/$ if $t\ge\mu(X)\/$. Hence, the behavior of the weight $w\/$ in
$\big[\mu(X),\infty\big)\/$ is irrelevant. We will assume, without loss of generality, that the weight $w\/$
vanishes in 
$\big[\mu(X),\infty\big)\/$ if $\mu(X)<\infty\/$. Observe that, in this way, we have that  $w\in
L^1(\Rm)\/$ if
$\mu(X)<\infty\/$.

(ii) If $w\notin L^1(\Rm)\/$ then $\lim_{t\to\infty}f^\ast(t)=0\/$ if $f\in\Lambda^{p,q}(w)\/$
($p<\infty\/$).

(iii) Simple functions with finite support are in $\Lambda^{p,q}(w)\/$. If $w\in L^1\/$ then
$L^\infty\subset\Lambda^{p,q}(w)\/$ and every simple function is in  $\Lambda^{p,q}(w)\/$.

(iv) Observe that $\lox{p,p}{}{w}=\Lambda^p(w),\,0<p\le\infty\/$.
\eobs

\beje\label{Eje2.3}

(i) In the case $w=1\/$ we have, for $0<p,q\le\infty,\;\Lambda^{p,q}_X(1)=L^{p,q}(X)\/$. Here
$W(t)=t,\;t\ge0\/$.

(ii) If $0<p,q<\infty,\;\Lambda^q_X(t^{q/p-1})=L^{p,q}(X)\/$ with  equality of \lq\lq norms\rq\rq. In this
case
$W(t)={p\over q}t^{q/p},\;t\ge0\/$.

(iii) If $w=\chi_{(0,1)}\/$, the space $\Lambda^1_X(w)=\Lambda^1_X\big(\chi_{(0,1)}\big)\/$ contains  
$L^\infty(X)\/$ and $W(t)=t,\;0\le t<1\/$. In this case the functional $\|\cdot\|_{\Lambda^1}\/$ is a norm
and the space is a 
\lq\lq  Banach function space" (see Definition~\ref{def4.3}).

(iv) If $X=\N^\ast=\{0,1,2,\dots\}\/$ and  we consider the counting measure, then measurable functions in
$X\/$ are sequences 
$f=\big(f(n)\big)_{n\ge0}\subset\C\;$ and
$$
\|f\|_{\Lambda_X^p(w)}=\bigg(\sum_{n=0}^\infty \big(f^\ast(n)\big)^p\Omega_n\bigg)^{1/p},
$$ 
where for each 
$n=0,1,2,\dots,\;\Omega_n=\displaystyle\int_n^{n+1} w(s)\,ds=W(n+1)-W(n)\/$. Thus, $\Lambda_X^p(w)\/$
depends only on the sequence of positive numbers $\Omega=(\Omega_n)_{n=0}^\infty\/$
 and it is usually denoted as 
$d(\Omega,p)\/$. In the weak-type case,
$$
\|f\|_{\Lambda_X^{p,\infty}(w)}=\sup_{n\ge0}\bigg(\sum_{k=0}^n\Omega_k\bigg)^{1/p}f^\ast(n),
$$ 
and we will use the symbol $d^\infty(\Omega,p)\/$ to denote
$\Lambda_X^{p,\infty}(w)\/$. Given $d(\Omega,p)\/$ or $d^\infty(\Omega,p)\/$
we will assume that $\Omega\/$ is always a sequence of positive numbers.

It is clear that both $d(\Omega,p)\/$ and $d^\infty(\Omega,p)\/$ are always contained in
$\ell^\infty(\N^\ast)\/$ and if 
$\Omega\in \ell^1(\N^\ast)\/$ we have in fact that $d(\Omega,p)=d^\infty(\Omega,p)
=\ell^\infty(\N^\ast)\/$ (with
equivalent norms). The only interesting case is hence  $\Omega\notin\ell^1(\N^\ast)\/$. 
In this case, the
space is contained in 
$c_0(\N^\ast)\/$, the space of sequences that converge to $0\/$ and, for each $f\in 
d(\Omega,p)\/$ ($f\in
d^\infty(\Omega,p)\/$), 
$\big(f^\ast(n)\big)_{n\ge0}\/$ is the nonincreasing rearrangement of the sequence
$\big(f(n)\big)_{n\ge0}\/$.

If the sequence $\Omega=(\Omega_n)_n\/$ of the previous example is decreasing and $p\ge1\/$, 
 $d(\Omega,p)\/$
is a Banach space. This case has been studied by several authors (\cite{Ga}, \cite{CH}, 
\cite{Al}, \cite{Po},
\cite{NO}, \cite{AM2}, etc.). The case  $\Omega\/$ increasing is considered in \cite{AEP}. 
 We will assume no a
priori conditions on $\Omega$.
\eeje

The following lemma (see \cite{CS1}) will be very useful.

\blem
\label{le: II.2.4}
 Let $0<p<\infty\/$ and $v\ge0\/$ be a measurable function in $\R^+\/$. Let $V(r)=\displaystyle\int_0^r
v(s)\,ds,\,0\le r<\infty\/$. Then, for every decreasing function $f\/$ we have that 
$$
\int_0^\infty f^p(s) v(s)\,ds=\int_0^\infty pt^{p-1}V\big(\lambda_f(t)\big)\,dt.
$$
\elem

The following result gives several equivalent expressions for the functional 
 $\|\cdot\|_{\lox{p,q}{X}{w}}$.  In particular, we see that it only depends on  $W\/$.

\bpro
\label{pro: II.2.5}
For $0<p,q<\infty\/$ and  $f\/$ measurable in $X\/$,
\begin{enumerate}
\item[(i)] $\|f\|_\lox{p,q}{X}{w}=\bigg(\displaystyle\int_0^\infty pt^{q-1}W^{q/p}
\big(\lambda_f(t)\big)\,dt\bigg)^{1/q},$
\item[(ii)] $\|f\|_\lox{p}{X}{w}=\bigg(\displaystyle\int_0^\infty pt^{p-1}W\big(\lambda_f(t)\big)\,dt
\bigg)^{1/p},$
\item[(iii)] $\|f\|_\loinx{p}{X}{w}=\sup_{t>0}tW^{1/p}\big(\lambda_f(t)\big)=\sup_{t>0}
f^\ast(t)W^{1/p}(t).$
\end{enumerate}
\epro

\bdem

(i) Since every function and its decreasing rearrangement have the same distribution function (see
\cite{BS}) we have that,
$$
W\big(\lambda_f(t)\big)=W\big(\lambda_{f^\ast}(t)\big)=\int_0^{\lambda_{f^\ast}(t)}w(s)\,ds
=\int_{\{f^\ast>t\}}w(s)\,ds=\lambda^w_{f^\ast}(t).
$$
By Lemma~\ref{le: II.2.4} and since  $\lambda^w_{f^\ast}\/=(f^\ast)^\ast_w\/$, 
we obtain that
\begin{eqnarray*}
\bigg(\int_0^\infty pt^{q-1}W^{q/p}\big(\lambda_f(t)\big)\,dt\bigg)^{1/q}&=&\bigg(\int_0^\infty 
pt^{q-1}\big(\lambda^w_{f^\ast}(t)\big)^{q/p}\,dt\bigg)^{1/q}\\
&=&\bigg(\int_0^\infty t^{q/p-1}\int_0^{(f^\ast)^\ast_w(t)}qs^{q-1}\,ds\,dt\bigg)^{1/q}\\ &
=&\bigg(\int_0^\infty\big((f^\ast)^\ast_w(t)\big)^q
t^{q/p-1}\,dt\bigg)^{1/q}\\ &=&\|f\|_{\lox{p,q}{X}{w}}.
\end{eqnarray*}

(ii) It is an immediate consequence of  (i).

(iii) For the first inequality we observe that 
$$
\|f\|_{\loinX}=\|f^\ast\|_{L^{p,\infty}(w)}=\sup_{t>0}t\big(\lambda^w_{f^\ast}(t)\big)^{1/p}
=\sup_{t>0}tW^{1/p}\big(\lambda_f(t)\big).
$$ 
The second inequality for characteristic functions is trivial and, by monotonicity, the
general case follows.$\qquad\qed$
\edem

\bobs
\label{obs: II.2.6}

(i)  If we compare  \ref{pro: II.2.5}   (i) and  \ref{pro: II.2.5} (ii) we see that, for
$q<\infty\/$, 
$\|f\|_{\lox{p,q}{X}{w}}=\|f\|_\lox{q}{X}{\tilde{w}}\/$ where
$\tilde{w}(t)=W^{q/p-1}(t)w(t),\;0<t<\mu(X)\/$. Therefore, every  Lorentz space as defined here
reduces to 
$\loX\/$ and its weak version $\loinX\/$.

(ii) From \ref{pro: II.2.5} (iii), we deduce that  $\lox{p,\infty}{X}{w}=
\Lambda^{q,\infty}_X\big((q/p)\tilde{w}\big)\/$ for $0<p,q<\infty\/$.
\eobs

For the spaces $L^{p,\infty}(X)\/$ it is known that the quasi-norm $\|f\|_{p,\infty}\/$ 
is, for every $q<p\/$, equivalent to the functional 
$$
\sup_{E\subset X}\big\|f\chi_E\big\|_q\,\mu(E)^{1/p-1/q}.
$$ This is the so-called  Kolmogorov condition
(see e.g.\ \cite{GR}). An analogous version for  
$\Lambda^{p,\infty}(w)\/$ also holds.

\bpro
 If $0<q<p<\infty\/$ and $f\in\M(X)\/$,
$$
\|f\|_{\Lambda^{p,\infty}_X(w)}\le\sup_{E\subset X}\big\|f\chi_E\big\|_{\Lambda^q_X(w)}W
\big(\mu(E)\big)^{1/p-1/q}\le\bigg({p\over
p-q}\bigg)^{1/q}\,\|f\|_{\Lambda^{p,\infty}_X(w)},
$$ 
where the supremum is taken over all measurable sets $E\subset X\/$.
\epro

\bdem
Let 
$$
S=\sup_{E\subset X}\big\|f\chi_E\big\|_{\Lambda^q_X(w)}W\big(\mu(E)\big)^{1/p-1/q}.
$$ 
To show the first inequality,  let  $E=\{|f|>t\}\/$, $t>0\/$.
Then,
\begin{eqnarray*}
S&\ge&\big\|f\chi_E\big\|_{\Lambda^q_X(w)}W\big(\mu(E)\big)^{1/p-1/q}\\ &
=&\bigg(\int_0^\infty\big(\big(f\chi_E\big)^\ast(s)\big)^q
w(s)\,ds\bigg)^{1/q}W\big(\mu(E)\big)^{1/p-1/q}\\ 
&\ge&\bigg(t^q\int_0^{\mu(E)} w(s)\,ds\bigg)^{1/q}
W\big(\mu(E)\big)^{1/p-1/q}\\
&=&t\,W\big(\mu(E)\big)^{1/p}=t\,W\big(\lambda_f(t)\big)^{1/p},
\end{eqnarray*} 
and taking the supremum in $t>0\/$
we get $\|f\|_{\Lambda^{p,\infty}_X(w)}\le S\/$ (c.f.
\ref {pro: II.2.5} (iii)).

To prove the second inequality, for each  $f\in\M(X),\,E\subset X\/$, let $a=
\|f\|_{\Lambda^{p,\infty}_X(w)}W\big(\mu(E)\big)^{-1/p}\/$. Then, by \ref{pro: II.2.5} (iii),
\begin{eqnarray*}
\big\|f\chi_E\big\|^q_{\Lambda^q_X(w)}&=&\int_0^\infty qt^{q-1}W\big(\lambda_{f\chi_E}(t)\big)\,dt
\\ &=&\int_0^a
qt^{q-1}W\big(\lambda_{f\chi_E}(t)\big)\,dt+\int_a^\infty qt^{q-1}W\big(\lambda_{f\chi_E}(t)
\big)\,dt\\ &\le &W\big(\mu(E)\big)\int_0^a qt^{q-1}\,dt+\int_a^\infty q
t^{q-1}{\|f\|^p_{\Lambda^{p,\infty}_X(w)}\over t^p}\,dt\\ &=&{p\over p-q}\,
\|f\|^q_{\Lambda^{p,\infty}_X(w)}W\big(\mu(E)\big)^{1-q/p}.
\end{eqnarray*}
 Hence,
$$
\big\|f\chi_E\big\|_{\Lambda^q_X(w)}W\big(\mu(E)\big)^{1/p-1/q}\le\bigg({p\over p-q}\bigg)^{1/q}
\,\|f\|_{\Lambda^{p,\infty}_X(w)}.\quad\qed
$$
\edem

In the following proposition we state some elementary properties for these spaces (see \cite{BS}).

\bpro
\label{pro: II.2.8}
 For $0<p,q\le\infty\/$ and $f,g,f_k,\;k\ge1,\/$ measurable functions in  $X\/$, we have that:
\begin{enumerate}
\item[(i)]  $|f|\le|g|\Rightarrow\|f\|_{\lox{p,q}{X}{w}}\le\|g\|_{\lox{p,q}{X}{w}}\/$,

\item[(ii)] $\|tf\|_{\lox{p,q}{X}{w}}=|t|\|f\|_{\lox{p,q}{X}{w}},\quad t\in\C\/$,

\item[(iii)]  $0\le f_k\le f_{k+1}\limto{k}f\/$ a.e. $\Rightarrow\|f_k\|_{\lox{p,q}{X}{w}}
\limto{k}\|f\|_{\lox{p,q}{X}{w}}\/$,

\item[(iv)]  $\big\|\liminf_k|f_k|\big\|_{\lox{p,q}{X}{w}}\le\liminf_k\|f_k\|_{\lox{p,q}{X}{w}}\/$,

\item[(v)]  $\lox{p,q_0}{X}{w}\subset\lox{p,q_1}{X}{w}$ continuously, $0<q_0\le q_1\le\infty\/$,

\item[(vi)]  If\ $W\big(\mu(X)\big)<\infty\/$ then, for $0<p_0<p_1\le\infty\/$, 
we have that  $\lox{p_1,q}{X}{w}\subset\lox{p_0,r}{X}{w},\;0<r\le\infty\/$ continuously,

\item[(vii)]  $\chi_E\in\lox{p,q}{X}{w}\/$ if $\mu(E)<\infty\/$.
\end{enumerate}
\epro

The following property connects the norm convergence $\|\cdot\|_{\lox{p,q}{X}{w}}\/$ 
with the convergence  in measure  and it is related to the  \lq\lq completeness\rq\rq\ property of
our spaces.  

\bpro
\label{pro: II.2.9}
Assume that  $W>0\/$ in $(0,\infty)\/$, let  $\Lambda=\lox{p,q}{X}{w}\/$ be
a Lorentz space and  let $(f_n)_n\/$ be a sequence of measurable functions in  $X\/$.
\begin{enumerate}
\item[(i)] If   $\lim_{m,n}\|f_m-f_n\|_{\Lambda}=0\/$ then $(f_n)_n\/$  is a Cauchy sequence in measure
and there exists 
$f\in
\M(X)\/$ such that
$\lim\|f-f_n\|_{\Lambda}=0\/$.

\item[(ii)] If  $f\in \M(X)\/$ and  $\lim_{n}\|f-f_n\|_{\Lambda}=0\/$ then $(f_n)_n\/$ 
converges to $f\/$ in measure and there exists a partial $(f_{n_k})_k\/$ convergent to
$f\/$ a.e.
\end{enumerate}
\epro

\bdem
The case  $q<p=\infty\/$ is  trivial and in  the case $p=q=\infty,\;\Lambda^{p,q}
=L^\infty\/$ and the result is already known. If $p<\infty\/$ it is immediate,
by Proposition~\ref{pro: II.2.5}, that
$$
W\big(\lambda_f(r)\big)\le C_{p,q}{\|f\|_{\lox{p,q}{X}{w}}^p\over r^p},\quad r>0,
\quad0<q\le\infty.
$$ 
Using  the hypothesis of  (i) we obtain then that
$W\big(\lambda_{f_n-f_m}(r)\big)\limto{m,n}0\/$,   for every $r>0\/$, which (since $W>0\/$) implies
$\lambda_{f_n-f_m}(r)\limto{m,n}0,\;r>0\/$, that is, $(f_n)_n\/$ is a  Cauchy sequence 
in measure. We know  that this implies the convergence in measure of $(f_n)_n\/$ 
to some measurable function $f\/$ and the existence of a partial $(f_{n_k})_k\/$ converging to $f\/$
a.e. By Proposition~\ref{pro: II.2.8} (iv) we have that  $\|f-f_n\|_\Lambda
\le\liminf_{k}\|f_{n_k}-f_n\|_\Lambda\/$ and, thus, 
$\lim_n\|f-f_n\|_{\Lambda}=0\/$.

The proof of  (ii) is analogue.$\qquad\qed$
\edem

The functional $\|\cdot\|_{\Lambda}\/$ is not, in general, a quasi-norm and, in fact,  
$\Lambda\/$ could not even be a vector space. The following lemma characterizes the
quasi-normability of these spaces, which, as we will see, only depends on  the
weight $w\/$ and on the measure space $X\/$.

\blem
\label{le: II.2.10}
 If $0<p<\infty$ and $0<q\le\infty\/$, the space $\Lambda^{p,q}_X(w)\/$ 
is quasi-normed if and only if
\bequ\label{delta2}
0<W\big(\mu(A\cup B)\big)\le C\big( W\big(\mu(A)\big)+W\big(\mu(B)\big) \big),
\eequ for every pair of measurable sets $A,B\subset X\/$ with $\mu(A\cup B)>0\/$.
\elem

\bdem
Sufficiency:  The hypothesis implies that  $W\big(\mu(A)\big)>0\/$ if $\mu(A)>0\/$. 
If $\|f\|_{\Lambda^{p,q}}=0\/$, by  
Proposition~\ref{pro: II.2.5}, we have that
$W\big(\lambda_f(t)\big)=0,\;t>0,\/$ and hence, $\lambda_f(t)=0\/$ for every $t>0\/$, 
that is $f=0\/$ a.e. It remains to show the quasi-triangular inequality and it   suffices to prove
it for positive functions.  Let   $0\le
f,g\in\Lambda^{p,q}(w)\/$  and $t>0\/$. Then $\{f+g>t\}\subset\{f>t/2\}\cup\{g>t/2\}\/$  and by
hypothesis 
$$
W\big(\lambda_{f+g}(t)\big)\le C\big(W\big(\lambda_f(t/2)\big)+W\big(\lambda_g(t/2)\big)\big).
$$ 
Since $C\/$ does not depend on  $t\/$, it satisfies that (Proposition~\ref{pro: II.2.5})  
$$
\|f+g\|_{\Lambda^{p,q}(w)}\le C_{p,q}\big(\|f\|_{\Lambda^{p,q}(w)}+\|g\|_{\Lambda^{p,q}(w)}
\big).
$$ 
Necessity: If $A,B\/$ are two measurable sets with
$\mu(A\cup B)>0,\;\chi_{A\cup B}\le\chi_A+\chi_B\/$ and since  $\Lambda^{p,q}(w)\/$ 
is quasi-normed, we have that (Proposition~\ref{pro: II.2.5}),
\begin{eqnarray*}
 0<C_{p,q}W^{1/p}\big(\mu(A\cup B)\big)&=&\|\chi_{A\cup B}\|_{\Lambda^{p,q}(w)}\\ 
&\le&\|\chi_A+\chi_B\|_{\Lambda^{p,q}(w)}\\ &\le&
C\big(\|\chi_A\|_{\Lambda^{p,q}(w)}+\|\chi_B\|_{\Lambda^{p,q}(w)}\big)\\ &=&C^\prime_{p,q}
\big(W^{1/p}\big(\mu(A)\big)+W^{1/p}\big(\mu(B)\big)\big)
\end{eqnarray*}
 which is equivalent to the condition of the statement.$\quad\qed$
\edem

The previous result motivates the following definition.

\bdefi
 Let $w\/$ be a weight in $\Rm\/$. We write $W\in\Delta_2(X)\/$\index{$\Delta_2(X)$}
(or $W\in\Delta_2(\mu)\/$)\index{$\Delta_2(\mu)$} if $W\/$ satisfies (\ref{delta2}). If
$X=\R\/$ we shall simply write $W\in\Delta_2\/$\index{$\Delta_2$}.
\edefi

Therefore, Lemma~\ref{le: II.2.10} tells us that, for every  $0<p<\infty,\;0<q\le\infty\/$,
$$
\Lambda^{p,q}_X(w)\hbox{\rm is quasi-normed\ }\Leftrightarrow W\in\Delta_2(X).
$$

\bpro 
\label{pro: II.2.12}
The following conditions are equivalent:
\begin{enumerate}
\item[(i)] $W\in\Delta_2\/$,

\item[(ii)] $W(2r)\le CW(r),\;r>0\/$,

\item[(iii)] $W(t+s)\le C\big(W(t)+W(s)),\;t,s>0\/$,
\end{enumerate}

\noindent and in any of these cases, $W(t)>0,\;t>0\/$.
\epro

This proposition (whose proof is trivial) shows that $W\in\Delta_2$ is something easy to check. The condition
$W\in\Delta_2\/$ is sufficient to have the quasi-normability of 
$\Lambda^{p,q}_X(w)\/$,  independent of the measure space $X\/$ and  this is the content of
the following theorem. The second part was proved in \cite{CS1} and the proof of (i) is immediate.

\bteo 
\label{te: II.2.13}
Let $0<p<\infty,\,0<q\le\infty\/$.
\begin{enumerate}
\item[(i)] If $W\in\Delta_2,\;\Lambda^{p,q}(w)\/$ is quasi-normed.

\item[(ii)]  If $X\/$ is nonatomic, $\lox{p,q}{X}{w}\/$ is quasi-normed if and only if 
$W\in\Delta_2\/$.
\end{enumerate}
\noindent That is,  $\Delta_2\subset\Delta_2(X)\/$ for every $X\/$ and if  $X\/$ 
is nonatomic, $\Delta_2=\Delta_2(X)\/$ (we are  assuming that  $w\/$ is zero outside of 
$\big(0,\mu(X)\big))\/$.

\eteo

\bobs

For the case $p=q=\infty\/$, it  is obvious that the following properties are equivalent:
\begin{enumerate}
\item[(i)] $\Lambda^\infty_X(w)$  is quasi-normed.
\item[(ii)] $\Lambda^\infty_X(w) =L^\infty(X)$ (with equality of norms). 
\item[(iii)] $\big(\Lambda^\infty_X(w),\|\cdot\|_{\Lambda^\infty_X(w)}\big)$
is a Banach space.
\item[(iv)] $ W\big(\mu(A)\big)>0$ if $\mu(A)>0,\;A\subset X.$
\end{enumerate}
Any of these conditions hold if  $W\in\Delta_2(X)\/$.
\eobs

Let us recall the following definition that we shall use quite often in what follows.

\bdefi
\label{def: II.2.15}
A resonant measure space\index{resonant space} $X\/$  (\cite{BS}) 
is a $\sigma$-finite space  which is nonatomic or it is a union (at most  countable)  of atoms
with equal measure. 
\edefi

If $X\/$ is nonatomic we  have already characterized (Theorem~\ref{te: II.2.13}) 
when $\Lambda_X\/$ is quasi-normed. The following theorem completes the characterization for the
case 
$X\/$ is a resonant measure space. Its proof is immediate.

\bteo
\label{te: II.2.16}
Let $X\/$ be an atomic measure space whose atoms have  measure  $b>0\/$. 
Then $W\in\Delta_2(X)\/$ if and only if 
$$
W\big(2nb)\le C\,W(nb),\qquad n\ge1.
$$ 
In particular, if $0<p<\infty\/$  the spaces  $d(\Omega,p)\/$ and $d^\infty(\Omega,p)\/$
(Examples~\ref{Eje2.3}(iv)) are quasi-normed if and only if
$$
\sum_{n=1}^{2N}\Omega_{n-1}\le C\sum_{n=1}^N\Omega_{n-1},\quad N=1,2,\dots
$$
\eteo

\section{Quasi-normed Lorentz spaces}
\label{se:cuasi}

In this section $(X,\mu)\/$ will be   an arbitrary
$\sigma$-finite measure space and the  Lorentz spaces will be defined on the measure space   $X\/$.
We shall study the topology and some elementary properties of the quasi-normed 
Lorentz spaces. As we know,  the condition to have that 
$\Lambda_X(w)\/$ is quasi-normed is  $W\in\Delta_2(X)\/$ (Lemma~\ref{le: II.2.10}).

As in every quasi-normed space, we consider in  $\Lambda\/$ 
the topology induced by the quasi-norm $\|\cdot\|_\Lambda\,:\;G\subset\Lambda\/$ is an open set
if for every function 
$f\in G\/$ there exists 
$r>0\/$ such that  $\big\{g\in\Lambda\,:
\,\|f-g\|_\Lambda<r\big\}\subset G\/$. We know that there exists in $\Lambda\/$ a translation
invariant distance  $d\/$
  and an exponent  $p\in(0,1]\/$ so that 
\bequ
\label{II.3.1}
d(f,g)\le\|f-g\|_\Lambda^p\le 2 d(f,g),\quad f,g\in\Lambda. 
\eequ
This is a common fact in every quasi-normed space (see \cite{BL}). Inequality (\ref{II.3.1}) tells
us that the topology of 
$\Lambda\/$ coincides with the one induced by the metric $d\/$ and hence,
 $\Lambda\/$ is a metric space. From Proposition~\ref{pro: II.2.9}  it can be also deduced that it
is complete.
$\Lambda\/$ is then a quasi-Banach space. These and other facts are summarized in the following
theorem. Its proof is obvious using the previous considerations and Proposition~\ref{pro: II.2.9}
(see also  \cite{BL} for property (iii) and \cite{BS} for  (iv)).
\bigbreak
\bteo
\label{te: II.3.2}
\ 
Every quasi-normed  Lorentz space $\Lambda\/$ is quasi-Banach. 
In particular $\Lambda\/$ is an  F-space (topological vector space which is 
metrizable with a translations invariant measure and complete). Every function in  $\Lambda\/$ is
finite a.e. and if 
$(f_n)_n\subset\Lambda\/$,
\begin{enumerate}
\item[(i)]
$(f_n)_n$ is a Cauchy sequence, if and only if  $\lim_{m,n}\|f_n-f_m\|_\Lambda=0\/$ and then it is also a
Cauchy sequence in measure.

\item[(ii)] $\Lambda-\lim f_n=f\/$ if and only if  $\lim_{n}\|f-f_n\|_\Lambda=0\/$ and then  $(f_n)_n\/$
converges to 
$f\/$ in measure  and there exists a partial that converges to 
$f\/$ a.e.

\item[(iii)] If  $p=\log(2)/\log(2C)\/$ where $C\/$ is the  constant of 
the quasi-norm $\|\cdot\|_\Lambda\/$
in $\Lambda\/$, there exists a $p$-norm in $\Lambda\/$ equivalent to 
$\|\cdot\|_\Lambda^p\/$. If $\sum_n\|f_n\|_\Lambda^p<\infty\/$ the  series $\sum_n f_n\/$ 
is convergent in $\Lambda\/$.

\item[(iv)] If $\Lambda_X\subset\widetilde{\Lambda}_X\/$ and both are quasi-normed  Lorentz  
spaces, the embedding is continuous.

\item[(v)] If $F\/$ is another quasi-normed space,  a linear operator $T:\Lambda\to F\/$ 
is continuous if and only if   $\sup_{\|f\|_\Lambda\le1}\|Tf\|_F<\infty\/$.
\end{enumerate}
\eteo

\subsection{Absolutely continuous norm} 

We now study the equivalent property to the dominated convergence theorem  of the   $L^p\/$
spaces: If 
$\lim_n f_n(x)=f(x)\hbox{\ a.e.\ }x\/$ and 
$|f_n|\le g\in L^p\/$ then  $\|f_n-f\|_p\limto{n}0\/$. In general, in a Banach function   space
(on  
$X\;\sigma$-finite)  a function $g\/$ as before  is said to have   an  absolutely
continuous norm. A space in which every function has    absolutely
continuous norm satisfies the dominated convergence theorem (see \cite{BS}).  

\bdefi
 Let $\big(E,\|\cdot\|\big),\;E\subset\M(X),\/$ be a quasi-normed
space. A function $f\in E\/$   is said to have     absolutely
continuous norm\index{absolutely
continuous norm} if
$$
\lim_{n\to\infty}\big\|f\chi_{A_n}\big\|=0,
$$ 
for every decreasing sequence of 
 measurable sets $(A_n)_n\/$ with $\chi_{A_n}\to0\/$ a.e. If every function in 
$E\/$ has this property, we say that  $E\/$ has   an absolutely
continuous norm.
\edefi

The connection of this property with the dominated convergence theorem is clear from the following
proposition.
Its proof is 
(except on some minor modifications) identical to the one in  \cite{BS} (pp. 14--16) for
 Banach function spaces and we omit it.

\bpro
If $\Lambda\/$ is a quasi-normed  Lorentz space and   $f\in\Lambda\/$, the following statements are 
equivalent:
\begin{enumerate}
\item[(i)] $f\/$ has  absolutely continuous norm.

\item[(ii)] $\lim_n\big\|f\chi_{E_n}\big\|_\Lambda=0\/$ if $(E_n)_n\/$ is a sequence of measurable sets
with
$\chi_{E_n}\to0\/$ a.e.

\item[(iii)] $\lim_n\|f_n\|_\Lambda=0\/$ if $|f_n|\le|f|\/$ and $\lim_n f_n=0\/$ a.e.

\item[(iv)] $\lim_n\|g-g_n\|_\Lambda=0\/$ if $|g_n|\le|f|\/$ and $\lim_n g_n=g\/$ a.e.
\end{enumerate}
\epro

The following result shows that, except in the special case $\mu(X)=\infty\/$ and  $w\in L^1\/$, 
the spaces  $\loX\/$ have   absolutely continuous norm.

\bteo
\label{te: II.3.5} 
Let $0<p<\infty\/$ and  $\Lambda^p_X(w)\/$ be a quasi-normed space. Then,
\begin{enumerate}
\item[(i)] If $\mu(X)<\infty,\;\Lambda^p_X(w)\/$ has   absolutely continuous norm.

\item[(ii)] If $\mu(X)=\infty,\;\Lambda^p_X(w)\/$ has   absolutely continuous norm if and only if 
 $w\notin L^1\/$.
\end{enumerate}
\eteo

\bdem
 (i) Let us assume  that $\mu(X)<\infty\/$ and let us  see that if $0\le f\in\loX\/$ 
and  $(g_n)_n\/$ is a sequence of measurable functions  satisfying
$$
0\le g_n\le f\in\Lambda,\quad g_n\to0\hbox{\rm \  a.e.},
$$ 
then   $\lim_n\|g_n\|_\Lambda=0\/$.
Since the space is of finite measure, the pointwise convergence implies the convergence in measure 
  and we have that  $\lim_n\lambda_{g_n}(t)=0,\;0<t<\infty\/$. 
In particular $W\big(\lambda_{g_n}(t)\big)\to0,\;t>0\/$. Since
$W\big(\lambda_{g_n}(t)\big)\le W\big(\lambda_f(t)\big)\/$ and $f\in\Lambda^p\/$, 
from Proposition~\ref{pro: II.2.5}(ii) and the dominated convergence theorem, it follows that 
$\lim_n\|g_n\|_\Lambda=0\/$.

(ii) Assume that  $w\notin L^1,\;0\le f\in\loX\/$ and  $(g_n)_n\/$ is as above. 
The hypothesis on $w\/$ implies $\lambda_f(t)<\infty,\;t>0\/$ and the sets
$E_k=\{f\le1/k\}\/$ have complement of finite measure and, if $f_k=f\chi_{E_k}\/$ 
we have that $f_k^\ast\le1/k\/$ and also $f_k^\ast\le f^\ast\/$. By the 
dominated convergence theorem $\|f_k\|_\Lambda\to0\/$.
Therefore, given $\epsilon>0\/$ there exists a measurable set  $E\subset X\/$ with 
$\mu(X\setminus E)<\infty\/$ and such that 
$\|f\chi_E\|_\Lambda<\epsilon\/$. Then the functions  $g_n\chi_{X\setminus E}\/$ 
are as in  (i) and  $\|g_n\chi_{X\setminus E}\|_\Lambda\to0\/$.
We also have that 
$$
\limsup_n\|g_n\|_\Lambda\le C\limsup_n\big(\|g_n\chi_{X\setminus E}\|_\Lambda
+\|g_n\chi_E\|_\Lambda\big)\le C\|f\chi_E\|_\Lambda\le C\epsilon.
$$ 
And, since this is true for every  $\epsilon>0\/$, we obtain that $\lim_n\|g_n\|_\Lambda=0\/$.
This proves the sufficiency in (ii). To show the necessity, let us assume that  $\mu(X)=\infty,\;w\in
L^1\/$ and let us see that 
$\loX\/$ has not  absolutely continuous norm.  In this case, since  $X=\bigcup_{n=1}^\infty X_n\/$ 
with
$\mu(X_n)<\infty\/$ for every 
$n\/$, the sets 
$E_n=\bigcup_{k=n}^\infty X_k\/$ satisfy $\chi_{E_n}\to0\/$ a.e. and also, 
$\chi_{E_n}^\ast=1\/$. Thus, $\lim_n\|\chi_{E_n}\|_{\Lambda^p}=\|w\|_1^{1/p}>0\/$ and the function 
$1\in\loX\/$ does not have   absolutely continuous norm.$\quad\qed$
\edem

\bcor
\label{co: II.3.6} 
If $0<p<\infty\/$ and $\lox{p}{X}{w}\/$ is a quasi-normed space, 
every function in this space which vanishes in the complement of a set of finite measure, has 
absolutely continuous norm.
\ecor

The analogous question in the case of the weak-type space  $\loinX\/$ is solved 
 in the following theorem, for the case when  $X\/$ is a resonant measure space: 
these spaces do not have, except in the trivial case,   absolutely continuous norm.

\bteo
\label{te: II.3.7} 
If $X\/$ is a resonant measure space, $0<p<\infty\/$ and  $W\in \Delta_2(X),\;\loinX\/$
has    absolutely continuous norm if and only if  $X\/$ is a finite union of atoms.
\eteo

\bdem 
If $X\/$ is a finite union of atoms the proof is trivial.
Hence, we shall assume that this is not the case and we consider two possibilities:

(i) If  $X\/$ is nonatomic, since the function $W^{-1/p}\/$ is decreasing and continuous in 
$\R^+\/$,
there exists a measurable function
$f\ge0\/$ in $X\/$ such that
$f^\ast(t)=W^{-1/p}(t),\;0<t<\mu(X)\/$. Therefore
$\|f\|_{\Lambda^{p,\infty}(w)}=\sup_t W^{1/p}(t)f^\ast(t)=1\/$ and  $f\in\Lambda^{p,\infty}(w)\/$.
The sets
$E_n=\chi_{\{f>n\}},\;n=1,2,\dots\/$ form a decreasing sequence with  $\chi_{E_n}\to0\/$ a.e. and,
however,
$\|f\chi_{E_n}\|_{\Lambda^{p,\infty}}=1\/$ for every
$n\/$, since 
$$
\big(f\chi_{E_n}\big)^\ast(t)=f^\ast(t)\chi_{\big[0,\lambda_f(n)\big)}(t)
=W^{-1/p}(t)\chi_{\big[0,\lambda_f(n)\big)}(t),\quad0<t<\infty.
$$ 
Hence, $f\/$ does not have an
  absolutely continuous norm.

(ii) If $X=\bigcup_{n=0}^\infty X_n\/$ with $X_n\/$ an atom of measure  $b>0\/$ and  $X_n\cap
X_m=\emptyset,\;n\ne m,\/$ the function
$$
f=\sum_{n=0}^\infty W^{-1/p}\big((n+1)b\big)\chi_{X_n}
$$ 
is in $\Lambda^{p,\infty}_X(w)\/$ 
and has norm  $1\/$. The sets $E_N=\bigcup_{n=N}^\infty
X_n,\;N=1,2,\dots\/$ form a decreasing sequence with  $\chi_{E_n}\to0\/$ a.e. and,
 for each  $N\/$,
$$
\big(f\chi_{E_N}\big)^\ast(nb)=W^{-1/p}\big((n+N+1)b\big),\quad n=0,1,2,\dots
$$ 
Using that $f^\ast\/$ is constant on every interval $\big[nb,(n+1)b\big)\/$ and 
Theorem~\ref{te: II.2.16}, we obtain
\begin{eqnarray*}
\big\|f\chi_{E_N}\big\|_{\Lambda^{p,\infty}}&=&\sup_{t>0}W^{1/p}(t)\big(f\chi_{E_N}\big)^\ast(t)
\\ &=&\sup_{n\ge1}W^{1/p}(nb)W^{-1/p}\big((n+N)b\big)\\
&\ge&\big(W(Nb)W^{-1}(2Nb)\big)^{1/p}\\ &\ge& C^{-1/p}.
\end{eqnarray*}
Then,  
$\lim_N\big\|f\chi_{E_N}\big\|_{\Lambda^{p,\infty}}\neq0\/$ and 
$\Lambda^{p,\infty}_X(w)\/$ does not have absolutely continuous norm. $\quad\qed$
\edem

\bdefi
$$
L^\infty_0(X)\index{$L^\infty_0(X)$}=\big\{f\in L^\infty(X)\,:\,\mu\big(\{f\neq0\}\big)<\infty\big\}.
$$
\edefi

We can now state a positive partial result for the spaces  $\Lambda^{p,\infty}\/$.

\bpro
\label{pro: II.3.9} 
If $0<p<\infty\/$ and $W\in\Delta_2(X)\/$,
every function in  $L^\infty_0(X)\/$ has  absolutely continuous norm in
$\Lambda^{p,\infty}_X(w)\/$.
\epro

\bdem
 If $f\in L^\infty_0\/$, let  $Y=\{f\neq0\}\subset X\/$. 
Then, if $(A_n)_n\/$ is a sequence of measurable sets with  $\chi_{A_n}\to0\/$
a.e., the functions $\big(f\chi_{A_n}\big)_n\/$ are zero in the complement of  $Y\/$ 
and $\lim_n f(x)\chi_{A_n}(x)=0\hbox{\ a.e.\ }x\in Y\/$. Since  $Y\/$ has finite measure, it
follows from    Egorov's theorem,  that the convergence of the previous 
functions is quasi-uniform. Thus,  since  $\epsilon>0\/$, there exists a set
$Y_\epsilon\subset Y\/$ of measure less than
$\epsilon\/$ and such that $f\chi_{A_n}\to0\/$ uniformly in $X\setminus Y_\epsilon\/$. Therefore,
there exists 
$n_0\in\N\/$ such that
$$
\big|f(x)\chi_{A_n}(x)\big|<\epsilon,\quad x\in X\setminus Y_\epsilon,\;n\ge n_0.
$$

Let $n\ge n_0\/$. Then:

(i) if  $t\ge\|f\|_\infty\/$,
$$
tW^{1/p}\big(\lambda_{f\chi_{A_n}}(t)\big)=tW^{1/p}(0)=0,
$$

(ii) if $\epsilon\le t\le\|f\|_\infty\/$,
$$
tW^{1/p}\big(\lambda_{f\chi_{A_n}}(t)\big)\le tW^{1/p}\big(\mu(Y_\epsilon)\big)
\le\|f\|_\infty W^{1/p}(\epsilon),
$$

(iii) if $0\le t<\epsilon\/$,
$$
tW^{1/p}\big(\lambda_{f\chi_{A_n}}(t)\big)\le \epsilon W^{1/p}\big(\mu(Y)\big)
$$ 
and by Proposition~\ref{pro: II.2.5} (iii) we get that
$$
\big\|f\chi_{A_n}\big\|_{\Lambda^{p,\infty}_X(w)}\le
\max\big\{\|f\|_\infty W^{1/p}(\epsilon),\epsilon W^{1/p}\big(\mu(Y)\big)\big\},\qquad n\ge n_0,
$$
that is,
$$
\limsup_n \big\|f\chi_{A_n}\big\|_{\Lambda^{p,\infty}_X(w)}\le 
\max\big\{\|f\|_\infty W^{1/p}(\epsilon),\epsilon W^{1/p}\big(\mu(Y)\big)\big\}.
$$ 
Since this inequality holds for every $\epsilon>0\/$ and  $\lim_{t\to0}W(t)=0\/$, we conclude that 
$$
\lim_n\big\|f\chi_{A_n}\big\|_{\Lambda^{p,\infty}_X(w)}=0.\qquad\qed
$$
\edem

\subsection{Density  of the  simple functions and    $L^\infty\/$}

Let us now study the density  of the  simple functions and also of  $L^\infty\/$ 
 in the quasi-normed  Lorentz spaces. This question is solved in a positive way for the spaces
$\Lambda^p\/$ (except when
$\mu(X)=\infty,\;w\in L^1\/$).  The behavior of the density in the spaces  
$\Lambda^{p,\infty}\/$ is more irregular. It is not true here that the simple functions are dense and,
in fact,  not even $L^\infty\/$ is always dense in  $\Lambda^{p,\infty}\/$. 

\bdefi
Let us denote by  $\S=\S(X)\/$\index{$\S(X)$}  the class of simple functions in $X\/$. That is
$$
\S=\big\{f\in\M(X)\,:\,\hbox{card}\big(f(X)\big)<\infty\big\}.
$$
$\SS\/$\index{$\SS(X)$} will be the simple functions  with support in a set of finite measure:
$$
\SS=\SS(X)=\big\{f\in\S\,:\,\mu\big(\{f\neq0\}\big)<\infty\big\}.
$$
\edefi

It is clear that $\SS\subset\LL\subset\Lambda\/$ for every  Lorentz space $\Lambda\/$.

\blem
\label{lem: II.3.11}
 If $f\in\M(X)\/$ there exists a sequence $(s_n)_n\subset\S\/$ satisfying
\begin{enumerate}
\item[(i)] $\lim_n s_n(x)=f(x),\;x\in X$,

\item[(ii)] $\big(|s_n(x)|\big)_n$ is an increasing sequence and $|s_n(x)|\le|f(x)|$, for every $x\in X$,

\item[(iii)] $|f-s_n|\le|f|,\;n\in\N.$
\end{enumerate}
Moreover, if $f\/$ is bounded then $s_n\to f\/$ uniformly.
 In particular $\S\/$ is dense in  $L^\infty\/$.
\elem

\bteo
\label{te: II.3.12}
If $\Lambda(w)\/$ is a quasi-normed Lorentz space,
$$
\SS\subset\LL\subset\Lambda(w)
$$ 
and  $\SS\/$ is dense in $\LL\/$ with the topology of $\Lambda\/$. If $w\in L^1\/$,
$$
\S\subset L^\infty\subset\Lambda(w)
$$ 
and $\S\/$ is dense in $L^\infty\/$ but, 
in this case, $\LL\/$ is not dense in  $L^\infty\/$, if $\mu(X)=\infty\/$ (always with the topology
of $\Lambda\/$).
\eteo

\bdem
The embeddings are trivial. On the other hand, if $f\in\LL\/$ 
and $E=\{f\neq0\},\;\mu(E)<\infty\/$ and if $(s_n)_n\subset\S\/$ is the sequence of the previous
lemma, since $|s_n|\le|f|\/$, these functions have also support in  $E\/$ 
and  $(s_n)_n\subset\SS\/$. Besides $\|f-s_n\|_{L^\infty}\to0\/$ and thus,
$$
\|f-s_n\|_{\Lambda}\le\big\|\|f-s_n\|_{L^\infty}\chi_E\big\|_{\Lambda}
=\|f-s_n\|_{L^\infty}\|\chi_E\|_{\Lambda}\to0.
$$ 
This proves that $\SS\/$ is dense in 
$\LL\/$. If $w\in L^1,\;L^\infty\subset\Lambda\/$ and this embedding is continuous
(Theorem~\ref{te: II.3.2}). Since $\S\/$ is dense in $\big(L^\infty,\|\cdot\|_{L^\infty}\big)\/$
(Lemma~\ref{lem: II.3.11}), we have that $\S\/$ is dense in $L^\infty\/$ with the topology of 
$\Lambda\/$.  In this case $1\in L^\infty\subset\Lambda\/$ and if  $g\in\LL,\;|g-1|=1\/$ in
a set of infinite measure if  $\mu(X)=\infty\/$. Thus,  $(1-g)^\ast\ge 1\/$ and  $\|1-g\|_{\Lambda}
\ge C_{\Lambda}>0\/$. Therefore,  $1\/$ cannot be a limit in 
$\Lambda\/$ of functions in  $\LL\/$ and this set is not dense in $L^\infty.\quad\qed$
\edem

\bteo\label{te: II.3.13}
 Let  $0<p<\infty\/$ and $W\in\Delta_2(X)\/$. Then,
\begin{enumerate}
\item[(i)] if  $w\notin L^1,\;\LL\/$ is dense in $\Lambda^p(w)\/$, 

\item[(ii)] if  $w\in L^1,\;L^\infty\/$ is dense in $\Lambda^p(w)\/$. But, in this case, 
$\LL(X)\/$ is not dense in $\Lambda^p_X(w)\/$, if $\mu(X)=\infty\/$.
\end{enumerate}
\eteo

\bdem
Let us see that $L^\infty\/$ ($\LL\/$ in the case $w\notin L^1\/$) 
is dense in $\Lambda^p\/$. If $f\in\Lambda^p\/$,  Proposition~\ref{pro: II.2.5} tells us that
$\lim_{t\to\infty}\lambda_f(t)=0\/$. If for each  $n=1,2,\dots\/$ we define $f_n
=f\chi_{\{|f|\le n\}}\/$, since $f-f_n=
f\chi_{\{|f|>n\}},$ then $(f-f_n)^\ast=f^\ast\chi_{\big[0,\lambda_f(n)\big)}\/$  and we have that 
$\lim_n(f-f_n)^\ast(t)=0,\;t>0\/$. On the other hand $(f-f_n)^\ast\le f^\ast\/$ and by the dominated
convergence theorem, we get 
$\|f-f_n\|_{\Lambda^p}\to0\/$.  This proves that 
$L^\infty\cap\Lambda^p\/$  is dense in  $\Lambda^p\/$. If $w\in L^1\/$,  the first of these spaces
is 
$L^\infty\/$ and we are done. If $w\notin L^1\/$
 we only have to see that every function in  $L^\infty\cap\Lambda^p\/$ is limit (in
$\Lambda^p\/$) of functions in  $\LL\/$. But if $g\in L^\infty\cap\Lambda^p\/$ 
and  $w\notin L^1,\;\lim_{t\to\infty}g^\ast(t)=0\/$ and the functions
$g_n=g\chi_{\{|g|>g^\ast(n)\}}\in\LL,\;n=1,2,\dots,\/$ satisfy $(g-g_n)^\ast\le 
g^\ast(n)\to0\/$ and since  $(g-g_n)^\ast\le g^\ast\/$, the dominated
convergence theorem shows that $\|g-g_n\|_{\Lambda^p}\to0\/$.

It only remains to see that $\LL\/$ is not dense in $\Lambda^p(w)\/$ if $w\in L^1\/$ and 
$\mu(X)=\infty\/$. But, in this case, we know by  Theorem~\ref{te: II.3.12}, that $\LL\/$
is not dense in  $L^\infty\/$ and, hence, it can neither be dense in $\Lambda^p\supset
L^\infty.\quad\qed$
\edem

In the weak-type case we have the following result of density.

\bteo
If $X\/$ is a resonant measure space, $0<p<\infty\/$ and $\Lambda^{p,\infty}_X(w)\/$ 
is a quasi-normed space,
\begin{enumerate}
\item[(i)] $\LL\/$ is dense in  $\Lambda^{p,\infty}\cap L^\infty\/$  if and only if  $\mu(X)<\infty\/$ 
and then they coincide (the density is considered with respect to the topology of 
$\Lambda^{p,\infty}\/$).

\item[(ii)] $\Lambda^{p,\infty}\cap L^\infty\/$ is dense in $\Lambda^{p,\infty}\/$ 
if and only if  $X\/$ is atomic and then they coincide.
\end{enumerate}
\noindent Hence, in both cases, if these two spaces do not coincide, the smaller one is not dense in
the other.
\eteo

\bdem
 (i) If $\mu(X)<\infty\/$ both spaces coincide and there is nothing to prove. If $\mu(X)=\infty\/$
and 
$X\/$ is not atomic, we can take   $f\in L^\infty\cap\Lambda^{p,\infty}\/$ with
$$
f^\ast=W^{-1/p}(1)\chi_{[0,1)}+W^{-1/p}\chi_{[1,\infty)}.
$$ 
If $g\in\LL\/$ and  $E=\{g=0\}\/$ 
then $b=\mu(X\setminus E)<\infty\/$ and 
$$
(f-g)^\ast(t)\ge\big((f-g)\chi_E\big)^\ast(t)=\big(f\chi_E\big)^\ast(t)\ge f^\ast(b+t)
=W^{-1/p}(t+b),\quad t>1,
$$ 
and we have that,
$$
\|f-g\|_{\Lambda^{p,\infty}}=\sup_{t}W^{1/p}(t)(f-s)^\ast(t)\ge \sup_{t>1}{W^{1/p}(t)
\over W^{1/p}(t+b)}.
$$ 
If $w\in L^1\/$ the above supremum is  $1\/$ and, in the opposite case, since
$W\in\Delta_2,\;W(t+b)\le C\big(W(t)+W(b)\big)\/$ (cf. Theorem~\ref{te: II.2.13} and
Proposition~\ref{pro: II.2.12})  we obtain that
$\|f-g\|_{\Lambda^{p,\infty}}\ge C^{-1/p}\/$ and  $f\/$ cannot be the limit of functions in $\LL\/$.
This shows that $\LL\/$ is not dense in  $L^\infty
\cap\Lambda^{p,\infty}\/$ in the nonatomic case. If $\mu(X)=\infty\/$ and  $X\/$ is an atomic
resonant measure space,
$X=\bigcup_{n=0}^\infty X_n\/$ with every $X_n\/$ an atom of measure  $b>0\/$ and  $X_n\cap
X_m=\emptyset,\;n\ne m,\/$ and the function
$$
f=\sum_{n=0}^\infty W^{-1/p}\big((n+1)b\big)\chi_{X_n}
$$ 
is in  $\Lambda_X^{p,\infty}(w)
\cap L^\infty\/$ and   has norm $1\/$. If $g\in\LL\/$, then its support has finite measure (a finite
union 
$X_{n_1}\cup\dots\cup X_{n_k}\/$ of atoms) and there exists $N\in\N\/$
such that $g=0\/$ in $\bigcup_{n=N}^\infty X_n\/$. Hence,
$(f-g)^\ast(nb)\ge f^\ast\big((n+N)b\big)=W^{-1/p}\big((n+N+1)b\big),\;n=0,1,2,\dots,\/$ 
proceeding as in the proof of  Theorem~\ref{te: II.3.7}  we can conclude that 
$\|f-g\|_{\Lambda^{p,\infty}}\ge C>0\/$. Thus,  $f\/$ is not the limit of functions in $\LL\/$ 
and consequently this set is not dense in 
$\Lambda^{p,\infty}\cap L^\infty\/$ and (i) is proved.

(ii) As before, if $X\/$ is atomic both spaces coincide and there is nothing to be proved.  If
$X\/$ is nonatomic, there exists $0\le f\in\M(X)\/$ with $f^\ast=W^{-1/p}\/$ in $[0,\mu(X))\/$.
Thus, 
$\|f\|_{\Lambda^{p,\infty}}=1\/$ and $f\in\Lambda^{p,\infty}\setminus
 L^\infty\/$. If $g\in L^\infty(X)\/$ with $\|g\|_\infty=b,$ then $|f-g|\ge|f|-b$
and, hence, $(f-g)^\ast(t)\ge f^\ast(t)-b=W^{-1/p}(t)-b\/$ if $0\le t<\lambda_f(b)=t_0\/$.
Therefore,
$$
\|f-g\|_{\Lambda^{p,\infty}}=\sup_t W^{1/p}(t)(f-g)^\ast(t)\ge\sup_{0<t<t_0}\big(1-b
W^{1/p}(t)\big)=1,
$$ 
since $t_0=\lambda_f(b)>0\/$. We conclude that $f\/$ is not the limit in $\Lambda^{p,\infty}\/$ of
functions in 
$L^\infty(X)\/$. This proves that 
$L^\infty\cap\Lambda^{p,\infty}\/$  is not dense in  $\Lambda^{p,\infty}\/$  in this case.
$\quad\qed$
\edem

\section{Duality}
\label{se:dual}

In this section $(X,\mu)\/$ will be a $\sigma$-finite measure space.  
$d(\Omega,p)\/$  and its weak version   $d^\infty(\Omega,p)\/$ will be the 
Lorentz spaces on  $\N^\ast\/$ introduced in  Example~\ref{Eje2.3}(iv).
Here $\Omega=(\Omega_n)_{n=0}^\infty\/$ is a sequence of positive numbers $\Omega_0\neq0\/$. 

We shall study the dual spaces and the associate spaces of the 
Lorentz spaces on $X\/$ with special attention to the case when $X\/$ is a resonant measure space
(Definition~\ref{def: II.2.15}).
We shall describe the associate space and we shall deduce a necessary  and sufficient condition to
have that the dual and the associate spaces  coincide. In some cases, we shall also give a
description of the Banach envelope  of 
$\Lambda\/$, which is of interest if this space is not a normed space. In our study, we include
 the sequence Lorentz spaces $d(\Omega,p)\/$ with  $\Omega\/$ arbitrary. As far as we know, only the
Banach case has been previously studied  (that is, the case $\Omega\/$ decreasing).

We shall introduce the associate space generalizing the definition that   can be found in 
\cite{BS} in the context of  Banach function spaces.

\bdefi
\label{def:II.4.1}
 If $\|\cdot\|:\M(X)\to[0,\infty]\/$ is a  positively homogeneous functional 
and  $E=\big\{f\in\M(X)\,:\,\|f\|<\infty\big\}\/$, we define the
 associate \  \lq\lq norm\rq\rq
$$
\|f\|_{E^\prime}=\sup\bigg\{\int_X|f(x)g(x)|\,d\mu(x)\,:\,\|g\|\le1,\;g\in E\bigg\},\quad f\in\M(X).
$$ 
The associate space\index{associate space} of  $E\/$ is then 
$E^\prime=\big\{f\in\M(X)\,:\,\|f\|_{E^\prime}<\infty\big\}\/$.
\edefi

\bobs
\label{Obs: II.4.2} 
The following properties are immediate:

(i)  $\|\cdot\|_{E^\prime}\/$ is subadditive and positively homogeneous and, if $E\/$ 
contains the characteristic functions of sets of finite measure,
$\big(E^\prime,\|\cdot\|_{E^\prime}\big)\/$ is a normed space.

(ii) If $(E,\|\cdot\|)\/$ is a lattice ($\|f\|\le\|g\|\/$, if $|f|\le|g|\/$) 
then  $E^\prime\/$ is also a lattice. If $\|\cdot\|\/$ has the  Fatou property, the same happens to 
$\|\cdot\|_{E^\prime}\/$.

(iii) If we denote by  $E^{\prime\prime}=(E^\prime)^\prime\/$ we have that
$E\subset E^{\prime\prime}\/$ and   $\|f\|_{E^{\prime\prime}}\le\|f\|\/$ for every  $f\in\M(X)\/$.
\eobs

\bdefi
\label{def4.3}
A  Banach function space (BFS)\index{Banach function space} $E\/$ in $X\/$  is a subspace of   $\M(X)\/$
defined by 
$E=\{f:\|f\|<\infty\}\/$ where
$\|\cdot\|=\|\cdot\|_E\/$ is a norm (called \lq\lq Banach function norm") 
that satisfies the following properties, for $f,g,f_n\in E,\;A\subset X\/$ measurable (see \cite{BS}):
\begin{enumerate}
\item[(i)]  $\|f\|\le\|g\|\/$ if $|f|\le|g|\/$,

\item[(ii)]  $0\le f_n\le f_{n+1}\to f\Rightarrow\|f_n\|\to\|f\|\/$,

\item[(iii)]  $\chi_A\in E\/$ if $\mu(A)<\infty\/$ and $E\neq\{0\}\/$,

\item[(iv)]  $\displaystyle\int_A|f(x)|\,d\mu(x)\le C_A\|f\|\/$ if $\mu(A)<\infty\/$.
\end{enumerate}
\edefi

The result that   follows establishes that the associate space  $\Lambda^\prime_X\/$ 
of a quasi-normed Lorentz space is a   Banach function space whenever  $X\/$ is a resonant measure
space.

\bteo
If $\Lambda_X\/$ is a quasi-normed  Lorentz space on  a resonant measure space $X\/$, the associate
space 
$\Lambda^\prime_X\/$  is a  rearrangement invariant Banach function space\index{rearrangement invariant}.
For every 
$f\in\Lambda^\prime_X\/$,
\bequ\label{asso}
\|f\|_{\Lambda^\prime_X}=\sup\bigg\{\int_0^\infty f^\ast(t) g^\ast(t)\,dt\,:\,\|g\|_{\Lambda_X}\le1\bigg\}.
\eequ
Moreover, a function $f\in\M(X)\/$ is in $\Lambda^\prime_X\/$  if and only if
$\displaystyle\int_X|f(x)g(x)|\,d\mu(x)<\infty\/$ for every 
$g\in\Lambda_X\/$ and  $\Lambda^\prime_X\ne\{0\}\/$  if and only if  $\Lambda_X\subset
L^1_{\hbox{\decp loc}}(X)\/$.
\eteo

\bdem
Let  us first prove  (\ref{asso}). If $f,g\in\M(X)$, we have
that $\displaystyle\int_X|f(x)g(x)|\,d\mu(x)\le\displaystyle\int_0^\infty  f^\ast(t) g^\ast(t)\,dt\/$ and
thus,
$$
\|f\|_{\Lambda^\prime_X}\le\sup\bigg\{\int_0^\infty f^\ast (t)
g^\ast(t)\,dt\,:\,\|g\|_{\Lambda_X}\le1\bigg\}.
$$
To prove the other inequality, we observe that  for
$g\in\Lambda_X\/$, with  $\|g\|_{\Lambda_X}\le1\/$ we have, since $X\/$ is resonant,
\begin{eqnarray*} 
\int_0^\infty f^\ast(t) g^\ast(t)\,dt&=&\sup_{h^\ast=g^\ast}\int_X|f(x)h(x)|\,d\mu(x)\\
&\le&\sup_{\|h\|_{\Lambda_X}\le1}\int_X|f(x)h(x)|\,d\mu(x)
=\|f\|_{\Lambda^\prime_X}.
\end{eqnarray*} 
Hence, the functional 
$\|\cdot\|_{\Lambda^\prime_X}\/$ is rearrangement invariant. 

Let us now prove that  $f\in\Lambda^\prime_X\/$ if and only if 
$\displaystyle\int_X|f(x)g(x)|\,d\mu(x)<\infty\/$  for every $g\in\Lambda_X\/$. The necessity is a
consequence of the definition of 
$\Lambda^\prime_X\/$. The sufficiency follows from the closed graph theorem (see
\cite{Ru2})  since under the hypothesis, the linear operator $T_f(g)=fg,\;T_f:\Lambda_X\to
L^1(X)\/$ is well defined and  both are  F-spaces continuously embedded in  $\M(X)\/$ 
(from which one can immediately see that the graph is closed). $T_f\/$ is then continuous and this
proves    (Theorem~\ref{te: II.3.2}
(v)) that $\displaystyle\int_X|f(x)g(x)|\,d\mu(x)\le C\|g\|_{\Lambda_X}\/$  and thus 
$f\in\Lambda^\prime_X\/$.

Let us now see that  $\Lambda^\prime_X\/$ is a   Banach function space. That 
$\|\cdot\|_{\Lambda^\prime_X}\/$  is a norm is immediate. Properties 
(i) and  (ii) in the definition of BFS (Definition~\ref{def4.3}) are also trivial.  Property  (iv)
of that definition is also very easy, since if  $\mu(A)<\infty\/$ and $f\in\Lambda^\prime_X\/$,
$$
\int_A|f(x)|\,d\mu(x)\le\|\chi_A\|_{\Lambda_X}\|f\|_{\Lambda^\prime_X}.
$$ 
It only remains to prove (iii). That is that 
$\chi_A\in\Lambda^\prime_X\/$ if $\mu(A)<\infty\/$ and
$\Lambda^\prime_X\neq\{0\}\/$. In the atomic case  it is trivial and if $X\/$  is not atomic we have
that, if
$0\neq f\in\Lambda^\prime_X\/$ there exists $t>0\/$ such that, if 
$E=\{|f|>t\},$ then $\mu(E)>0\/$. Hence, $t\chi_E\le|f|\/$ and   $\chi_E\in\Lambda^\prime_X\/$. That is, if 
$\Lambda^\prime_X\neq\{0\}\/$ there exists 
$\chi_E\in\Lambda^\prime_X\/$ with $b=\mu(E)>0\/$. Since $\|\cdot\|_{\Lambda^\prime_X}\/$ is
rearrangement invariant and monotone, we have that  
$\chi_A\in\Lambda^\prime\/$ for every measurable set  $A\/$ with $\mu(A)\le b\/$. If
$\infty>\mu(A)>b\/$,  the above is also true since $A=\bigcup_{n=1}^N A_n\/$ with
$\mu(A_n)\le b\/$ and $\|\chi_A\|_{\Lambda^\prime_X}\le\sum_n\|\chi_{A_n}\|_{\Lambda^\prime_X}<\infty\/$.

It only remains to prove that  $\Lambda^\prime_X\ne\{0\}\/$  if and only if $\Lambda_X\subset
L^1_{\hbox{\decp loc}}(X)\/$. The necessity is immediate since 
$\Lambda^\prime_X\ne\{0\}\/$ and  $\mu(A)<\infty\/$ then we have already seen
that $\chi_A\in\Lambda^\prime_X\/$ and, in particular
$\displaystyle\int_A|f(x)|\,d\mu(x)<\infty,\;f\in\Lambda_X\/$. On the other hand if
$\Lambda_X\subset L^1_{\hbox{\decp  loc}}(X)\/$ and  $0<\mu(A)<\infty\/$ we have that
$\displaystyle\int_X|f(x)\chi_A(x)|\,d\mu(x)<\infty\/$ for every  $f\in\Lambda_X\/$ and that implies
$\chi_A\in\Lambda^\prime_X\/$ and $\Lambda^\prime_X\neq\{0\}.\quad\qed$
\edem

\bobs
 If in the previous theorem, we omit the condition that   $X\/$
is a resonant measure space   we  can still prove that 
$\Lambda^\prime_X\/$ is a  Banach space  with monotone norm and with the Fatou property. 
The convergence in $\Lambda^\prime_X\/$ implies the convergence in measure on sets of finite measure
 (same for the  Cauchy sequences) and it is also true that  $f\in\Lambda^\prime_X\/$ if and only if 
$\displaystyle\int_X|f(x)g(x)|\,d\mu(x)<\infty\/$ for every
$g\in\Lambda_X\/$. The last statement of the theorem is not true  (in general) in this case. 
\eobs

The result that follows describes the associate space of  $\Lambda_X\/$  in the case $X\/$
nonatomic. But first, we need to define the Lorentz space $\Gamma\/$.\index{Lorentz space!Gamma}

\bdefi  
\label{def: II.4.6}
If $0<p\le\infty\/$ we define
$$
\Gamma^p_X(w)\index{$\Gamma^p_X(w)$}=\bigg\{f\in\M(X)\,:\,\|f\|_{\Gamma^p_X(w)}=\bigg(\int_0^\infty
(f^{\ast\ast}(t))^pw(t)\,dt\bigg)^{1/p}<\infty\,\bigg\}.
$$ 
The weak-type version of the previous space is
$$
\Gamma^{p,\infty}_X(w)\index{$\Gamma^{p,\infty}_X(w)$}=\bigg\{f\in\M(X)\,:\,\|f\|_{\Gamma^{p,\infty}_X(w)}=\sup_{t>0}
W^{1/p}(t)f^{\ast\ast}(t)<\infty\,\bigg\}.
$$
The last definition can be extended in the following way. If $\Phi\/$  is an arbitrary function in 
$\R^+\/$ we write
$$
\Gamma^{p,\infty}_X(\di\Phi)\index{$\Gamma^{p,\infty}_X(\di\Phi)$}=\bigg\{f\in\M(X)\,:\,\|f\|_{\Gamma^{p,\infty}_X(\dip\Phi)}
=\sup_{t>0}\Phi^{1/p}(t)f^{\ast\ast}(t)<\infty\,\bigg\}.
$$
\edefi

We can always assume that the function $\Phi\/$ in the previous definition is increasing since
otherwise it can be  substituted by 
$\widetilde{\Phi}(t)=\sup_{0<s<t}\Phi(s),\;t>0,\/$ which is increasing and satisfies 
$\|f\|_{\Gamma^{p,\infty}_X(\dip\Phi)}=\|f\|_{\Gamma^{p,\infty}_X(\dip\widetilde{\Phi})}\/$.

 Condition (ii) of the following result  is a direct consequence of   Sawyer's formula 
stated in  Theorem~\ref{te: I.5.7} while  (i) was proved in 
\cite{CS2}.

\bteo
\label{te: II.4.7}
 Let $X\/$  be a nonatomic measure space and let  $w\/$ be an arbitrary weight in  $\R^+\/$.
\begin{enumerate}
\item[(i)] If $0<p\le1\/$, then
$$
\Lambda_X^p(w)^\prime=\Gamma^{1,\infty}_X(\di\Phi)\qquad\hbox {(with equal norms),}
$$
where
$\Phi(t)=tW^{-1/p}(t),\;t>0\/$.

\item[(ii)] If $1<p<\infty\/$ and $f\in\M(X)\/$, then
\begin{eqnarray*}
\|f\|_{\Lambda_X^p(w)^\prime}&\approx&\bigg(\int_0^\infty\bigg({1\over W(t)}\int_0^t
f^\ast(s)\,ds\bigg)^{p^\prime}w(t)\,dt\bigg)^{1/p^\prime}+{\displaystyle\int_0^\infty f^\ast(t)\,dt\over
W^{1/p}(\infty)}\\ &\approx&\bigg(\int_0^\infty\bigg({1\over W(t)}\int_0^t
f^\ast(s)\,ds\bigg)^{p^\prime-1}f^\ast(t)\,dt\bigg)^{1/p^\prime}.
\end{eqnarray*}

\item[(iii)] If $0<p<\infty\/$, then 
$$
\Lambda_X^{p,\infty}(w)^\prime=\Lambda^1(W^{-1/p})\qquad\hbox{(with equal norms).}
$$
\end{enumerate}
\eteo

\bdem
 Since  $X\/$ is nonatomic and  $\sigma$-finite, every decreasing function in 
$[0,\mu(X))\/$ equals  a.e. to the decreasing rearrangement of a function in 
$\M(X)\/$. Besides,   $X\/$ is resonant  and hence,
$$
\|f\|_{\lox{p,q}{X}{w}^\prime}=\sup_{g\in\lox{p,q}{X}{w}}{\displaystyle\int_X|f(x)g(x)|\,d\mu(x)\over\|g\|_{\lox{p,q}{X}{w}}
}=\sup_{g\downarrow}{\displaystyle\int_0^\infty
g(t)f^\ast(t)\,dt\over\|g\|_{L^{p,q}(w)}}.
$$ 
The first case (i) is solved applying Corollary~\ref{cor: I.2.14}
(with $p_1=1,\,T=\hbox{\rm Id}\/$) to the regular class  $L=\Fd\/$  or  
Theorem 2.12 in \cite{CS1}. The second case  (ii) is an immediate consequence of Theorem~\ref{te: I.5.7}
(E. Sawyer). Finally, (iii) corresponds to  $q=\infty\/$ and we have that,
$$
\|f\|_{\lox{p,\infty}{X}{w}^\prime}=\sup\bigg\{\int_0^\infty f^\ast(t)
g(t)\,dt\,:\,\|g\|_{L^{p,\infty}}(w)=1,\;g\downarrow\bigg\}.
$$ 
Now, if 
$g\downarrow,\;\|g\|_{L^{p,\infty}(w)}=\sup_t W^{1/p}(t)g(t)\/$ and  $\|g\|_{L^{p,\infty}}(w)=1\/$
implies $g\le W^{-1/p}\/$, and therefore
$$
\|f\|_{\lox{p,\infty}{X}{w}^\prime}\le\int_0^\infty f^\ast(t) W^{-1/p}(t)\,dt=\|f\|_{\lox{1}{}{W^{-1/p}}}.
$$
On the other hand $W^{-1/p}\/$  is decreasing, 
$\big\|W^{-1/p}\big\|_{L^{p,\infty}}(w)=1\/$ and we have the equality. $\quad\qed$
\edem

\bobs
\label{obs: II.4.8}

(i) If $p>1\/$ and  $\tilde{w}(t)=t^{p^\prime}W^{-p^\prime}(t)w(t),\;t>0\/$, then
(ii) of the previous theorem can be stated of the following way:
\begin{eqnarray*}
\lox{p}{X}{w}^\prime&=&\Gamma^{p^\prime}_X(\tilde{w}),\qquad\hbox{\rm if\ }w\notin L^1,\\
\lox{p}{X}{w}^\prime&=&\Gamma^{p^\prime}_X(\tilde{w})\cap L^1(X),\qquad\hbox{\rm if\ }w\in L^1.
\end{eqnarray*}
It is assumed that the norm  in the intersection space is the maximum of the sum  of both norms  and
in these equalities we are assuming the equivalence of the norms.

(ii) If $w\notin L^1,\;p>1\/$ and $\tilde{w}\/$ (as above defined) is in $B_{p^\prime}\/$,  the  Hardy
operator $A\/$ satisfies the boundedness $A:L^{p^\prime}_{\hbox{\decp  dec}}(\tilde{w})\to
L^{p^\prime}(\tilde{w})\/$ and then,
$$
\|f\|_{\Lambda^p_X(w)^\prime}\approx\|f\|_{\Lambda_X^{p^\prime}(\tilde{w})}\,,\quad f\in\M(X).
$$
It is easy to see that this condition on the weight    $w\/$ is equivalent to 
\bequ
\label{II.4.9}
\bigg(\int_0^r\bigg({W(t)\over t}\bigg)^{-p^\prime}w(t)\,dt\bigg)^{1/p^\prime}W^{1/p}(r)\ge
Cr,\quad r>0, 
\eequ
which is the opposite inequality to the condition   $w\in
B_{p,\infty}=B_p\/$ (Theorem~\ref{te: I.6.5}). Since one of the embedding always holds, it follows
that condition  (\ref{II.4.9})  is necessary  and sufficient  (in the case
$w\notin L^1\/$) to have the identity 
$$
\lox{p}{X}{w}'=\Lambda^{p^\prime}_X(\tilde{w})\qquad\hbox{\rm (with equivalent norms).}
$$

(iii) If $1<p<\infty\/$ the space $\Lambda_X^{p^\prime}(w^{1-p^\prime})\/$ is embedded in 
$\lox{p}{X}{w}^\prime\/$ since, by  H\"older's inequality,
\begin{eqnarray*} 
\int_X |f(x)g(x)|\,d\mu(x)&\le&\int_0^\infty (f^\ast(t) w^{1/p}(t))(g^\ast(t)
w^{-1/p}(t))\,dt\\
&\le&\|f\|_{\Lambda^p(w)}\|g\|_{\Lambda^{p^\prime}(w^{1-p^\prime})}.
\end{eqnarray*} 
If $\tilde{w}\/$ is as before,  we have that 
\bequ
\label{II.4.10}
\lox{p^\prime}{X}{w^{1-p^\prime}}\subset\lox{p}{X}{w}^\prime\subset
\Gamma^{p\prime}_X(\tilde{w})\subset\Lambda^{p^\prime}_X(\tilde{w}). 
\eequ

(iv)  The  comments made after    Definition~\ref{def: II.4.6}
  tell  us that, for every 
$0<p\le1\/$,
$$
\Lambda_X^p(w)^\prime=\Gamma^{1,\infty}_X(\di\Phi_p)\qquad\hbox{\rm (with equal  norms),}
$$
 where
$\Phi_p(t)=\sup_{0<s<t}sW^{-1/p}(s),\;t>0\/$. Note that
$$
\Phi_p(t)={t\over\inf_{0<s<t}{t\over s}W^{1/p}(s)}={t\over\inf_{s>0}\max\big(1,{t\over
s}\big)W^{1/p}(s)}={t\over W_p(t)},
$$ 
where
$W_p(t)=\inf_{s>0}\max\big(1,{t\over s}\big)W^{1/p}(s),\;t>0\/$. It is easy to check that this
function is quasi-concave (\cite{BS}). In fact, it is the biggest quasi-concave function majorized
by  
$W^{1/p}\/$ and consequently it is called the greatest concave minorant of  $W^{1/p}\/$ (see
\cite{CPSS}).  Hence, we can write,
$$
\|f\|_{\Lambda^p_X(w)^\prime}=\sup_{t>0}{1\over W_p(t)}\int_0^t f^\ast(s)\,ds,\qquad f\in\M(X).
$$
\eobs

As a first consequence of Theorem~\ref{te: II.4.7} we obtain the  characterization of the weights 
$w\/$  for which $\Lambda^\prime(w)=\{0\}\/$.

\bteo
\label{te: II.4.11} 
If $X\/$ is nonatomic:
\begin{enumerate} 
\item[(i)] If $0<p\le1,\;(\Lambda^p_X(w))^\prime\ne\{0\}\Leftrightarrow\sup_{0<t<1}{\displaystyle t^p\over 
\displaystyle W(t)}<\infty.$
\item[(ii)] If $1<p<\infty,\;(\Lambda^p_X(w))^\prime\ne\{0\}\Leftrightarrow\displaystyle
\int_0^1\bigg({\displaystyle t\over \displaystyle W(t)}\bigg)^{p^\prime-1}\,dt<\infty.$
\item[(iii)] If $0<p<\infty,\;(\Lambda^{p,\infty}_X(w))^\prime\ne\{0\}\Leftrightarrow\displaystyle \int_0^1
{\displaystyle 1\over \displaystyle W^{1/p}(t)}\,dt<\infty.$
\end{enumerate} 
\eteo

\bdem
 Since $\Lambda^\prime\/$ is a BFS, it is not identically zero if and only if 
it contains the functions $\chi_E\/$ with $\mu(E)<\infty\/$, that is  if and only if
$\|\chi_E\|_{\Lambda^\prime}<\infty,\;\mu(E)<\infty\/$. The conditions are then obtained 
applying Theorem~\ref{te: II.4.7}.~$\quad\qed$
\edem

\bobs
Since $W\/$ is continuous in $\R^+$,  the condition of the previous theorem only depends on the
local behavior of  $w\/$ at  $0\/$.
\eobs

\bdefi
A weight $w\/$ in $\R^+\/$ is called  regular\index{regular weight} (see \cite{Re}) if it satisfies 
$$
{W(t)\over t}\le C\,w(t),\qquad t>0,
$$ 
with $C>0\/$ independent of $t\/$. A sequence of positive numbers  $(\Omega_n)_{n=0}^\infty\/$ is
said, analogously, to be regular if\index{regular sequence}
$$
{1\over n+1}\sum_{k=0}^n\Omega_k\le C\,\Omega_n,\qquad n=0,1,2,\dots
$$
\edefi

Every increasing sequence   and every power sequence    are regular. In
\cite{Re} it is proved that if  $w\/$ is regular and decreasing,  the spaces
$\lox{p^\prime}{X}{w^{1-p^\prime}}\/$ and $\lox{p}{X}{w}^\prime\/$ coincide. In the following
result, we extend this to the case of an arbitrary  weight $w\/$.

\bteo
\label{te: II.4.14} 
Let  $1<p<\infty\/$ and  $X\/$ be nonatomic. Then: 
\begin{enumerate}
\item[(i)] If $w\notin L^1,\;\lox{p}{X}{w}^\prime=\lox{p^\prime}{X}{w^{1-p^\prime}}\/$
if and only if there exists  $C>0\/$ such that, for $r>0$,
\bequ\label{equiv}
\int_0^r w^{1-p^\prime}(t)\,dt\le C\bigg(r^{p^\prime}W^{1-p^\prime}(r)+\int_0^r
t^{p^\prime}W^{-p^\prime}(t)w(t)\,dt\bigg).
\eequ

\item[(ii)] If $w\/$ is regular,
$$
\lox{p^\prime}{X}{w^{1-p^\prime}}=\lox{p}{X}{w}^\prime=\Gamma^{p\prime}_X(\tilde{w})
=\Lambda^{p^\prime}_X(\tilde{w}),
$$
where
$\tilde{w}(t)=t^{p^\prime}W^{-p^\prime}(t)w(t),\;t>0\/$.

\item[(iii)] If  $w\/$ is increasing, $\lox{p}{X}{w}^\prime=\lox{p^\prime}{X}{w^{1-p^\prime}}\/$ with
equality of norms.
\end{enumerate}
\eteo

\bdem
 (i) We have observed in  Remark~\ref{obs: II.4.8} that we always have  
$\lox{p^\prime}{X}{w^{1-p^\prime}}\subset\lox{p}{X}{w}^\prime=\Gamma^{p^\prime}_X(\tilde{w})\/$.
Thus, the equality of these two spaces is equivalent to the embedding 
$$
\Gamma^{p^\prime}_X(\tilde{w})\subset\lox{p^\prime}{X}{w^{1-p^\prime}},
$$ 
which   is also equivalent to 
the opposite inequality for the   Hardy operator,
$$
\|g\|_{L^{p^\prime}(w^{1-p^\prime})}\le C\|Ag\|_{L^{p^\prime}(\tilde{w})},\qquad g\d.
$$ 
A necessary  and sufficient condition for it  can be found in \cite{CPSS}  or in \cite{Ne1}. This
condition is 
$$
\int_0^r w^{1-p^\prime}(t)\,dt\le C\bigg(\widetilde{W}(r)+r^{p^\prime}\int_r^\infty
t^{-p^\prime}\tilde{w}(t)\,dt\bigg),
$$ 
which is (\ref{equiv}), since 
$\displaystyle\widetilde{W}(r)=\int_0^r t^{p^\prime}W^{-p^\prime}(t)w(t)\,dt\/$,     and 
$$
r^{p^\prime}\int_r^\infty t^{-p^\prime}\tilde{w}(t)\,dt=r^{p^\prime}\int_r^\infty
W^{-p^\prime}(t)w(t)\,dt={1\over p^\prime-1}r^{p^\prime}W^{1-p^\prime}(r).
$$

(ii) If  $w\/$ is regular we have that
$$
w^{1-p^\prime}(t)=w(t)w^{-p^\prime}(t)\le Cw(t)t^{p^\prime}W^{-p^\prime}(t)=C\tilde{w}(t).
$$ 
Hence, 
$\Lambda^{p^\prime}(\tilde{w})\subset\Lambda^{p^\prime}(w^{1-p^\prime})\/$  and by
Remark~\ref{obs: II.4.8} (iii), we conclude the equality of the four spaces in the statement.

(iii) If $f\in\M(X)\/$ and  $w\/$ is increasing, the function $g_0=(f^\ast w^{-1})^{p^\prime-1}\/$
is decreasing in  $\R^+\/$ and we have that 
\begin{eqnarray*} 
\|f\|_{\Lambda^p(w)^\prime}&=&\sup_{g\d}{\displaystyle\int_0^{\infty} f^\ast(t)
g(t)\,dt\over\bigg(\displaystyle\int_0^{\infty} g^p(t)
w(t)\,dt\bigg)^{1/p}}\\
&\ge&{\displaystyle\int_0^{\infty} f^\ast(t)
g_0(t)\,dt\over\bigg(\displaystyle\int_0^{\infty} g_0^p(t) w(t)\,dt\bigg)^{1/p}}
=\|f\|_{\Lambda^{p^\prime}(w^{1-p^\prime})}.
\end{eqnarray*} 
Since the opposite inequality is always true  (Remark~\ref{obs: II.4.8}  (iii)), we get
the result.$\qquad\qed$
\edem

The following result describes the biassociate space of  $\Lambda^p\/$ in the case $p\le1\/$. See
also \cite{CPSS}.

\bteo
\label{te: II.4.15} 
Let $X\/$ be nonatomic and  let  $w\/$ be a weight in  $\R^+\/$. Let $0<p\le1\/$ and $W_p\/$
the greatest concave minorant of  $W^{1/p}\/$ (Remark~\ref{obs: II.4.8} (iv)). Then there exists a decreasing
weight 
$\tilde{w}_p\/$ with  ${1\over2}\widetilde{W}_p\le W_p\le\widetilde{W}_p\/$ and such that 
$\lox{p}{X}{w}^{\prime\prime}=\lox{1}{X}{\tilde{w}_p}\/$. Moreover,
$$
{1\over2}\|\cdot\|_{\lox{1}{X}{\tilde{w}_p}}\le\|\cdot\|_{\lox{p}{X}
{w}^{\prime\prime}}\le\|\cdot\|_{\lox{1}{X}{\tilde{w}_p}}.
$$
\eteo

\bdem
Let us note that  $\|g\|_{ \Lambda^p(w)^\prime }=\sup_{t>0}W^{-1/p}(t)\displaystyle\int_0^t
g^\ast(s)\,ds\/$ (see Theorem~\ref{te: II.4.7}) and this norm is less than or equal to 
$1\/$, if and only if   $\displaystyle\int_0^t g^\ast(s)\,ds\le W^{1/p}(t),\;t>0\/$. By Lemma~\ref{le:
II.2.4} and  Proposition~\ref{pro: II.2.5}, we obtain, for $f\in\M(X)\/$,
\begin{eqnarray*}
\|f\|_{ \Lambda^p(w)^{\prime\prime} }&=&\sup\bigg\{\int_0^\infty f^\ast(t) g^\ast(t)\,dt\,:\,
\|g\|_{\Lambda^p(w)^\prime}\le1\bigg\}\\ &=&\sup\bigg\{\int_0^\infty f^\ast(t)
\tilde{w}(t)\,dt\,:\,\tilde{w}\d,\;\int_0^t \tilde{w}(s)\,ds\le W^{1/p}(t),\;\forall t>0\bigg\}\\
&=&\sup\bigg\{\int_0^\infty\widetilde{W}\big(\lambda_f(t)\big)\,dt\,:\,\widetilde{W}\in\Upsilon
\bigg\},
\end{eqnarray*}
where
$\Upsilon=\big\{\widetilde{W}\,:\,\tilde{w}\d,\,\widetilde{W}\le W^{1/p}\big\}\/$. Every function
in 
$\Upsilon\/$ is quasi-concave and it is majorized by   $W^{1/p}\/$. Hence, 
$\widetilde{W}\le W_p\/$ for every $\widetilde{W}\in\Upsilon\/$.  On the other hand $W_p\/$  is
quasi-concave and there exists a concave function   
$\widetilde{W}_p(t)=\displaystyle\int_0^t\tilde{w}_p(s)\,ds,\,t>0,\/$ with ${1\over2}\widetilde{W}_p\le
W_p\le\widetilde{W}_p\/$ (see \cite{BS}). In particular $\tilde{w}_p\d\/$ and
${1\over2}\widetilde{W}_p\in\Upsilon\/$. It follows that
\begin{eqnarray*}
{1\over2}\|f\|_{\lox{1}{X}{\tilde{w}_p}}&=&\int_0^\infty{1\over2}\widetilde{W}_p
\big(\lambda_f(t)\big)\,dt\le\|f\|_{
\Lambda^p(w)^{\prime\prime}
}\\ 
&\le&\int_0^\infty\widetilde{W}_p\big(\lambda_f(t)\big)\,dt=\|f\|_{\lox{1}{X}{\tilde{w}_p}}.
\qquad\qed
\end{eqnarray*}
\edem

It is convenient now to recall  the definition of the discrete Hardy operator  $A_d\/$
that has already appeared in (\ref{dischardy}). This operator acts on sequences $f=\big(f(n)\big)_n\/$  by 
$$
A_d f(n)={1\over n+1}\sum_{k=0}^n f(k),\qquad n=0,1,2,\dots
$$ 
We can also describe the associate space when  $X\/$ is  resonant  totally atomic.  If
$\mu(X)<\infty,\;\lox{p,q}{X}{w}=L^\infty(X)\/$ (with equivalent norms) and in this case it is not
of interest. If
$\mu(X)=\infty\/$ there are only, up to isomorphisms, two spaces: $d(\Omega,p)\/$ and
$d^\infty(\Omega,p)\/$ (Example~ \ref{Eje2.3} (iv)) and the most interesting case is when 
$\Omega\notin\ell^1\/$. Up to now, the space
$d(\Omega,p)^\prime\/$ had been identified only in some special cases. In \cite{Al} it is solved
for 
$p\ge1\/$ with $\Omega\/$  decreasing and regular. In \cite{Po} and \cite{NO} the case 
$0<p<1\/$ and 
$\Omega\d\/$ is studied, while  in \cite{AEP} the associate is described under the condition 
$\Omega\u\/$  unbounded (this seems to be the unique reference where the sequence
$\Omega\/$ is not decreasing). Here, we shall identify these   spaces in the most general
case:  

\bteo
\label{te: II.4.16} 
Let $W_n=\sum_{k=0}^n\Omega_k,\;n=0,1,2,\dots\/$ Then, for every 
$f=\big(f(n)\big)_{n=0}^\infty\subset\C\/$,
\begin{enumerate}
\item[(i)] If $0<p\le1\/$,
$$
\|f\|_{d(\Omega,p)^\prime}=\sup_{n\ge0}\,W_n^{-1/p}\sum_{k=0}^n f^\ast(k).
$$

\item[(ii)] If $1<p<\infty\/$ and $\Omega\notin\ell^1\/$,
$$
C_1 \|f\|_{d(\Omega,p)^\prime}\le\bigg(\sum_{n=0}^\infty\big(A_d
f^\ast(n)\big)^{p^\prime}\widetilde{\Omega}_n\bigg)^{1/p^\prime}\le C_2
\|f\|_{d(\Omega,p)^\prime},
$$ 
with
$\widetilde{\Omega}_0=W_0^{1-p^\prime},\;\widetilde{\Omega}_n=(n+1)^{p^\prime}
\big(W_{n-1}^{1-p^\prime}-W_n^{1-p^\prime}\big),\;n\ge1\/$,
and the constants 
$C_1$ and $C_2\/$ depending only on  $p\/$.

\item[(iii)] If $0<p<\infty\/$,
$$
\|f\|_{d^\infty(\Omega,p)^\prime}=\sum_{n=0}^\infty f^\ast(n)W_n^{-1/p}=\|f\|_{d(W^{-1/p},1)}.
$$
\end{enumerate}
\eteo

\bdem
By definition
$$
\|f\|_{d(\Omega,p)^\prime}=\sup_g{\sum_{n=0}^\infty
f(n)g(n)\over\|g\|_{d(\Omega,p)}}=\sup_{g\d}{\sum_{n=0}^\infty
f^\ast(n)g(n)\over\bigg(\sum_{n=0}^\infty g(n)^p\Omega_n\bigg)^{1/p}}.
$$ 
Then  (i) can be directly deduced  applying the first part of  Theorem~\ref{te: I.6.11}. Applying the second
expression of the same theorem in  (ii), we obtain, for $p>1\/$,
$$
\|f\|_{d(\Omega,p)^\prime}^{p^\prime}\approx\int_0^\infty\bigg({\widetilde{V}(t)\over\widetilde{W}(t)}
\bigg)^{p^\prime}\tilde{w}(t)\,dt,
$$
where, $\tilde{v}=\sum_{n=0}^\infty f^\ast(n)\chi_{[n,n+1)},$ $\tilde{w}=\sum_{n=0}^\infty
\Omega_n\chi_{[n,n+1)},$ $\widetilde{V}(t)=\displaystyle\int_0^t\tilde{v}(s)\,ds,$ $\widetilde{W}(t)=
\displaystyle\int_0^t
\tilde{w}(s)\,ds,\;t>0\/$.
Being  $f^\ast\/$ decreasing we have, for $n\ge1$ and  $t\in[n,n+1)\/$,
$$
{1\over
2}\big(f^\ast(0)+\dots+f^\ast(n)\big)\le\widetilde{V}(t)\le\big(f^\ast(0)+\dots+f^\ast(n)\big).
$$
On the other hand
\begin{eqnarray*}
\int_n^{n+1}{\tilde{w}(t)\over\widetilde{W}^{p^\prime}(t)}\,dt&=&\int_0^1\Omega_n\big(\Omega_0+\dots+
\Omega_{n-1}+\Omega_n
t\big)^{-p^\prime}\\
&=&{W_{n-1}^{1-p^\prime}-W_n^{1-p^\prime}\over
p^\prime-1}=C_p{\widetilde{\Omega}_n\over(n+1)^{p^\prime}}.
\end{eqnarray*}
Hence, 
\begin{eqnarray*}
\|f\|_{d(\Omega,p)^\prime}^{p^\prime}&\approx&\int_0^1\bigg({\widetilde{V}(t)\over\widetilde{W}(t)}
\bigg)^{p^\prime}\tilde{w}(t)\,dt+\sum_{n=1}^\infty\int_n^{n+1}\bigg({\widetilde{V}(t)\over\widetilde{W}(t)}
\bigg)^{p^\prime}\tilde{w}(t)\,dt\\
&\approx& f^\ast(0)^{p^\prime}\Omega_0^{1-p^\prime}+\sum_{n=1}^\infty
\bigg(\sum_{k=0}^n f^\ast(k)\bigg)^{p^\prime}{\widetilde{\Omega}_n\over(n+1)^{p^\prime}}\\
&=&\sum_{n=0}^\infty\big(A_d f^\ast(n)\big)^{p^\prime}\widetilde{\Omega}_n.
\end{eqnarray*}

It only remains to prove  (iii). Note that
\begin{eqnarray*}
\|f\|_{d^\infty(\Omega,p)^\prime}&=&\sup\bigg\{\sum_{n=0}^\infty
f^\ast(n)g^\ast(n)\,:\,\|g\|_{d^\infty(\Omega,p)}\le1\bigg\}\\ &=&\sup\bigg\{\sum_{n=0}^\infty
f^\ast(n)g^\ast(n)\,:\,g^\ast(k)\le W_k^{-1/p},\;\forall k\ge0\bigg\},
\end{eqnarray*}
and this supremum is attained at the sequence 
$W^{-1/p}=\big(W_n^{-1/p}\big)_n\/$. Hence,
$$
\|f\|_{d^\infty(\Omega,p)^\prime}=\sum_{n=0}^\infty f^\ast(n)W_n^{-1/p}.\qquad\qed
$$
\edem

\bobs
\label{obs: II.4.17}

(i) Both  $d(\Omega,p)\/$ and $d^\infty(\Omega,p)\/$  are continuously embedded in 
$\ell^\infty\/$ and, thus,  $\ell^1\subset d(\Omega,p)^\prime\/$ (and also $\ell^1\subset
d^\infty(\Omega,p)^\prime\/$) for $0<p<\infty\/$. In particular, the associate space of these
spaces is never trivial.

(ii) The argument used in Remark~\ref{obs: II.4.8} (iii) is also useful here and,
for
$1<p<\infty\/$,
$$
d(\Omega^{1-p^\prime},p^\prime)\subset d(\Omega,p)^\prime,
$$ 
and 
$\|f\|_{d(\Omega,p)^\prime}\le\|f\|_{d(\Omega^{1-p^\prime},p^\prime)}\/$ for every sequence 
$f\/$. If $\Omega\u\/$ we can  use the argument in the proof of   
Theorem~\ref{te: II.4.14} (iii) and conclude (as in \cite{AEP}) 
that the two previous spaces (and their norms) are equal.
\eobs

In \cite{Al}, Allen proves, in the case $\Omega\d,\;p>1\/$, that
$d(\Omega,p)^\prime=d(\Omega^{1-p^\prime},p^\prime)\/$ if $\Omega\/$ is regular. The following
theorem extends this result to the general case.

\bteo
\label{te: II.4.18}
 Let $1<p<\infty\/$ and $\Omega\notin\ell^1\/$. Then
$d(\Omega,p)^\prime=d(\Omega^{1-p^\prime},p^\prime)\/$ if and only if 
$$
\sum_{k=0}^n\Omega_k^{1-p^\prime}\le C\sum_{k=0}^n\big(A_d\Omega(k)\big)^{1-p^\prime},\qquad
n=0,1,2,\dots
$$ 
In particular
$d(\Omega,p)^\prime=d(\Omega^{1-p^\prime},p^\prime)\/$ if $\Omega\/$ is regular.
\eteo

\bdem
The equality of the two   spaces is equivalent (see Remark~\ref{obs: II.4.17}) to the
embedding 
$d(\Omega,p)^\prime\subset d(\Omega^{1-p^\prime},p^\prime)\/$ that, by Theorem~\ref{te:
II.4.16} (ii), holds if and only if the inequality 
\bequ
\label{II.4.19} 
\bigg(\sum_{n=0}^\infty g(n)^{p^\prime}\Omega_n^{1-p^\prime}\bigg)^{1/p^\prime}\le
C\bigg(\sum_{n=0}^\infty
\big(A_dg(n)\big)^{p^\prime}\widetilde{\Omega}_n\bigg)^{1/p^\prime} 
\eequ
holds for every positive and decreasing sequence $g=\big(g(n)\big)_n\/$. Here
$\widetilde{\Omega}\/$ is the sequence defined by $\widetilde{\Omega}_0=\Omega_0^{1-p^\prime}\/$,
$$
{\widetilde{\Omega}_n\over(n+1)^{p^\prime}}=\bigg(\sum_{k=0}^{n-1}\Omega_k\bigg)^{1-p^\prime}-
\bigg(\sum_{k=0}^n\Omega_k\bigg)^{1-p^\prime},\qquad n=1,2,3,\dots
$$
The class of positive and decreasing sequences  is regular (Definition~\ref{def: I.2.4}) in $\N^\ast\/$ and
we can use Theorem~\ref{twoop} with
$T_0=A_d\/$ (which is order continuous, linear and positive) to characterize  (\ref{II.4.19}).
The condition  is obtained \lq\lq applying" the inequality to the sequences 
$(1,1,\dots,1,0,0,\dots)\/$ (decreasing characteristic functions in $\N^\ast\/$) and it is
equivalent to 
\begin{eqnarray*}
\sum_{k=0}^n\Omega_k^{1-p^\prime}\lesssim
\Omega_0^{1-p^\prime}&+&\sum_{k=1}^n(k+1)^{p^\prime}\big(W_{k-1}^{1-p^\prime}-W_k^{1-p^\prime}\big)\\
&+&(n+1)^{p^\prime}\sum_{k=n+1}^\infty\big(W_{k-1}^{1-p^\prime}-W_k^{1-p^\prime}\big),
\end{eqnarray*}
$n=1,2,\dots\/$, with $W_k=\sum_{j=0}^k\Omega_j\/$. The second term is equivalent, for
$n\ge1,\/$ to the expression
\begin{eqnarray*}
\Omega_0^{1-p^\prime}+\sum_{k=1}^n(k+1)^{p^\prime}(W_{k-1}^{1-p^\prime}-W_k^{1-p^\prime})
+(n+1)^{p^\prime}W_n^{1-p^\prime}&\approx&\sum_{k=0}^{n-1}\bigg({W_k\over
k+1}\bigg)^{1-p^\prime}\\ &\approx&\sum_{k=0}^n\bigg({W_k\over k+1}\bigg)^{1-p^\prime}.
\end{eqnarray*}
And  the condition of the theorem is obtained. 

The second assert is immediate, since every regular weight trivially satisfies  the given
condition.$\quad\qed$
\edem

In \cite{NO}, M. Nawrocki and A. Orty\'{n}ski prove, for the case $\Omega\/$ decreasing, that the
associate of 
$d(\Omega,p),\;p\le1\/$, is $\ell^\infty\/$ if the condition 
$\sum_{k=0}^n\Omega_k\ge C(n+1)^p\/$ holds. We now extend this result proving that this condition
also works for the  general case. Moreover, this condition is also necessary.

\bteo
\label{te: II.4.20} 
Let $0<p\le1\/$. Then  $d(\Omega,p)^\prime\subset\ell^\infty\/$,  and these two spaces are equal, if
and only if there exists   
$C>0\/$ such that 
$$
\sum_{k=0}^n\Omega_k\ge C(n+1)^p,\qquad n=0,1,2,\dots
$$
\eteo

\bdem
By Theorem~\ref{te: II.4.16} (i), we have that, with $W_n=\sum_{k=0}^n\Omega_k\/$,
$$
\|f\|_{d(\Omega,p)^\prime}\ge W_0^{-1/p}f^\ast(0)=W_0^{-1/p}\|f\|_{\ell^\infty}.
$$
Therefore, 
$d(\Omega,p)^\prime\subset\ell^\infty\/$.

If the inequality of the theorem holds, 
$$
W_n^{-1/p}\sum_{k=0}^nf^\ast(k)={n+1\over W_n^{1/p}}A_d f^\ast(n)\le C^{-1/p}A_d f^\ast(n)\le
C^{-1/p}\|f\|_{\ell^\infty},
$$ 
for all  $n=0,1,2,\dots\/$ 
Hence,
$\|f\|_{d(\Omega,p)^\prime}=\sup_n W_n^{-1/p}\sum_{k=0}^nf^\ast(k)\le C^{-1/p}\|f\|_{\ell^\infty}\/$
and we have that  $\ell^\infty= d(\Omega,p)^\prime\/$.

Conversely, if $\ell^\infty\subset d(\Omega,p)^\prime\/$, then $1\in d(\Omega,p)^\prime\/$ and
$\sup_n{n+1\over W_n^{1/p}}=\|1\|_{d(\Omega,p)^\prime}$
$=C<\infty\/$, which  implies the result.$\qquad\qed$
\edem

The following result describes the biassociate of  $d(\Omega,p)\/$ in the case $p\le1\/$ in an
analogous way as we did in the nonatomic case  (Theorem~\ref{te: II.4.15}). We omit here the
proof.

\bteo
\label{te: II.4.21} 
Let $0<p\le1\/$ and let us denote by $W_n=\sum_{k=0}^n\Omega_k,\;n=0,1,2,\dots\/$
Then there exists a decreasing sequence 
$\widetilde{\Omega}=(\widetilde{\Omega}_n)_n\/$ satisfying
$$
{1\over2}\sum_{k=0}^n\widetilde{\Omega}_k\le\inf_{m\ge0}\max\bigg(1,\,{n+1\over
m+1}\bigg)W^{1/p}_m\le\sum_{k=0}^n\widetilde{\Omega}_k,
$$ 
and such that
$d(\Omega,p)^{\prime\prime}=d\big(\widetilde{\Omega},1\big)\/$ with equivalent norms. More
precisely,
$$
{1\over2}\|\cdot\|_{d(\widetilde{\Omega},1)}\le\|\cdot\|_{d(\Omega,p)^{\prime\prime}}
\le\|\cdot\|_{d(\widetilde{\Omega},1)}\;.
$$
\eteo

Let us study now the topologic dual  and its connection with the associate space.  

\bdefi
If $(E,\|\cdot\|)\/$ is a quasi-normed space,  we define the dual\index{dual space} 
$E^\ast\/$ in the usual way:
$$
E^\ast=\{u:E\to\C\,:\,u\hbox{\rm \ linear and continuous}\}.
$$ 
If $u\in E^\ast\/$ we denote by 
$\|u\|=\|u\|_{E^\ast}=\sup\big\{|u(f)|\,:\,\|f\|\le1,\;f\in E\big\}\in[0,\infty)\/$.
\edefi

The  dual $E^\ast\/$  of a quasi-normed space $\big(E,\|\cdot\|_E\big)\/$  is a Banach space.
If $E^\ast\/$ separates points, every $f\in E\/$ can be identified with the linear and continuous
form 
$\tilde{f}:E^\ast\to\C\/$  such that $\tilde{f}(u)=u(f),\;u\in E^\ast\/$. We have then a continuous
injection,
$$
\big(E,\|\cdot\|_E\big)\hookrightarrow\big(E^{\ast\ast},\|\cdot\|_{E^{\ast\ast}}\big),
$$ 
and the constant of this embedding  is less than or equal to $1\/$ since, for $f\in E\/$,
$$
\|f\|_{E^{\ast\ast}}=\|\tilde{f}\|_{(E^\ast,\|\cdot\|_{E^\ast})^\ast}=\sup_{u\in
E^\ast}{\big|u(f)\big|\over\|u\|_{E^\ast} }\le\|f\|_E.
$$

The  Mackey topology\index{Mackey topology} in   $E\/$ associated to the dual pair  $(E,E^\ast)\/$, is
defined as having as a local basis the convex envelope of the balls in  $E\/$ (see \cite{Ko}). It is the
finest locally convex topology in 
$E\/$ having 
$E^\ast\/$ as its topological dual. The completion 
$\widetilde{E}\/$ with this topology is called   Mackey completion or also  
Banach envelope of 
$E\/$. In a quasi-normed space whose dual separates points, this topology corresponds to the one
induced in     $E\/$ by 
$\big(E^{\ast\ast},\|\cdot\|_{E^{\ast\ast}}\big)\/$ and the  Mackey completion  is then the closure
of $E\/$ in
$\big(E^{\ast\ast},\|\cdot\|_{E^{\ast\ast}}\big)\/$.

In our case $E=\Lambda\/$ is a quasi-normed  Lorentz space. If $f\in\Lambda^\prime\/$
the application $u_f:\Lambda\to\C\/$ such that $u_f(g)=\displaystyle\int_X f(x)g(x)\,d\mu(x)\/$ is obviously
linear and continuous with norm equal to    $\|f\|_{\Lambda^\prime}\/$. Thus,  $\Lambda^\prime\/$
is isometrically isomorphic to a subspace of   
$\Lambda^\ast\/$ and, in fact, we shall identify the functions of  $\Lambda^\prime\/$
as linear and continuous forms in  $\Lambda\/$:
$$
\Lambda^\prime\subset\Lambda^\ast.
$$ 
On the other hand, $\Lambda\/$ is continuously embedded in 
$\big(\Lambda^{\prime\prime},\|\cdot\|_{\Lambda^{\prime\prime}}\big)\/$ (Remark~\ref{Obs:
II.4.2}), and therefore (if $\Lambda^\ast\/$ separates points) we are working with three independent
topologies in
$\Lambda\/$ and two embeddings with constant $1\/$:
\begin{eqnarray*}
\big(\Lambda,\|\cdot\|_{\Lambda}\big)&\hookrightarrow&\big(\Lambda^{\ast\ast},\|
.\|_{\Lambda^{\ast\ast}}\big),\\
\big(\Lambda,\|\cdot\|_{\Lambda}\big)&\hookrightarrow&\big(\Lambda^{\prime\prime},
\|\cdot\|_{\Lambda^{\prime\prime}}\big).
\end{eqnarray*}

\bpro
\label{pro: II.4.23} 
If $\Lambda\/$ is a quasi-normed  Lorentz space whose dual separates points,
$$
\big(\Lambda,\|\cdot\|_{\Lambda}\big)\hookrightarrow
\big(\Lambda,\|\cdot\|_{\Lambda^{\ast\ast}}\big)\hookrightarrow
\big(\Lambda,\|\cdot\|_{\Lambda^{\prime\prime}}\big),
$$
and  $\|f\|_{\Lambda^{\prime\prime}}\le\|f\|_{\Lambda^{\ast\ast}}
\le\|f\|_{\Lambda},\;f\in\Lambda\/$.
If in addition, $\Lambda^\prime=\Lambda^\ast\/$, then
$\Lambda^{\prime\prime}\/$ is isometrically identified with a subspace of
$\Lambda^{\ast\ast}\/$. In particular, 
$\|f\|_{\Lambda^{\prime\prime}}=\|f\|_{\Lambda^{\ast\ast}},$ for every $ f\in\Lambda,\/$  in this case. 
\epro

\bdem
If $f\in\Lambda\/$, then $f\in\Lambda^{\ast\ast}\/$ and
$$
\|f\|_{\Lambda}\ge\|f\|_{\Lambda^{\ast\ast}}=
\sup_{u\in\Lambda^\ast}{|u(f)|\over\|u\|_{\Lambda^\ast}}
\ge\sup_{g\in\Lambda^\prime}{\int|fg|\over\|g\|_{\Lambda^\prime}}=\|f\|_{\Lambda^{\prime\prime}}
.
$$
 
If $\Lambda^\prime=\Lambda^\ast\/$, every continuous and linear form $u\in\Lambda^\ast\/$ is of the
form 
$u(f)=u_g(f)=\displaystyle\int_X f(x)g(x)\,d\mu(x)\/$ with $g\in\Lambda^\prime\/$. Also, to every function 
$f\in\Lambda^{\prime\prime}\/$ we can associate  the  linear form 
$\tilde{f}:\Lambda^\prime=\Lambda^\ast\to\C\/$ defined by
$\tilde{f}(g)=\tilde{f}(u_g)=\displaystyle\int_X f(x)g(x)\,d\mu(x)\/$ and with norm 
$$
\|\tilde{f}\|_{\Lambda^{\ast\ast}}=\sup_{u_g\in\Lambda^\ast}
{\big|\tilde{f}(u_g)\big|\over\|u_g\|_{\Lambda^\ast}}=
\sup_{g\in\Lambda^{\prime}}{\int|fg|\over\|g\|_{\Lambda^{\prime}}}
=\|f\|_{\Lambda^{\prime\prime}}.\quad\qed
$$
\edem

We shall now study, among other, the two following questions: (i) To characterize the weights  $w\/$
for which 
$\Lambda(w)^\ast=\{0\}\/$. (ii)  When is 
$\Lambda^\prime=\Lambda^\ast\/$? As in the case of the Lebesgue space  $L^p\/$,  
Radon-Nikodym  theorem is the key for the following   result.

\bpro
\label{pro: II.4.24} 
Let $\Lambda\/$  be a quasi-normed  Lorentz space on  $X\/$. If
$u\in\Lambda^\ast\/$, there exists a unique function 
$g\in\Lambda^\prime\/$ with $\|g\|_{\Lambda^\prime}\le\|u\|\/$, satisfying
$$
u(f)=\int_X f(x)g(x)\,d\mu(x),\quad f\in\LL.
$$
\epro

\bdem
Let us first assume that  $\mu(X)<\infty\/$ (in this case $\S\subset
L^\infty\subset\Lambda\/$). For every  $E\subset X\/$ let
$$
\sigma(E)=u(\chi_E)\in\C.
$$ 
Then $\sigma\/$ is a complex measure in  $X\/$ (the fact that it is 
 $\sigma$-additive is a consequence of the fact that  $\mu(X)<\infty\/$), absolutely continuous with
respect to the $\sigma$-finite measure    $\mu\/$ of $X\/$, since  $\mu(E)=0\/$
implies
$\|\chi_E\|_{\Lambda}=0\/$ and thus  $u(\chi_E)=0\/$. By Radon-Nikodym theorem, there exists a
function 
$g\in L^1(X)\/$ so that 
$\sigma(E)=\displaystyle\int_E g(x)\,d\mu(x)\/$  for every measurable set  $E\subset X\/$. In particular
(since
$u\/$ is linear) $u(f)=\displaystyle\int_X f(x)g(x)\,d\mu(x)\/$ for $f\in\S\subset\Lambda\/$. This also holds
if
$f\in L^\infty(X)\subset\Lambda\/$ since there exists a sequence 
$(s_n)_n\subset\S\subset\Lambda\/$ converging to  $f\/$ in $\Lambda\/$  and also uniformly (see the
proof of Theorem~\ref{te: II.3.12}).

If $\mu(X)=\infty\/$, the previous argument is true in any subset  
$Y\subset X\/$ of finite measure: we define a complex measure $\sigma_Y\/$ in $Y\/$ and we obtain
the existence of a function 
$g_Y\in L^1(Y)\/$ with $u(f)=\displaystyle\int_X f(x)g(x)\,d\mu(x)\/$ for every $f\in L^\infty(X)\/$
supported in 
$Y\/$. If $Y_1,Y_2\/$ are two sets of finite measure,  necessarily   $\displaystyle\int_E
g_{Y_1}(x)\,d\mu(x)=\displaystyle\int_E g_{Y_2}(x)\,d\mu(x)\/$ for every measurable set  $E\subset Y_1\cap
Y_2\/$ and it follows that 
$g_{Y_1}(x)=g_{Y_2}(x)\/$ a.e.
$x\in Y_1\cap Y_2\/$. Hence, since $X$ is $\sigma$-finite we can assure the existence of a function
$g\in L^1_{\hbox{\decp loc}}(X)\/$ such that  $u(f)=\displaystyle\int_X f(x)g(x)\,d\mu(x)\/$ for every $f\in
L^\infty_0(X)\/$.

To show that  $\|g\|_{\Lambda^\prime}\le\|u\|\/$, let $\alpha\in\M(X)\/$ with $|\alpha|=1,\;\alpha
g=|g|\/$ and let also  $(X_n)_n\/$ be an increasing sequence of sets of finite measure with 
$\bigcup_n X_n=X\/$. If $f\in\Lambda\/$, we consider  $f_n=|f|\chi_{\big\{|f|\le n\big\}\cap X_n}\/$.
Then $0\le f_n\le f_{n+1}\to|f|\/$ and every  $f_n\/$ is bounded and with support in a set of
finite measure. Then,
\begin{eqnarray*}
\int_X|f(x)g(x)|\,d\mu(x)&=&\lim_n\bigg|\int_X\alpha(x) f_n(x) g(x)\,d\mu(x)\bigg|\\
&=&\lim_n\big|u(\alpha
f_n)\big|\le\|u\|\lim_n\|f_n\|_{\Lambda}=\|u\|\|f\|_{\Lambda},
\end{eqnarray*}
from which 
$g\in\Lambda^\prime\/$ and  $\|g\|_{\Lambda^\prime}\le\|u\|\/$.

Finally, the uniqueness of  $g\/$ is also clear, since if  $g_1\/$  is another function with
the same properties than $g\/$, we have that 
$\displaystyle\int_E g(x)\,d\mu(x)=\displaystyle\int_E g_1(x)\,d\mu(x)=u(\chi_E)\/$ for every measurable
set  $E\subset X\/$ of finite measure, which implies that  $g=g_1.\quad\qed$
\edem

\bobs
 If $\mu(X)=\infty\/$ and $w\in L^1\/$ every simple function, regardless of the measure of its support, is
in 
$\Lambda_X(w)\/$. We can then define, for every linear form $u\in\Lambda^\prime\/$,
$$
\sigma(E)=u(\chi_E),\quad E\subset X\hbox{\ measurable}.
$$ 
but this set function, defined
on the whole  $\sigma$-algebra of $X\/$,  is not $\sigma$-additive in general. 
To have this property, we need that, for every family of disjoint measurable sets
$(E_n)_n\/$, the functions
$\chi_{\bigcup_1^N E_n}\/$ have to converge to
$\chi_{\bigcup_1^\infty E_n}\/$ in $\Lambda$. Contrary to what happens if $\mu(X)<\infty$, this is not true
in general in this case. For example, in
$\Lambda^1_{\R}\big(\chi_{(0,1)}\big)\/$ the sets $E_n=(n,n+1)\/$ give a counterexample.

This explains why we need the restriction  \lq\lq support of finite measure" in the previous
proposition.
\eobs

\bcor
 If $\Lambda\/$  is a quasi-normed Lorentz space on
$X,\;\Lambda^\prime=\{0\}\/$ if and only if every functional  $u\in\Lambda^\ast\/$ is zero on 
$\LL(X)\/$. In particular, if $X\/$ is resonant, $\Lambda_X^\ast\/$ separates points if and only if 
$\Lambda_X^\prime\neq\{0\}\/$.
\ecor

\bcor
\label{cor: II.4.27} 
If $\Lambda\/$  is a quasi-normed Lorentz space  and  $f\in\Lambda\/$
has  absolutely continuous norm, then 
$\|f\|_{\Lambda^{\ast\ast}}=\|f\|_{\Lambda^{\prime\prime}}\/$.
\ecor

\bdem
Since  $f\/$ is  pointwise limit  of a sequence $(f_n)_n\/$ in $\LL\/$ with
$|f_n|\le|f|\/$ and it has absolutely continuous norm, we have that  $f=\lim_n f_n\/$  in
 $\Lambda\/$
and, by continuity,  (Proposition~\ref{pro: II.4.23}), we also have the convergence in 
$\Lambda^{\prime\prime}\/$ and in 
$\Lambda^{\ast\ast}\/$. Hence, it is sufficient to prove the result assuming that  $f\in\LL\/$. By 
Proposition~\ref{pro: II.4.24}, for every  $u\in\Lambda^\ast\/$ there exists $g\in\Lambda^\prime\/$
with
$\|g\|_{\Lambda^{\prime}}\le\|u_g\|_{\Lambda^\ast}=\|u\|_{\Lambda^\ast}\/$ and such that
$u(f)=u_g(f)=\displaystyle\int_X f(x)g(x)\,d\mu(x)\/$. Then, 
$$
\|f\|_{\Lambda^{\ast\ast}}=\sup_{u\in\Lambda^\ast}{|u(f)|
\over\|u\|_{\Lambda^\ast}}=\sup_{g\in\Lambda^\prime}{|u_g(f)|
\over\|u_g\|_{\Lambda^\ast}}=\sup_{g\in\Lambda^\prime}{\displaystyle\int|f(x)g(x)|\,d\mu(x)\over\|g\|_{\Lambda^\prime}}
=\|f\|_{\Lambda^{\prime\prime}}.\quad\qed
$$
\edem

We can give now a representation of the    dual of  $\Lambda^p_X\/$  for an arbitrary 
$\sigma$-finite measure space $X$.

\bteo
If $0<p<\infty\/$ and $W\in\Delta_2(X)\/$,
$$
\Lambda^p_X(w)^\ast=\Lambda^p_X(w)^\prime\oplus\Lambda^p_X(w)^s\,,
$$ 
where
$$
\Lambda^p_X(w)^s=\big\{u\in\Lambda^p(w)^\ast\,:\,|u(f)|\le C\lim_{t\to\infty}f^\ast(t),\forall
f\in\Lambda^p_X(w)\big\}.
$$ 
This last subspace is formed by functionals that are zero 
on the functions whose support is of finite measure and it is not zero only if 
$\mu(X)=\infty\/$ and 
$w\in L^1\/$.
\eteo

\bdem
That every functional in  $\Lambda^p(w)^s\/$  is zero on functions with support of finite measure
is immediate from the definition, since for such a function $f\/$ we have that 
$\lim_{t\to\infty}f^\ast(t)=0\/$. This also tells us that 
$\Lambda^p(w)^\prime\cap\Lambda^p(w)^s=\{0\}\/$ since if  $u=u_g\/$ is in the previous intersection
we have that $\displaystyle\int_X f(x)g(x)\,d\mu(x)=u(f)=0\/$ for every  $f\in\Lambda^p\/$ with support of
finite measure. Then it follows that $\displaystyle\int_X f(x)g(x)\,d\mu(x)=0\/$ for every 
$f\in\Lambda^p\/$ and hence 
$u=0\/$. Let us see now the decomposition $\Lambda^\ast=\Lambda^\prime+\Lambda^s\/$. To this end, let 
$u\in\Lambda^p(w)^\ast\/$. By Proposition~\ref{pro: II.4.24}, there exists $g\in\Lambda^p(w)^\prime\/$ such
that the continuous linear functional 
$u_g(f)=\displaystyle\int_X f(x)g(x)\,d\mu(x)\/$ coincides with $u\/$ in $\LL\/$. If $f\in\Lambda^p\/$ and 
$Y=\{f\neq0\}\/$ has finite measure, then $u(f)=u_g(f)\/$ since $u$ and $u_g\/$ are continuous linear forms
on 
$\Lambda^p_Y(w)\/$  and coincide in $\LL(Y)\/$, which is dense in  $\Lambda^p_Y(w)\/$ (Theorem~\ref{te:
II.3.13}).  Let
$u^s=u-u_g\/$. Then $u^s\/$ is zero on the functions of  $\Lambda^p\/$ supported in a set of finite
measure. Therefore, if
$f\/$ any function in 
$\Lambda^p\/$ and  
$a=\lim_{t\to\infty}f^\ast(t)\/$ and $f_n=f\chi_{\{|f|\le a+1/n\}},\;n=1,2,\dots,\/$ we have 
$u^s(f)=u^s(f_n)\/$ for every $n\/$ (since $f-f_n\/$  has support of finite measure). Thus,
$$
|u^s(f)|=|u^s(f_n)|\le\|u^s\|\|f_n\|_{\Lambda^p},\quad n=1,2,\dots
$$ 
We now have two cases: (i)
$w\notin L^1\/$. Then $a=0\/$ and  $\Lambda^p(w)\/$ has absolutely continuous norm
(Theorem~\ref{te: II.3.5}). Since $|f_n|\le|f|\/$ and  $|f_n|\to0\/$ a.e. we get that
$\|f_n\|_{\Lambda^p}\to0\/$ and 
$|u^s(f)|=0=a\/$. (ii) $w\in L^1\/$. Then
$\|f_n\|_{\Lambda^p}\le\|f_n\|_\infty\|w\|_1^{1/p}\le(a+1/n)\|w\|_1^{1/p}\/$ and we have that
$|u^s(f)|\le\|u^s\|\|w\|_1^{1/p}a=Ca=C\lim_{t\to\infty}f^\ast(t)\/$. In any of these two cases,
there exists $C\in(0,+\infty)\/$ (independent of $f\/$) such that  $|u^s(f)|\le
C\lim_{t\to\infty}f^\ast(t)\/$. That is,  $u^s\in\Lambda^p(w)^s\/$ and we have 
$\Lambda^\ast=\Lambda^\prime\oplus\Lambda^s\/$.

It remains to prove that  $\Lambda^p(w)^s\neq\{0\}\/$ if and only if  $\mu(X)=\infty\/$ and $w\in
L^1\/$. If
$\mu(X)<\infty\/$ or if  $w\notin L^1\/$ every function
$f\in\Lambda^p(w)\/$ satisfies $\lim_{t\to\infty}f^\ast(t)=0\/$ and therefore, for
$u\in\Lambda^s\/$, we get $u(f)=0,$ for every $ f\in\Lambda^p\/$. That is $u=0\/$ or equivalently
$\Lambda^s=\{0\}\/$. If $\mu(X)=\infty\/$ and $w\in L^1\/$, the functional
$p(f)=\lim_{t\to\infty}f^\ast(t),\;p:\Lambda^p\to[0,+\infty),\/$ is a seminorm:
\begin{eqnarray*} 
(i)&\quad &p(\lambda f)=|\lambda p(f)|\quad\hbox{\rm (obvious)},\\ (ii)&\quad&
p(f+g)=\lim_{t\to\infty}(f+g)^\ast(t)\le\lim_t\big(f^\ast(t/2)+g^\ast(t/2)\big)=p(f)+p(g).
\end{eqnarray*} 
Since
$p(1)=1\/$ (the constant function $1\/$ is in
$\Lambda^p(w)\/$ in this case) $p\/$ is not zero and there exists a nonzero linear form  $u\/$
on 
$\Lambda^p\/$ satisfying
\bequ
\label{II.4.29}
|u(f)|\le p(f)=\lim_{t\to\infty}f^\ast(t),\quad\forall f\in\Lambda^p 
\eequ
(see \cite{Ru2} 
for example). In particular, $u\/$ is continuous since 
$$
\|f\|_{\Lambda^p(w)}=\bigg(\int_0^\infty(f^\ast(t))^p w(t)\,dt\bigg)^{1/p}\ge p(f)\|w\|_1^{1/p},
\quad\forall
f,
$$
and it follows that $|u(f)|\le p(f)\le C\|f\|_{\Lambda^p},$ for every $ f\/$. Finally,   
(\ref{II.4.29}) tells us that  $u\in\Lambda^s\/$ and hence $\Lambda^s\neq\{0\}.\quad\qed$
\edem

We shall state two immediate consequences of the last theorem that solve, in the case
$\Lambda=\Lambda^p\/$, the two questions about duality we were looking for.

\bcor
\label{cor: II.4.30}
 Let $0<p<\infty\/$ and  $W\in\Delta_2(X)\/$.
\begin{enumerate}
\item[(i)] If $\mu(X)<\infty\/$ or $w\notin L^1\/$,
$$
\Lambda^p_X(w)^\prime=\Lambda^p_X(w)^\ast.
$$ 
In particular
$d(\Omega,p)^\ast=d(\Omega,p)^\prime\/$  if $\Omega\notin\ell^1\/$ and $d(\Omega,p)\/$ is
quasi-normed.

\item[(ii)] If $\mu(X)=\infty\/$ and $w\in L^1\/$,
$$
\Lambda^p_X(w)^\prime\varsubsetneq\Lambda^p_X(w)^\ast.
$$  
In particular
$\Lambda^p_X(w)^\ast\neq\{0\}\/$ in this case.
\end{enumerate}
\ecor

As a consequence of this last result, Theorems~\ref{te: II.4.16}, \ref{te: II.4.18}, \ref{te:
II.4.20} and Remark~\ref{obs: II.4.17},  remain true (in the case $d(\Omega,p)\/$
quasi-normed,
$\Omega\notin\ell^1\/$) if we substitute  $d(\Omega,p)^\prime\/$ by
$d(\Omega,p)^\ast\/$. 

\bcor
 Let $X\/$ be a nonatomic measure space and let $W\in\Delta_2\/$.
\begin{enumerate}
\item[(i)] If $0<p\le1,\;\Lambda^p_X(w)^\ast=\{0\}\/$ if and only if  the two following conditions hold:
\begin{eqnarray*} 
i.1)&\quad&\mu(X)<\infty\hbox{\ or\ else\ }\mu(X)=\infty\hbox{\ and\ }w\notin L^1, \hbox{\ and }\\
i.2)&\quad&\sup_{0<t<1}{t^p\over W(t)}=\infty.
\end{eqnarray*}

\item[(ii)] If $1<p<\infty,\;\Lambda^p_X(w)^\ast=\{0\}\/$ if and only if it holds:
\begin{eqnarray*}  
ii.1)&\quad&\mu(X)<\infty\hbox{\ or\ else\ }\mu(X)=\infty\hbox{\ and }w\notin L^1, \hbox{\ and }\\
ii.2)&\quad&\int_0^1\bigg({t\over W(t)}\bigg)^{p^\prime-1}\,dt=\infty.
\end{eqnarray*} 
\end{enumerate}
\ecor

The spaces $d(\Omega,p)\/$ have been studied up to now in the case $\Omega\d\/$. As we shall see,
this condition asserts (if
$p\ge1\/$) that
$\|\cdot\|_{d(\Omega,p)}\/$ is a norm. In the case $\Omega\d,p<1\/$ there has been interest in
finding the     Banach envelope of $d(\Omega,p)\/$ (see for example  \cite{Po} and 
 \cite{NO}). The following theorem solves this question in the general case of $w\notin
L^1\/$.

\bteo
\label{te: II.4.32} Let $0<p\le1\/$. 
\begin{enumerate}
\item[(i)] Let $X\/$ be nonatomic and let us assume that  
\begin{enumerate}
\item[(a)] $ W\in\Delta_2,$ 
\item[(b)] $ w\notin L^1\hbox{\ or\ }\mu(X)<\infty,$
\item[(c)] $\sup_{0<t<1}{t^p\over  W(t)}<\infty.$
\end{enumerate}
Then, there exists a decreasing weight $\tilde{w}\/$
with  
$$
\widetilde{W}(t)\approx\inf_{s>0}\max(1,t/s)W^{1/p}(s),t>0,\/
$$ 
and such that the 
 Mackey topology in  $\Lambda^p_X(w)\/$ is the one induced by the norm
$\|\cdot\|_{\Lambda^1_X(\tilde{w})}\/$. The Banach envelope of  $\Lambda^p(w)\/$ is
$\Lambda^1(\tilde{w})\/$ if $\tilde{w}\notin L^1\/$, and it is the space
$$
\Lambda^1(\tilde{w})_0=\big\{f\in\Lambda^1(\tilde{w})\,:\,\lim_{t\to\infty}f^\ast(t)=0\big\}
$$ 
otherwise.

\item[(ii)] If $\Omega\notin\ell^1\/$ and the space $d(\Omega,p)\/$ is quasi-normed, there exists a
decreasing sequence
$\widetilde{\Omega}=(\widetilde{\Omega}_n)_n\/$ such that
$$
\sum_{k=0}^n\widetilde{\Omega}_k\approx\sup_{m\ge0}\max\bigg(1,\,{n+1\over
m+1}\bigg)\bigg(\sum_{k=0}^m\Omega_k\bigg)^{1/p},\;n=0,1,2,\dots,
$$ 
and the norm 
$\|\cdot\|_{d(\widetilde{\Omega},1)}\/$  induces the  Mackey topology in $d(\Omega,p)\/$. If
$\widetilde{\Omega}\notin\ell^1\/$   the Mackey completion of 
$d(\Omega,p)\/$ is   $d(\widetilde{\Omega},1)\/$, and this space is  $c_0=c_0(\N^\ast)\/$ otherwise.
\end{enumerate}
\eteo

\bdem
By Corollary~\ref{cor: II.4.30} the dual and the associate space of 
$\Lambda=\Lambda^p_X(w)\/$ ($\Lambda=d(\Omega,p)\/$ in the case (ii)) coincide, and moreover
$\Lambda^p(w)^\prime=\Lambda^p(w)^\ast\neq\{0\}\/$  by Theorem~\ref{te: II.4.11}. It follows then
from Proposition~\ref{pro: II.4.23}
that the  Mackey topology in 
$\Lambda\/$ (the topology of the  bidual norm)  is the one induced by the norm of the biassociate
space, and since $\Lambda^{\prime\prime}\/$ is complete, the Mackey completion 
$\widetilde{\Lambda}\/$ is contained in  $\Lambda^{\prime\prime}\/$. By Theorem~\ref{te: II.4.15} 
(Theorem~\ref{te: II.4.21} in the case (ii)),
$\Lambda^{\prime\prime}=\Lambda^1(\tilde{w})\/$ ($\Lambda^{\prime\prime}=d(\widetilde{\Omega},1)\/$ resp.).
If
$\tilde{w}\notin L^1\/$ ($\widetilde{\Omega}\notin\ell^1\/$),  Theorem~\ref{te: II.3.13}  tells us
that 
$\Lambda\/$ is dense in  $\Lambda^{\prime\prime}\/$ (since $\LL\subset\Lambda\/$) and hence
$\widetilde{\Lambda}=\Lambda^{\prime\prime}\/$. If, on the contrary $\tilde{w}\in L^1\/$, 
$\Lambda^1(\tilde{w})_0\/$ is closed in 
$\Lambda^1(\tilde{w})\/$ and it contains  the space  $\Lambda^p(w)\/$. If $f\in\Lambda^1(\tilde{w})_0\/$, the
functions 
$f_n=f\chi_{\{|f|>f^\ast(n)\}}\/$ have support in a set of finite measure and converge to $f\/$
in $\Lambda^1(\tilde{w})\/$. Each of these functions $f_n\/$ has absolutely continuous norm
(Corollary~\ref{co: II.3.6}) and thus   $f_n\/$ is limit in 
$\Lambda^1(\tilde{w})\/$ of function in  $L^\infty_0\/$. It follows that this   space is dense 
in $\Lambda^1(\tilde{w})_0\/$ and since  $L^\infty_0\subset\Lambda^p(w)\/$, we conclude that the 
Banach envelope of 
$\Lambda^p(w)\/$ is $\Lambda^1(\tilde{w})_0\/$. In the atomic case,
$\widetilde{\Omega}\in\ell^1\/$ implies
$\Lambda^{\prime\prime}=d(\widetilde{\Omega},1)=\ell^\infty\/$ and 
$\widetilde{\Lambda}=c_0.\qquad\qed$
\edem

\bobs
Conditions  (a) and  (c) in  (i) of the previous theorem are  natural, since if 
$W\notin\Delta_2\/$ there is no topology in
$\Lambda^p(w)\/$ and, if condition (c) fails, $\Lambda^\ast=\{0\}\/$ (Theorem~\ref{te: II.4.11}) and it
does not make sense to consider the  Mackey topology in  $\Lambda^p(w)\/$  with \lq\lq respect to the dual
pair 
$\big(\Lambda^p(w),\Lambda^p(w)^\ast\big)\/$\rq\rq. The conditions on  $\Omega\/$  in 
(ii) are not at all restrictive, since if   
$\Omega\in\ell^1,\;d(\Omega,1)=\ell^\infty\/$ and the  Mackey completion is, in this case, the same
space
$d(\Omega,p)=\ell^\infty\/$.
\eobs

N. Popa proves in  \cite{Po}, for the case $\Omega\d\/$, that the Mackey completion of 
$d(\Omega,p)\/$ ($0<p<1\/$) is $\ell^1\/$ if and only if  $d(\Omega,p)\subset\ell^1\/$. By Theorem~
\ref{te: II.4.32} we can extend this result to the  general case.

\bcor
 Let $0<p\le1,\;\Omega\notin\ell^1\/$. Then the  Banach envelope of 
$d(\Omega,p)\/$ is $\ell^1\/$ if and only if 
$d(\Omega,p)\subset\ell^1\/$.
\ecor

\bdem
If $d(\Omega,p)\subset\ell^1\/$ then
$\ell^\infty=\big(\ell^1\big)^\prime\subset d(\Omega,p)^\prime\/$ and it follows (Theorem~
\ref{te: II.4.20}) that these last two spaces are equal, and hence  
$d(\Omega,p)^{\prime\prime}=\ell^1\/$. Since
$d(\Omega,p)^{\prime\prime}=d(\widetilde{\Omega},1)\/$ for some sequence 
$\widetilde{\Omega}\/$ (Theorem~\ref{te: II.4.21}), it follows that 
$\widetilde{\Omega}\notin\ell^1\/$ and, by the previous theorem, the  Banach envelope of 
$d(\Omega,p)\/$ is
$d(\Omega,p)^{\prime\prime}=\ell^1\/$.

The converse is immediate.$\qquad\qed$
\edem

Let us now study the weak-type spaces $\Lambda^{p,\infty}\/$. First we shall consider the
nonatomic case $X\/$.

\bteo
\label{te: II.4.35} 
If $0<p<\infty,\;W\in\Delta_2\/$ and $X\/$ is nonatomic, then 
$\Lambda_X^{p,\infty}(w)^\prime=\Lambda_X^{p,\infty}(w)^\ast\/$ if and only if 
these two spaces are zero.
\eteo

\bdem
The sufficiency is obvious. To see the necessity, we only need to prove that 
$\loinx{p}{X}{w}^\prime\neq\loinx{p}{X}{w}^\ast\/$ if the first of these spaces is not  trivial. 

But if $\loinx{p}{X}{w}^\prime\neq\{0\}\/$, the function $W^{-1/p}\/$ is locally integrable in
$[0,\infty)\/$ (Theorem~\ref{te: II.4.11}) and we can define the seminorm 
$$
H(f)=\limsup_{t\to0}{\displaystyle\int_0^t f^\ast(s)\,ds\over\displaystyle\int_0^t
W^{-1/p}(s)\,ds},\qquad f\in\Lambda^{p,\infty}.
$$ 
If $f\in\loinx{p}{X}{w}\/$ we have that 
$f^\ast(t)\le\|f\|_{\Lambda^{p,\infty}}W^{-1/p}(t)\/$, and therefore $H\/$ is well defined and it is
continuous. Moreover, it is not zero because there exists   $f_0\in\Lambda^{p,\infty}\/$
with
$f_0^\ast=W^{-1/p}\/$ in $[0,\mu(X))\/$ and hence  $H(f_0)=1\/$. By  Hahn-Banach theorem, there exists a
nonzero $u\in\Lambda^\ast\/$  such that $|u(f)|\le H(f),$ for every $ f\/$. This linear form is not in 
$\Lambda^\prime\/$, because it is zero over all functions  $f\in\LL\/$:
$$
H(f)=\limsup_{s\to0}{\displaystyle\int_0^s f^\ast(t)\,dt\over\displaystyle\int_0^s
W^{-1/p}(t)\,dt}\le\|f\|_{L^\infty}\limsup_{s\to0}\displaystyle{s\over \displaystyle\int_0^s
W^{-1/p}(t)\,dt}=0.
$$ 
Thus
$\loinx{p}{X}{w}^\ast\neq\loinx{p}{X}{w}^\prime.\qquad\qed$
\edem

\bteo
 Let $0<p<\infty\/$ and  $d^\infty(\Omega,p)\/$ be quasi-normed. Let
$W_n=\sum_{k=0}^n\Omega_k$, $n=0,1,\dots\/$ Then, if
$\Omega\notin\ell^1\/$ and $W^{-1/p}\notin\ell^1\/$, we have that  $d^\infty(\Omega,p)^\ast\neq
d^\infty(\Omega,p)^{\prime}\/$.
\eteo

\bdem The seminorm 
$$
H(f)=\limsup_{n\to\infty}{\sum_{k=0}^n f^\ast(k)\over\sum_{k=0}^n W_k^{-1/p}},\qquad f\in
d^\infty(\Omega,p),
$$ 
is continuous and nonzero, and if 
$W^{-1/p}\notin\ell^1\/$, it is zero on all finite sequences. Thus, using and analogous argument
as in the previous proof,   there exists  $u\in d^\infty(\Omega,p)^\ast\setminus
d^\infty(\Omega,p)^\prime.\qquad\qed$
\edem

In the nonatomic cases, it remains to characterize the identity 
$\Lambda_X^{1,\infty}(w)^\ast=\{0\}\/$. We first need the following results.

\blem
\label{le: II.4.37} 
Let $W\in\Delta_2\/$. If $\mu(X)=\infty\/$ and $w\in L^1\/$, then
$\Lambda^{1,\infty}(w)^\ast\neq\{0\}\/$.
\elem

\bdem
By   hypothesis, we have that if
$H(f)=\lim_{t\to\infty}f^\ast(t),$ then $H:\Lambda^{1,\infty}(w)\to[0,\infty)\/$ 
is a continuous seminorm in 
$\Lambda^{1,\infty}(w)\/$ not identically zero (since  $1\in\Lambda^{1,\infty}(w)\/$ and
$H(1)=1\neq0\/$). By  Hahn-Banach theorem it follows that 
$\Lambda^{1,\infty}(w)^\ast\neq\{0\}.\quad\qed$
\edem

\blem
\label{le: II.4.38}  
Let $X\/$ be nonatomic and  $f\in\Lambda_X^{1,\infty}(w)\/$. If 
$\lim_{t\to\infty}f^\ast(t)=0\/$ there exists  $F\in\Lambda_X^{1,\infty}(w)\/$ satisfying:
\begin{enumerate}
\item[(i)] $ |f(x)|\le F(x)$ a.e. $x\in X$. 
\item[(ii)] $
F^\ast(t)=\|f\|_{\Lambda^{1,\infty}}W^{-1}(t),\quad 0<t<\mu(X).$
\end{enumerate}
\elem

\bdem
We can assume $\|f\|_{\Lambda^{1,\infty}}=1\/$. Let
$A=\{f\neq0\},\,a=\mu(A)\/$. Since $f^\ast(t)\limto{t\to\infty}0\/$, there exists a measure
preserving transformation 
$\sigma:A\to(0,a)\/$ such that  $|f(x)|=f^\ast\big(\sigma(x)\big)\/$ a.e. $x\in
A\/$ (see Theorem~II.7.6 in \cite{BS}). Define
$$
F_0(x)=W^{-1}\big(\sigma(x)\big)\chi_A(x),\quad x\in X.
$$ 
Since $\sigma\/$ is measure preserving,
$F_0^\ast(t)=W^{-1}(t),\;0<t<a\/$. Besides, since $f^\ast\le W^{-1}\/$, we have that 
$$
F_0(x)\ge f^\ast\big(\sigma(x)\big)=|f(x)|\quad\hbox{a.e.\ } x\in A.
$$ 
If $a=\mu(X)\/$ the function
we are looking for is clearly, $F=F_0\/$. If on the contrary
$a<\mu(X)\/$, we take  $0\le F_1\in\M(X\setminus A)\/$, with $F_1^\ast(t)=W^{-1}(a+t),\;0\le
t<\mu(X\setminus A)\/$. It is then immediate to check that the function  
$F=F_0+F_1\chi_{X\setminus A}\/$  satisfies the statement.$\quad\qed$
\edem

Let us define now the seminorm  $\NN\/$ whose properties (Proposition~\ref{Pro: II.4.40}) will be
very useful.

\bdefi
\label{def: II.4.39} 
If $(E,\|\cdot\|)\/$ is a quasi-normed space, we define the seminorm
$\NN:E\to[0,+\infty)\/$ by
$$
\NN(f)=\NN_E(f)=\inf\bigg\{\sum_{k=1}^K\|f_k\|\,:\,(f_k)_k\subset E,\,f=\sum_k f_k\bigg\},\quad
f\in E.
$$
\edefi

It can be proved  that  $\NN\/$ is the bidual norm:
$$
\NN(f)=\sup_{u\in\Lambda^\ast}{|u(f)|\over\|u\|_{\Lambda^\ast}}=\|f\|_{\Lambda^{\ast\ast}}.
$$ 

\bpro
\label{Pro: II.4.40}
 Let $\Lambda\/$ be a quasi-normed Lorentz space, and 
$\NN=\NN_\Lambda\/$. Then for  $f,g\in\Lambda\/$ we have,
\begin{enumerate}
\item[(i)] $\NN(f)\le\|f\|_\Lambda.$ 
\item[(ii)] $ |f|\le|g|\Rightarrow \NN(f)\le
\NN(g).$
\item[(iii)] $\NN(f)=\NN(|f|).$ 
\item[(iv)] $
\NN(f)=\inf\bigg\{\sum_{k=1}^K\|f_k\|_\Lambda\,:\,|f|\le\sum_k|f_k|\bigg\}.$
\end{enumerate} 
Moreover
$\Lambda^\ast=\{0\}\/$ if and only if  $\NN=0\/$.
\epro

\bdem
 (i) is immediate. To see (ii), if $g=\sum_k g_k\/$ we have that
$f=(f/g)g=\sum_k(f/g)g_k\/$,  hence
$\NN(f)\le\sum_k\big\|(f/g)g_k\big\|_\Lambda\le\sum_k\|g_k\|_\Lambda\/$ and thus 
$\NN(f)\le\NN(g)\/$. (iii) and (iv) are  corollaries of  (ii).

Let us now prove the last statement. If there exists  $0\neq u\in\Lambda^\ast\/$, we can find 
$f\in\Lambda\/$ with $u(f)\neq 0\/$. For every finite decomposition
$f=\sum_k f_k\/$ we then have that, $|u(f)|\le\sum_k|u(f_k)|\le\|u\|\sum_k\|f_k\|_\Lambda\/$. Hence
$\sum_k\|f_k\|_\Lambda\ge|u(f)|/\|u\|\/$ and taking the infimum, we obtain
$\NN(f)\ge|u(f)|/\|u\|>0\/$, that is, $\NN\neq0\/$. Conversely, since $\NN\/$ is a continuous
seminorm in 
$\Lambda\/$ (by  (i)), if 
$\NN\neq0\/$ there exists (Hahn-Banach theorem, \cite{Ru2}) a nonzero continuous linear form $u\/$  in  
$\Lambda.\quad\qed$
\edem

In the following result, proved by  A. Haaker in  \cite{Haa} for the case $X=\Rm\/$,
we give a necessary condition  to have that the  dual of
$\Lambda^{1,\infty}(w)\/$ is zero. As we shall see in  Theorem~\ref{te: II.4.45}, this condition is
also sufficient. Our proof is based on the proof of A.  Haaker.

\blem
\label{lem: II.4.41} 
Let $X\/$ be nonatomic and $W\in\Delta_2\/$. If
$\Lambda^{1,\infty}_X(w)^\ast=\{0\}\/$ and  $\epsilon>0\/$, there exists $n\in\N\/$ such that
$$
\int_{nt}^s W^{-1}(r)\,dr\le\epsilon\int_t^s W^{-1}(r)\,dr,\quad 0<t<nt<s\le\mu(X).
$$
\elem

\bdem
We can assume that  $W\/$ is strictly increasing (otherwise, we substitute this function  by
$W(t)(1+2t)/(1+t)\/$). Let $0\le f\in\Lambda^{1,\infty}(w)\/$ with $f^\ast=W^{-1}\/$ in
$\big(0,\mu(X)\big)\/$. By Proposition~\ref{Pro: II.4.40},
$\NN(f)=0\/$  and there exist positive functions $f_1,\dots,f_{n-1}\in\Lambda^{1,\infty}(w)\/$ 
with
$$
f\le\sum_k f_k\/,
$$ 
and such that, if  $a_k=\|f_k\|_{\Lambda^{1,\infty}},\;k=1,\dots,n-1\/$, then
$$
\sum_k a_k<\epsilon.
$$
By  Lemma~\ref{le: II.4.37}, $\mu(X)<\infty\/$  or $\mu(X)=\infty\/$ and $w\notin
L^1\/$. In any case $\lim_{t\to\infty}f_k^\ast(t)=0,\;k=1,\dots,n-1\/$. We can then take functions 
$g_1,\dots,g_{n-1}\in\Lambda^{1,\infty}(w)\/$ with $g_k\ge f_k\/$ a.e. and 
$g_k^\ast(t)=a_kW^{-1}(t),\;0<t<\mu(X)\/$ (using Lemma~\ref{le: II.4.38}). Summarizing,
\begin{eqnarray*}
(i)&\quad &f\le g_1+\dots+g_{n-1}\hbox{\ \ a.e.},\\ (ii)&\quad&
g_k^\ast(t)=a_kW^{-1}(t),\;0<t<\mu(X),\\ (iii)&\quad& a_1+\dots+a_{n-1}<\epsilon.
\end{eqnarray*}  
If  
$s,t>0\/$, with $nt<s<\mu(X)\/$, we define
\begin{eqnarray*} 
 F&=&\big\{x\in X\,:\,W^{-1}(s)<f(x)\le W^{-1}(t)\big\},\\ E_k&=&\big\{x\in
F\,:\,g_k(x)\le a_kW^{-1}(t)\big\},\quad k=1,\dots,n-1,\\ E&=&\bigcap_{k=1}^{n-1}E_k.
\end{eqnarray*}  
Since
$f^\ast=W^{-1}\/$ and this function  is strictly decreasing,
$$
\mu(F)=\lambda_{f^\ast}\big(W^{-1}(s)\big)-\lambda_{f^\ast}\big(W^{-1}(t)\big)=s-t.
$$ 
Moreover,
$$
\mu(F\setminus E_k)=\mu\big(\big\{x\in
F\,:\,g_k(x)>a_kW^{-1}(t)\big\}\big)\le\lambda_{g_k^\ast}\big(a_kW^{-1}(t)\big)=t,
$$ 
and
$\mu(F\setminus E)=\mu\big((F\setminus E_1)\cup\dots\cup(F\setminus E_{n-1})\big)\le(n-1)t\/$.
Hence,
$$
\mu(E)\ge\mu(F)-(n-1)t=s-nt.
$$
Since $E\subset F\/$, if $G=\big\{W^{-1}(s)<f\le
W^{-1}(s-\mu(E))\big\}\/$, the distribution function  of $f\chi_E\/$ majorizes 
$f\chi_G\/$ and we have that
\begin{eqnarray}
\label{II.4.42}
\int_E
f(x)\,d\mu(x)&=&\int_0^\infty(f\chi_E)^\ast(r)\,dr\ge\int_0^\infty(f\chi_G)^\ast(r)\,dr\\ \nonumber
&=&\int_{s-\mu(E)}^s
W^{-1}(r)\,dr\ge\int_{nt}^s W^{-1}(r)\,dr.
\end{eqnarray}
On the other hand,
\begin{eqnarray*}
\int_E f(x)\,d\mu(x)&\le&\sum_{k=1}^{n-1}\int_E
g_k(x)\,d\mu(x)\\
&\le&\sum_k\int_{E_k}g_k(x)\,d\mu(x)=\sum_k\int_0^\infty(g_k\chi_{E_k})^\ast(r)\,dr,
\end{eqnarray*}
but $\mu(E_k)\le\mu(F)=s-t\/$
and in this set
$g_k\le a_kW^{-1}(t)=g_k^\ast(t)\/$. Hence, if $H_k=\{g_k^\ast(t)\ge g_k>g_k^\ast(s)\}\/$, 
the distribution function of  $\chi_{H_k}\/$ majorizes   the   distribution function of
 $\chi_{E_k}\/$ (observe that 
$g_k^\ast=a_kW^{-1}\/$ is strictly decreasing) and we have that
$$
\int_0^{\infty}(g_k\chi_{E_k})^\ast(r)\,dr\le\int_t^s g_k^*(r)\,dr=a_k\int_t^s W^{-1}(r)\,dr\/,
$$
and hence,
$$
\int_E f(x)\,d\mu(x)\le\sum_k a_k\int_t^s W^{-1}(r)\,dr<\epsilon\int_t^s W^{-1}(r)\,dr.
$$ 
This and (\ref{II.4.42}) end the proof.$\quad\qed$  
\edem

\bdefi
A function  $f\in\M(X)\/$ is a step function\index{step function} if it is
 of the form
$$
f=\sum_{n=1}^\infty a_n\chi_{E_n},
$$ 
where $(a_n)_n\subset\C\/$ and $(E_n)_n\/$ is a sequence of  pairwise disjoint
measurable sets in  $X\/$.
\edefi

The proof of the following result is totally analogous to that of Lemma 7 in  \cite{Cw1}, where  the author
considers the case $W(t)=t^{1/p}\/$. 

\blem
\label{lem: II.4.44} 
Every function  $0\le f\in\Lambda^{1,\infty}_X(w)\/$ is majorized by a step  function 
$F\in\Lambda^{1,\infty}\/$.
\elem

Now, we can prove the following result.

\bteo
\label{te: II.4.45}
Let $X\/$ be nonatomic and  $W\in\Delta_2\/$. If $0<p<\infty\/$,
$$
\Lambda^{p,\infty}_X(w)^\ast=\{0\}
$$
if and only if 
\bequ
\label{II.4.46}
\lim_{t\to0^+}\,\sup_{0<r<\mu(X)}{W^{1/p}(tr)\over tW^{1/p}(r)}=0.
\eequ
\eteo 

\bdem
By Remark~\ref{obs: II.2.6}, we can assume $p=1\/$.

{\sl Sufficiency.} We assume that (\ref{II.4.46}) holds and we shall prove that 
$\NN=\NN_{\Lambda^{1,\infty}}=0\/$. By Proposition~\ref{Pro: II.4.40}  we have 
$\Lambda^{p,\infty}_X(w)^\ast=\{0\}$. Let then  $f\in\Lambda^{1,\infty}\/$ and let us see that 
$\NN(f)=0\/$. By Lemma~\ref{lem: II.4.44} and Proposition~\ref{Pro: II.4.40} (ii),  we can assume
that 
$f\/$ is a positive step function:
$$
f=\sum_n a_n\chi_{A_n},\quad (a_n)_n\subset(0,\infty),\;A_n\cap A_m=\emptyset,\;n\neq m. 
$$
Observe that either $\mu(X)<\infty\/$ or, by (\ref{II.4.46}),
$W(\infty)=\infty\/$. Therefore, $\lim_{t\to\infty}f^\ast(t)=0\/$, and hence
$\mu(A_n)<\infty\/$ for every $n=1,2,\dots\/$ Also by (\ref{II.4.46}), for every  $\epsilon>0\/$
there exists
$m\in\N\/$ such that
$$
2^m{W(2^{-m}r)\over W(r)}<\epsilon,\quad 0<r<\mu(X). 
$$
That is,
$$
W^{-1}(2^m t)<\epsilon\,2^{-m}W^{-1}(t),\quad 0<t<2^{-m}\mu(X). 
$$
Let us divide the sets  $A_n\/$ into $2^m\/$ pairwise disjoint measurable subsets
$(A_n^k)_{k=1}^{2^m}\/$  and with equal measure. For every
$k=1,2,\dots,2^m\/$, we have
$$
f_k=\sum_n a_n\chi_{A_n^k}. 
$$
Thus $f=\sum_k f_k\/$ and for  $0<t<2^{-m}\mu(X)\/$ ($f_k^\ast\/$
is zero outside this  interval) we have,
\begin{eqnarray*}
f_k^\ast(t)&=&f^\ast(2^m t)\le W^{-1}(2^m t)\|f\|_{\Lambda^{1,\infty}}\\
&<&\epsilon\,
2^{-m}W^{-1}(t)\|f\|_{\Lambda^{1,\infty}},\quad k=1,\dots,2^m. 
\end{eqnarray*}
Therefore,
$\|f_k\|_{\Lambda^{1,\infty}}\le\epsilon\,2^{-m}\|f\|_{\Lambda^{1,\infty}}\/$ and,
$$
\NN(f)\le\sum_k\|f_k\|_{\Lambda^{1,\infty}}\le\epsilon\,\|f\|_{\Lambda^{1,\infty}}. 
$$
Since this holds for every  $\epsilon>0\/$ we conclude that $\NN(f)=0\/$.

{\sl Necessity.} Let us assume now that $\Lambda^{1,\infty}_X(w)^\ast=\{0\}$. By
Lemma~\ref{lem: II.4.41}, there exists
$n\in\N\/$ such that 
\bequ
\label{II.4.47}
\int_{nt}^s W^{-1}(r)\,dr\le{1\over 2}\int_t^s W^{-1}(r)\,dr,\quad t,s>0,\;nt<s<\mu(X).
\eequ
To see 
(\ref{II.4.46}) we shall prove that for  $K\in\N\/$, 
$$
\sup_{0<r<\mu(X)}{W(tr)\over tW(r)}\le 4n\,2^{-K},
$$ 
if $0<t<n^{-2K}\/$. To this end, let us define for 
$0<r<\mu(X)\/$ 
$$
A_k=\int_{n^k tr}^r W^{-1}(x)\,dx,\quad k=0,1,\dots,K. 
$$
By  (\ref{II.4.47}),
$$
A_0\ge 2A_1\ge 4A_2\ge\dots\ge 2^K A_K
$$ 
and
$$
{1\over 2}\ge{A_1\over A_0}=1-{A_0-A_1\over A_0}\ge1-2^{-K}{A_0-A_1\over A_K}. 
$$
But
$A_0-A_1=\displaystyle\int_{tr}^{ntr} W^{-1}(x)\,dx\le(n-1)trW^{-1}(tr)\/$ and
$A_K\ge(1-n^Kt)rW^{-1}(r)\/$. Hence,
$$
{1\over 2}\ge1-2^{-K}{(n-1)tW(r)\over(1-n^K t)W(rt)}. 
$$
And then,
$$
{W(tr)\over tW(r)}\le 2^{1-K}{n-1\over 1-n^K t}\le 4n\,2^{-K}.\quad\qed
$$
\edem

\bobs 

(i) Condition (\ref{II.4.46}) appears for the first time in  \cite{Haa} where A. Haaker proves the
previous result for the case
$X=\Rm,\;w\notin L^1\/$. 

(ii) In the proof, we have seen that the condition of Lemma~\ref{lem: II.4.41} implies
(\ref{II.4.46}). Therefore such condition is also  equivalent to
$\Lambda^{1,\infty}(w)^\ast=\{0\}\/$.
\eobs

\bcor
 Let $0<p<\infty,\;X\/$  be nonatomic and    $\Lambda=\Lambda_X^{p,\infty}(w)\/$ be 
quasi-normed. Then:
\begin{enumerate}
\item[(i)] If $\lim_{t\to0^+}\,\sup_{0<r<\mu(X)}{W^{1/p}(tr)\over tW^{1/p}(r)}=0\/$,
$$
\{0\}=\Lambda^\prime=\Lambda^\ast.
$$

\item[(ii)] If the previous condition fails,  two cases can occur:
\begin{enumerate}
\item[(ii.1)] If $\displaystyle\int_0^1 W^{-1/p}(s)\,ds=\infty\/$,
$$
\{0\}=\Lambda^\prime\varsubsetneq\Lambda^\ast.
$$

\item[(ii.2)] If $\displaystyle\int_0^1 W^{-1/p}(s)\,ds<\infty\/$,
$$
\{0\}\neq\Lambda^\prime\varsubsetneq\Lambda^\ast.
$$
\end{enumerate}
\end{enumerate}
\ecor

\bdem
It is an immediate consequence of  Theorems \ref{te: II.4.11}, \ref{te: II.4.35}, and  \ref{te:
II.4.45}.$\qquad\qed$
\edem

\beje
\ \par
(i) We can apply the previous result to  identify the dual and associate space of the 
Lorentz spaces $L^{p,q}(X)\/$ ($X\/$ nonatomic), observing that 
$L^{p,q}=\Lambda^q\big(t^{q/p-1}\big)\/$ and $L^{p,\infty}=\Lambda^{p,\infty}(1)\/$. We obtain,
then, by other methods, results already known (see  \cite{Hu, CS, Cw1, Cw2}):
\medskip

(i.1) $L^{p,q}(X)^\prime=L^{p,q}(X)^\ast=\{0\}\/$, if $0<p<1,\;0<q<\infty\/$,

(i.2) $L^{1,q}(X)^\prime=L^{1,q}(X)^\ast=L^\infty(X)\/$, if $0<q\le1\/$,

(i.3) $L^{1,q}(X)^\prime=L^{1,q}(X)^\ast=\{0\}\/$, if $1<q<\infty\/$,

(i.4) $L^{p,q}(X)^\prime=L^{p,q}(X)^\ast=L^{p^\prime,\infty}(X)\/$, if $1<p<\infty,\;0<q\le1\/$,

(i.5) $L^{p,q}(X)^\prime=L^{p,q}(X)^\ast=L^{p^\prime,q^\prime}(X)\/$, if $1<p<\infty,\;1<q<\infty\/$,

(i.6) $L^{p,\infty}(X)^\prime=L^{p,\infty}(X)^\ast=\{0\}\/$, if $0<p<1\/$,

(i.7) $\{0\}=L^{1,\infty}(X)^\prime\varsubsetneq L^{1,\infty}(X)^\ast\/$,

(i.8) $\{0\}\neq L^{p,\infty}(X)^\prime\varsubsetneq L^{p,\infty}(X)^\ast\/$, if $1<p<\infty\/$.

\medskip (ii) If $w=\chi_{(0,1)},\;\Lambda^p_{\Rm}(w)\supset\Lambda^p_{\Rm}(1)=L^p(\Rm)\/$. 
$\Lambda^p_{\Rm}(w)^\ast\neq\{0\}\/$ for every 
$p\in(0,\infty)\/$ while, as for  $L^p,\;\Lambda^p_{\Rm}(w)^\prime\neq\{0\}\/$, if and
only if 
$p\ge1\/$. In the weak-type case, 
$\Lambda^{p,\infty}_{\Rm}(w)^\ast\neq\{0\}\/$ for every $p\in(0,\infty)\/$ and the associate space
is not zero only if $p>1\/$.
\eeje

\section{Normability}
\label{se:norma}

In this section $(X,\mu)\/$ is an arbitrary  $\sigma$-finite measure space. We shall study when 
$\lox{p,q}{X}{w}\/$ is a  Banach space in the sense that there exists a norm in
$\lox{p,q}{X}{w}\/$ equivalent to the functional $\|\cdot\|_{\lox{p,q}{X}{w}}\/$. In the nonatomic
case, some results are already known:  G.G. Lorentz (\cite{Lor},
1951) proved that, for
$p\ge1\/$, $\|\cdot\|_{\lox{p}{(0,l)}{w}}\/$ is a norm if and only if $w\/$ is decreasing (that is,
equal  a.e. to a decreasing  function) in $(0,l)\/$. A. Haaker (\cite{Haa}, 1970)
characterizes the normability of 
$\lox{p,\infty}{\Rm}{w},\,0<p<\infty\/$ (see also \cite{Sor}), and gives some partial results for 
$\lox{p}{\Rm}{w},\,p<\infty\/$. E. Sawyer (\cite{Saw2}, 1990) gives a necessary and sufficient
condition to have that $\loX\/$ is normed in the case
$1<p<\infty,\,X=\Rn\/$, while in  (\cite{CGS}, 1996) the authors solve the case
$0<p<\infty,\,\mu(X)=\infty,\;X\/$ nonatomic. 

In the case $X=\N^\ast\/$ (spaces $d(\Omega,p)\/$) we do not know any previous result concerning
normability. Except some isolated case, in all the references we know on the spaces $d(\Omega,p)\/$,
it is assumed that the sequence
$\Omega\/$ is decreasing which, in the case $p\ge1\/$, ensures that
$\|\cdot\|_{d(\Omega,p)}\/$ is already a norm.

Here, we solve the  general case for a resonant measure space $X\/$  unifying in this way 
all the previous results and including the unknown results for the spaces
$d(\Omega,p)\/$. In the atomic case we shall reduce the problem to  $\Rn\/$. 

Let us start with some general results.

\bteo
\label{te: II.5.1}
If $1\le p<\infty\/$ and $w\d\/$, then $\|\cdot\|_{\Lambda^p_X(w)}\/$ is a norm.
\eteo

\bdem
In the nonatomic case,
\begin{eqnarray*}
\|f+g\|_{\Lambda^p_X(w)}&=&\bigg(\int_0^\infty\big((f+g)^\ast(t)\big)^p
w(t)\,dt\bigg)^{1/p}\\
&=&\sup_{h^\ast=w}\bigg(\int_X|f(x)+g(x)|^p h(x)\,d\mu(x)\bigg)^{1/p},
\end{eqnarray*} 
and the 
triangular inequality is immediate. On the other hand if  $X\/$ is $\sigma$-finite,
we know (retract method) that there exists a nonatomic space  $\bar{X}\/$ and a linear
application 
$F:\M(X)\to\M(\bar{X})\/$ such that $F(f)^\ast=f^\ast\/$ (c.f. \cite{BS})  and the
result follows immediately from the atomic case.$\qquad\qed$
\edem

The following result gives sufficient conditions to have the normability.

\bteo
\label{te: II.5.2}
Let  $X$ be a $\sigma$-finite measure space. Then,
\begin{enumerate}
\item[(i)] if $1\le p<\infty\/$ and  $w\in B_{p,\infty},\;\Lambda^p_X(w)\/$ is normable,

\item[(ii)] if $0<p<\infty\/$ and $w\in B_p,\;\Lambda^{p,\infty}_X(w)\/$ is normable.
\end{enumerate}
\eteo

\bdem (i) If $p=1\/$ and $w\in B_{p,\infty}=B_{1,\infty}\/$, there exists  a decreasing
function
$\tilde{w}\/$ with $\widetilde{W}\approx W\/$ (\cite{CGS}). It follows that
$\|\cdot\|_{\Lambda^1(\tilde{w})}\/$ is a norm in $\Lambda^1(w)\/$ equivalent to the original
one. If
$p>1\/$ and $w\in B_{p,\infty}\/$,  the Hardy operator
$A\/$ satisfies $A:L^p_{\hbox{\decp dec}}(w)\to L^p(w)\/$  and it follows that the functional
$\|f\|=\|f^{\ast\ast}\|_{L^p(w)}\/$ is an equivalent norm to the    original quasi-norm.

(ii) If $w\in B_p\/$ we have that  $A:L^{p,\infty}_{\hbox{\decp dec}}(w)\to L^{p,\infty}(w)\/$ (see
\cite{Sor}) and the functional $\|f\|=\|f^{\ast\ast}\|_{L^{p,\infty}(w)}\/$ is an equivalent
norm.$\qquad\qed$
\edem 

\bobs
\label{obs: 5.3}
The functional $\|\cdot\|_{\Lambda^{p,\infty}}\/$ is not a norm, except in some trivial cases. In
fact, if 
$(X,\mu)\/$ is an arbitrary measure space and there exist two measurable sets $E\subset F\subset
X\/$ with
$1<a=W\big(\mu(F)\big)/W\big(\mu(E)\big)<\infty\/$ (where $w\/$ is an arbitrary weight in 
$\Rm\/$), we have that 
$\|f+g\|_{\Lambda^{p,\infty}}>\|f\|_{\Lambda^{p,\infty}}+\|g\|_{\Lambda^{p,\infty}}\/$ for
$f=a\chi_E+\chi_{F\setminus E},\;g=\chi_E+a\chi_{F\setminus E}\/$ (assuming without loss of
generality that
$\mu(F\setminus E)\le\mu(E)\/$). Then, the following statements are equivalent:
\medskip

(i)\ \  $\|\cdot\|_{\Lambda^{p,\infty}}\/$ is a norm.

(ii)\ \ $\|\cdot\|_{\Lambda^{p,\infty}}=C\|\cdot\|_{L^\infty}\/$.

(iii)\ \ The restriction of $W\/$ to the range of  $\mu\/$ is constant.
\eobs

\bteo
\label{te: II.5.4}
A Lorentz space $\Lambda\/$   is normable if and only if 
$\Lambda=\Lambda^{\prime\prime}\/$,  with equivalent  norms. In particular, every  normable Lorentz
space  $\Lambda\/$  is a  Banach function space with the norm
$\|\cdot\|_{\Lambda^{\prime\prime}}\/$.
\eteo

\bdem  The sufficiency is obvious. Conversely, if $\Lambda\/$ is
normable with a norm $\|\cdot\|\/$, it has to be equivalent to
$\|\cdot\|_{\Lambda^{\ast\ast}}\/$, (by  Hahn-Banach theorem), and by  Corollary~\ref{cor: II.4.27}, we
have 
$\|f\|_{\Lambda}\approx\|f\|_{\Lambda^{\prime\prime}}\/$ for every function  $f\/$ with absolutely
continuous norm in 
$\Lambda\/$. Since every $|f|\in\Lambda\/$ is a  pointwise limit of an increasing sequence of functions
in 
$\LL\/$  (which have absolutely continuous norm by Corollary~\ref{co: II.3.6}  and 
Proposition~\ref{pro: II.3.9}) and the functionals 
$\|\cdot\|_{\Lambda}\/$ and $\|\cdot\|_{\Lambda^{\prime\prime}}\/$  have the Fatou property, it
follows that 
$\|f\|_{\Lambda}\approx\|f\|_{\Lambda^{\prime\prime}}\/$ for every $f\in\Lambda\/$.  Finally, it is
immediate to prove that the norm  
$\|\cdot\|_{\Lambda^{\prime\prime}}\/$ satisfies  the properties of 
Definition~\ref{def4.3}.$\qquad\qed$
\edem 

\bobs
\label{obs: 5.5}
 If $\|\cdot\|_{\Lambda}\/$ is a norm, then
$(\Lambda,\|\cdot\|_{\Lambda})\/$ is a  Banach function space and
$\|\cdot\|_{\Lambda}=\|\cdot\|_{\Lambda^{\prime\prime}}\/$ (see \cite{BS}).
\eobs

From now on, we shall consider the case $X\/$ resonant. We shall see how the previous conditions
are also necessary (in the atomic case). As an  immediate consequence of the previous 
theorem   and the    representation theorem of  Luxemburg (\cite{BS}) we have the following result. 

\bcor
\label{cor: II.5.6}
Let $X\/$ be a nonatomic measure space, $0<p,q\le\infty\/$. If $\Lambda_X^{p,q}\/$ is
normable, so is $\Lambda^{p,q}_{(0,\mu(X))}\/$, and if 
$\|\cdot\|_{\Lambda_X^{p,q}}\/$ is a norm the same holds for  
$\|\cdot\|_{\Lambda_{(0,\mu(X))}^{p,q}}\/$.
\ecor

\blem
\label{lem: II.5.7} 
Let  $A\subset\R\/$ be a measurable set, and let $w\/$ be a weight with  $w=w\chi_{(0,|A|)}\/$.
Then,
\begin{enumerate}
\item[(i)]\ \ $\|\cdot\|_{\Lambda_A(w)}\/$ is a norm if and only if  $\|\cdot\|_{\Lambda_\R(w)}\/$ is a
norm,

\item[(ii)]\ \ $\Lambda_A(w)\/$ is normable if and only if  $\Lambda_{\R}(w)\/$ is normable.
\end{enumerate}
\elem

\bdem
Since $\Lambda_A(w)\subset\Lambda_\R(w)\/$ the  \lq\lq if\rq\rq\  implication is immediate. To see
the other part,  let us assume that
$$
\bigg\|\sum_{j=1}^n f_j\bigg\|_{\Lambda_A(w)}\le C\sum_j\|f_j\|_{\Lambda_A(w)},\quad\forall 
f_1,\dots,f_n\in\Lambda_A(w),
$$ 
and let us first see that we have an analogous inequality   substituting
$A\/$ by an arbitrary set  $Y\subset\R\/$ with $|Y|\le|A|$. By monotonicity
we can assume 
$|Y|<|A|\/$. Then there exists an open set  $G\supset Y\/$, with $b=|G|<|A|\/$. Let
$\tilde{A}\subset A\/$ with $|\tilde{A}|=b\/$ and let  $\sigma_1:\tilde{A}\to(0,b)\/$ be a measure
preserving transformation (Proposition~{7.4} in \cite{BS}). Since $G\/$ is a countable union of 
   intervals, it is easy to construct another measure preserving transformation  $\sigma_2:(0,b)\to
G\/$. Then
$\sigma=\sigma_2\circ\sigma_1:\tilde{A}\to G\/$ preserves the measure and for each 
$f_1,\dots,f_n\in\M(Y)\subset\M(G)\/$, the functions
$\tilde{f}_j=f_j\circ\sigma\in\M(A),\;j=1,\dots,n,\/$ satisfy $\tilde{f}_j^\ast=f_j^\ast,
$ for every $
j\/$ and
$(\tilde{f}_1+\dots+\tilde{f}_n)^\ast=(f_1+\dots+f_n)^\ast\/$. Hence
$$
\|f_1+\dots+f_n\|_{\Lambda_Y(w)}=\|\tilde{f}_1+\dots+\tilde{f}_n\|_{\Lambda_A(w)}\le
C\sum_j\|\tilde{f}_j\|_{\Lambda_A(w)}=C\sum_j\|f_j\|_{\Lambda_Y(w)}.
$$

If now  $f_1,\dots,f_n\/$ are arbitrary measurable functions in $\R\/$ and we choose
$Y\subset\big\{|f_1+\dots+f_n|\ge(f_1+\dots+f_n)^\ast(|A|)\big\}\/$, with
$|Y|=|A|\/$, since $w\/$ is supported in  $(0,|A|)\/$ we have that
$$
\|f_1+\dots+f_n\|_{\Lambda_\R(w)}=\big\|(f_1+\dots+f_n)\chi_Y\big\|_{\Lambda_\R(w)}=\bigg\|\sum_j
f_j\chi_Y\bigg\|_{\Lambda_Y(w)},
$$ 
and it follows that
$$
\|f_1+\dots+f_n\|_{\Lambda_\R(w)}\le C\sum_j\big\|f_j\chi_Y\big\|_{\Lambda_Y(w)}\le      
C\sum_j\|f_j\|_{\Lambda_\R(w)}.
$$
 
If $C=1\/$ we conclude that $\|\cdot\|_{\Lambda_\R(w)}\/$ is a norm. If $1<C<\infty\/$, the
previous inequality proves that 
$\Lambda_\R(w)\/$ is normable (the functional $\NN\/$ defined in \ref{def: II.4.39}
would be, for example, an equivalent norm).$\qquad\qed$
\edem 

The following result characterizes the normability in the nonatomic case.  

\bteo
\label{te: II.5.8}
 Let $(X,\mu)\/$ be a nonatomic measure space,
$0<p<\infty,\;w=w\chi_{(0,\mu(X))}\/$. Then:
\begin{enumerate}
\item[(i)]\ \ $\|\cdot\|_{\Lambda^p_X(w)}\/$ is a norm if and only if  $p\ge1,\;w\d\/$.

\item[(ii)]\ \ $\Lambda^p_X(w)\/$ is normable if and only if  $p\ge1,\;w\in B_{p,\infty}\/$.

\item[(iii)]\ \ $\|\cdot\|_{\Lambda^{p,\infty}_X(w)}\/$ is not a norm (except if  $\mu(X)=0\/$).

\item[(iv)]\ \ $\Lambda^{p,\infty}_X(w)\/$ is normable if and only if $w\in B_p\/$.
\end{enumerate}
\eteo

\bdem
That the conditions are sufficient has been already seen in  Theorems~\ref{te: II.5.1}, \ref{te: II.5.2},
and Remark~\ref{obs: 5.3}. To see the necessity, we use Corollary~\ref{cor: II.5.6}  and
 Lemma~\ref{lem: II.5.7} to conclude that we can substitute the space $X\/$ by  $\R\/$. Then (i)
follows from  \cite{CGS} and \cite{Lor}; (ii) from  \cite{CGS} and \cite{Saw2}; (iii) from
Remark~\ref{obs: 5.3}, and  (iv) follows from \cite{Sor}.$\qquad\qed$
\edem

We shall now study the normability of the sequence   Lorentz spaces $d(\Omega,p)\/$ and
$d^\infty(\Omega,p)\/$.  

\bteo
\label{te: II.5.9}
Let $\Omega=(\Omega_n)_{n=0}^\infty\subset[0,\infty)\/$.
\begin{enumerate}
\item[(i)]\ \ If $0<p<1,\;\|\cdot\|_{d(\Omega,p)}\/$ is a norm if and only if 
$\Omega=(\Omega_0,0,0,\dots)\/$.

\item[(ii)]\ \ If $1\le p<\infty,\;\|\cdot\|_{d(\Omega,p)}\/$ is a norm if and only if $\Omega\d\/$.
\end{enumerate}
\eteo

\bdem
The sufficiency in  (i) is obvious and in  (ii)   has already been proved in  Theorem~\ref{te:
II.5.1}. Let us see that if  $\|\cdot\|_{d(\Omega,p)}\/$ is a norm, necessarily 
$\Omega\d\/$. For $n\in\N$ and $t\in(0,1),\/$ let $f=(1,1,\dots,1,t,0,0,0,\dots)\/$  and let 
$g=(1,\dots,1,t,1,0,0\dots).\/$ 
Then, 
\begin{eqnarray*}
\|f\|_{d(\Omega,p)}&=&\|g\|_{d(\Omega,p)}=\big(\Omega_0+\dots+\Omega_n+t^p\Omega_{n+1}\big)^{1/p},
\\
\|f+g\|_{d(\Omega,p)}&=&\big(2^p\Omega_0+\dots+2^p\Omega_{n-1}+(1+t)^p(\Omega_n+
\Omega_{n+1})\big)^{1/p}.
\end{eqnarray*}
Form the inequality 
$\|f+g\|_{d(\Omega,p)}\le\|f\|_{d(\Omega,p)}+\|g\|_{d(\Omega,p)}\/$ it follows that 
$$
\Omega_{n+1}\le{2^p-(1+t)^p\over(1+t)^p-2^p t^p}\Omega_n,
$$ 
and letting  $t\/$ tend to  $1\/$
we obtain that  $\Omega_{n+1}\le\Omega_n\/$ and  (ii) is proved. 

Let us assume now that  $\|\cdot\|_{d(\Omega,p)}\/$ is a norm and  $p<1\/$. Then (Remark~\ref{obs:
5.5}) it coincides with the norm of the biassociate space. Then since  by  Theorem~\ref{te: II.4.16}
(i), $d({\Omega},p)'=d(\widetilde{\Omega},1)'$, where 
$\widetilde{\Omega}_0=\Omega_0^{1/p},\;\widetilde{\Omega}_n=\big(\sum_{k=0}^n
\Omega_k\big)^{1/p}-\big(\sum_{k=0}^{n-1}\Omega_k\big)^{1/p},\;n=1,2,\dots\/$, we obtain that
$\Vert\cdot\Vert_{d({\Omega},p)''}=\Vert\cdot\Vert_{d(\widetilde{\Omega},1)''}$. Applying this equality of
norms to 
$f=(1,t,0,0,\dots),\;t\in(0,1),\/$ we obtain that 
\begin{eqnarray*}
(\Omega_0+t^p\Omega_1)^{1/p}=\|f\|_{d(\Omega,p)}&=&\|f\|_{d(\widetilde{\Omega},1)^{\prime\prime}}
\le\|f\|_{d(\widetilde{\Omega},1)}\\
&=&\Omega_0^{1/p}+t\big((\Omega_0+\Omega_1)^{1/p}
-\Omega_0^{1/p}\big),
\end{eqnarray*}
and hence
$$
{(\Omega_0+t^p\Omega_1)^{1/p}-\Omega_o^{1/p}\over t}\le C<\infty. 
$$
If $\Omega_1\neq0\/$ the limit, when $t\to0\/$, of the left hand side is  $+\infty\/$, and thus 
$\Omega_1=0\/$ and,  since the sequence is decreasing, we finally obtain 
$\Omega_n=0,\;n\ge1.\qquad\qed$
\edem

Let us consider now the normability of  $d(\Omega,p)\/$.

\bteo
\label{te: II.5.10}
Let $\Omega=(\Omega_n)_{n=0}^\infty\subset[0,\infty)\/$ and let us denote
$W_n=\sum_{k=0}^n\Omega_k,\;n=0,1,2,\dots\/$
\begin{enumerate}
\item[(i)]\ \ If $0<p<1,\;d(\Omega,p)\/$ is normable if and only if  $\Omega\in\ell^1\/$.

\item[(ii)]\ \ $d(\Omega,1)\/$ is normable if and only if  
$$
{W_n\over n+1}\le C{W_m\over m+1},\qquad 0<m<n.
$$

\item[(iii)]\ \ If $1<p<\infty,\;d(\Omega,p)\/$ is normable if and only if 
$$
\sum_{k=0}^n{1\over W^{1/p}_k}\le C{n+1\over W^{1/p}_n}.\qquad n=0,1,2,\dots
$$
\end{enumerate}

That is, if $p<1,\;d(\Omega,p)\/$ is not  normable except in the trivial case $\Omega\in\ell^1\/$ 
and then 
$d(\Omega,p)=\ell^\infty\/$.
\eteo

\bdem (i) If $\Omega\in\ell^1,\;d(\Omega,p)=\ell^\infty\/$ and there is nothing to prove.
Conversely, if $d(\Omega,p)\/$ is normable,
$\|\cdot\|_{d(\Omega,p)^{\prime\prime}}\/$  is a norm in this space (Theorem~\ref{te: II.5.4}) which is
majorized by the norm of 
$d(\widetilde{\Omega},1)\/$, with
$\widetilde{\Omega}_0=\Omega_0^{1/p},\;\widetilde{\Omega}_n=\big(\sum_{k=0}^n\Omega_k\big)^{1/p}
-\big(\sum_{k=0}^{n-1}\Omega_k\big)^{1/p},\;n=1,2,\dots\/$
(see the proof of the previous theorem). We have then the inequality,
$$
\bigg(\sum_{n=0}^\infty g(n)^p\Omega_n\bigg)^{1/p}\le C\sum_{n=0}^\infty
g(n)\widetilde{\Omega}_n,\qquad g\d. 
$$
If $r=1/p>1\/$ the previous expression  is equivalent to 
$$
S=\sup_{g\d}{\sum_{n=0}^\infty g(n)\Omega_n\over\big(\sum_{n=0}^\infty
g(n)^r\widetilde{\Omega}_n\big)^{1/r}}<\infty. 
$$
By Theorem~\ref{te: I.6.11},
$$
S\approx\bigg(\int_0^\infty\bigg({W(t)\over\widetilde{W}(t)}\bigg)^{p/(1-p)}w(t)\,dt\bigg)^{1-p},
$$ 
with 
$w=\sum_{k=0}^\infty\Omega_k\chi_{[k,k+1)},\;W(t)=\displaystyle\int_0^t w(s)\,ds\/$ and analogously 
$\tilde{w}$ and $\widetilde{W}\/$. Using that $W\in \Delta_2\/$ (because $d(\Omega,p)\/$ is
quasi-normed), one can easily see that the above integrand is  comparable, in every  interval
$[n,n+1)\/$, to the function  ${w/ W}\/$. We have then
$$
\int_1^\infty {w(t)\over W(t)}\,dt=\log {W(\infty)\over W(1)}\lesssim S<\infty,
$$
which implies  $w\in L^1\/$ and  $\Omega\in\ell^1\/$.

(ii) and  (iii): Let $w(t)=\sum_{k=0}^\infty\Omega_k\chi_{[k,k+1)}(t),\;t>0\/$.
Then the condition of the statement implies  $w\in B_{p,\infty}\/$ (c.f. Theorem~\ref{te: I.6.5} and
Corollary~\ref{cor: I.6.18})  and by  Theorem~\ref{te: II.5.2}, $d(\Omega,p)=\Lambda^p_{\N^\ast}(w)\/$ is
normable. To see the converse, let 
$f\ge0\/$ be a decreasing sequence in 
$d(\Omega,p)\/$. If $n\ge1\/$, let us define $g_n=\sum_{j=-n}^n f_j\/$ where, for each
$j\in\N^\ast,\;f_j(k)=f(k+j)\chi_{\{k\ge -j\}}(k),\;k=0,1,2,\dots\/$ Then 
\bequ
\label{II.5.11}
\|g_n\|_{d(\Omega,p)}\le
C\sum_{j=-n}^n\|f_j\|_{d(\Omega,p)}\le(2n+1)C\|f\|_{d(\Omega,p)},\quad\forall n, 
\eequ
and on the other hand, we have  
$g_n\ge\big(f(0)+\dots+f(n)\big)\chi_{\{0,\dots,n\}}\/$ and hence
$\|g_n\|_{d(\Omega,p)}\ge\big(f(0)+\dots+f(n)\big)\big\|\chi_{\{0,\dots,n\}}\big\|_{d(\Omega,p)}
=\big(f(0)+\dots+f(n)\big)W_n^{1/p}\/$.
Combining this with  (\ref{II.5.11}) we obtain that
$$
A_d f(n)W_n^{1/p}\le 2C\|f\|_{d(\Omega,p)}=2C\|f\|_{\ell^p(\Omega)},\qquad n\in\N,\;f\d,
$$ 
which is equivalent to the boundedness of  $A_d:\ell^p_{\hbox{\decp
dec}}(\Omega)\to\ell^{p,\infty}(\Omega)\/$. Now, by Theorem~\ref{te: I.6.13}, the conditions of the
statement follow.$\qquad\qed$
\edem 

To finish this section, we characterize the normability in the weak  case.

\bteo
\label{te: II.5.12}
 Let $0<p<\infty,\;\Omega=(\Omega_n)_{n=0}^\infty\subset[0,\infty)\/$ and let us denote by 
 $W_n=\sum_{k=0}^n\Omega_k,\;n=0,1,2,\dots\/$ Then:
\begin{enumerate}
\item[(i)]\ \ $\|\cdot\|_{d^\infty(\Omega,p)}\/$ is  a norm if and only if  $\Omega=(\Omega_0,0,0,\dots).\/$

\item[(ii)]\ \ $d^\infty(\Omega,p)\/$ is normable if and only if 
$$
\sum_{k=0}^n{1\over W_k^{1/p}}\le C{n+1\over W_n^{1/p} },\qquad n=0,1,2,\dots
$$
\end{enumerate}
\eteo

\bdem
 (i) is a consequence of  Remark~\ref{obs: 5.3}. To see   (ii) it is enough to observe that the 
normability of  $d^\infty(\Omega,p)\/$ is equivalent to the boundedness of 
\bequ
\label{II.5.13}
A_d:\ell^{p,\infty}_{\hbox{\decp dec}}(\Omega)\to\ell^{p,\infty}(\Omega), 
\eequ
 since if 
(\ref{II.5.13}) holds, then $\|f\|=\|A_d f^\ast\|_{\ell^{p,\infty}(\Omega)}\/$ is a norm in 
$d^\infty(\Omega,p)\/$ equivalent to the original quasi-norm, and if  $d^\infty(\Omega,p)\/$ is
normable,  the same argument  used in the last part of the previous theorem 
(changing 
$\|f\|_{d(\Omega,p)}\/$ by
$\|f\|_{d^\infty(\Omega,p)})\/$ shows (\ref{II.5.13}). Applying now Theorem~\ref{te: I.6.16}
we conclude the proof.$\qquad\qed$
\edem

\section{Interpolation of operators}
\label{se:oper}

In this section $(X,\mu)$ and $(\bar{X},\bar{\mu})\/$ are $\sigma$-finite measure spaces, although this
will not be  necessary in the interpolation theorems  for the
spaces
$\lox{p,q}{X}{w}\/$. In some cases (Theorem~\ref{te: II.6.3}) no conditions on the weight 
$w\/$ will be required. In the more general result (Theorem~\ref{te: II.6.5}) we will assume
that the spaces involved are quasi-normed  ($W\in\Delta_2(X)\/$).
This  theorem is a generalization of  Marcinkiewicz theorem
adapted to the context of the   $\lox{p,q}{}{w}\/$ spaces. 

Finally, we shall see how some of the results of chapter~\ref{ch:characteristic} on boundedness of 
order continuous operators can be extended to the context of  Lorentz spaces.

We start by recalling some concepts connected with the real method of interpolation. In this section
we call functional lattice to any class $A\/$ formed by measurable functions
and defined by 
$$
A=\big\{f\,:\,\|f\|_A<\infty\big\},
$$ 
where $\|\cdot\|_A\/$ is a nonnegative functional that acts on measurable functions and satisfies
$\|rf\|_A=r\|f\|_A,\,r>0,\/$ and $\|f\|_A\le\|g\|_A\/$ if $|f(x)|\le|g(x)|\/$ a.e. $x\/$. In
particular the  Lorentz spaces $\lox{p,q}{X}{w}\/$ are of this type. In what follows
$A,B,A_0,A_1,B_0,B_1\/$ denote arbitrary functional lattices. We shall write $A\approx B\/$, if
$A=B\/$ and$\|\cdot\|_A\approx\|\cdot\|_B\/$. Let us now recall the definition of the $K$-functional  
associated to the pair $(A_0,A_1)\/$. For every function  $f\in A_0+A_1\/$ and 
$t>0\/$ 
$$ 
K(f,t,A_0,A_1)=\inf_{f=f_0+f_1}\big(\|f_0\|_{A_0}+t\|f_1\|_{A_1}\big),
$$ 
where the infimum extends over all possible decompositions $f=f_0+f_1,\;f_i\in A_i,\;i=0,1\/$.
We shall simply write
$K(t,f)\/$ if  it is clear which pair    $(A_0,A_1)\/$ we are working with. The main properties of 
$K\/$ are the following:

(i)\ \ \ \ $K(rf,t)=rK(f,t),\quad r>0\/$,

(ii)\ \ \ $|f|\le|g|\Rightarrow K(f,t)\le K(g,t)\/$,

(iii)\ \ $A_i\approx B_i,\;i=0,1\Rightarrow K(f,t,A_0,A_1)\approx K(g,t,B_0,B_1)\/$.

\noindent Properties (i) and (iii) are immediate while  (ii) can be seen in  \cite{KPS}. For each 
$0<\theta<1,\;0<q\le\infty,\/$ we define the following  \lq\lq norm\rq\rq\ in $A_0+A_1\/$:
$$
\|f\|_{(A_0,A_1)_{\theta,q}}=\bigg(\int_0^\infty\big(t^{-\theta}K(f,t,A_0,A_1)\big)^q\,{dt\over
t}\bigg)^{1/q},
$$ 
($\sup_t t^{-\theta}K(f,t,A_0,A_1)\/$ if 
$q=\infty\/$). It is then  natural to define the space
$$
(A_0,A_1)_{\theta,q}=\big\{f\in A_0+A_1\,:\,\|f\|_{(A_0,A_1)_{\theta,q}}<\infty\big\}. 
$$
We say that the operator $T\/$, defined in $A_0+A_1\/$,  and with values in $B_0+B_1\/$, is
quasi-additive if there exists $k>0\/$ such that
$$
\big|T(f+g)(x)\big|\le k\big(|Tf(x)|+|Tg(x)|\big)\quad\hbox{a.e.\ }x,
$$ 
for every pair of functions 
$f,g,f+g\in A_0+A_1\/$.

The fundamental result of this theory is the following.

\bteo
\label{te: II.6.1} Let $T\/$ be a quasi-additive operator defined in  $A_0+A_1\/$
and such that 
\begin{eqnarray*} 
&T&\,:\,A_0\longrightarrow B_0,\\ &T&\,:\,A_1\longrightarrow B_1.
\end{eqnarray*}
Then, for   
$0<\theta<1,\;0<q\le\infty,\/$ we have,
$$
T\,:\,(A_0,A_1)_{\theta,q}\longrightarrow(B_0,B_1)_{\theta,q}.
$$
\eteo

Our purpose now is to identify the space  $(A_0,A_1)_{\theta,q}\/$ or, equivalently, the functional 
$\|\cdot\|_{(A_0,A_1)_{\theta,q}}\/$ when 
$A_0,A_1\/$ are Lorentz spaces. As we shall see, $\|\cdot\|_{(A_0,A_1)_{\theta,q}}\/$ will be,
under appropriate conditions, equi\-valent to the  \lq\lq norm\rq\rq of a certain Lorentz space.

\bteo
\label{te: II.6.2}
 If $0<p<\infty\/$ and $f\in\M(X)\/$,
$$
K\big(f,t,\Lambda^p_X(w),L^\infty(X)\big)\approx K\big(f^\ast,t,L^p(w),L^\infty(w)\big),\qquad
t>0,
$$ 
with constants depending only on  $p\/$. In
particular, for  $0<\theta<1,\,0<q\le\infty\/$,
$$
\big(\Lambda^p_X(w),L^\infty(X)\big)_{\theta,q}\approx\Lambda^{\bar{p},q}_X(w)
$$ 
where,
$$\quad{1\over\bar{p}}={1-\theta\over p}.$$
\eteo

\bdem
 If $f=f_0+f_1\/$ with $f_0\in\Lambda^p_X(w),\,f_1\in L^\infty(X)\/$, we have that
$f^\ast(s)\le f_0^\ast(s)+f_1^\ast(0)=f_0^\ast(s)+\|f_1\|_{L^\infty(X)},\;s>0\/$. Hence,
$$
K\big(f^\ast,t,L^p(w),L^\infty(w)\big)\le 
\|f_0^\ast\|_{L^p(w)}+t\|f_1\|_{L^\infty(X)}=\|f_0\|_{\Lambda^p_X(w)}+t\|f_1\|_{L^\infty(X)}.
$$
Taking the infimum over all decompositions
$f=f_0+f_1\in\Lambda^p_X(w)+L^\infty(X)\/$ we obtain,
$$
K\big(f^\ast,t,L^p(w),L^\infty(w)\big)\le K\big(f,t,\Lambda^p_X(w),L^\infty(X)\big).
$$

To prove the converse inequality, if $f\in\M(X),\,t>0\/$ let $a=(f^\ast)^\ast_w(t^p)\/$ and let 
$$
f_0=\bigg(f-a{f\over|f|}\bigg)\chi_{\{|f|>a\}},\qquad f_1=f-f_0. 
$$
Then
$(f^\ast_0)^\ast_w=\big((f^\ast)^\ast_w-a\big)\chi_{[0,t^p)}\/$ while 
$f^\ast_1\le a\/$. Since $f=f_0+f_1\/$ we have that
\begin{eqnarray*} 
K\big(f,t,\loX,L^\infty(X)\big)&\le&\|f_0\|_\loX+t\|f_1\|_{L^\infty(X)}\\
&\le&\|f^\ast_0\|_{L^p(w)}+ta\\ &=& \big\|(f^\ast_0)^\ast_w\big\|_p+ta\\
&=&\bigg(\int_0^{t^p}\big((f^\ast)^\ast_w(s)-a\big)^p\,ds
\bigg)^{1/p}+
\bigg(\int_0^{t^p}a^p\,ds\bigg)^{1/p}\\
&\le& C_p\bigg(\int_0^{t^p}\big((f^\ast)^\ast_w(s)\big)^p\,ds\bigg)^{1/p}.
\end{eqnarray*} 
Since the last expression is equivalent to $K\big(f^\ast,t,L^p(w),L^\infty(w)\big)\/$ (see
\cite{BL}), the first part of the theorem is proved.

The second part is an immediate consequence of the previous one. To see this, observe that the \lq\lq
norm\rq\rq in
$\big(\Lambda^p_X(w),L^\infty(X)\big)_{\theta,q}\/$ is defined using the $K$-functional and we 
have 
$$
\|f\|_{\big(\Lambda^p_X(w),L^\infty(X)\big)_{\theta,q}}=\|f^\ast\|_{\big(L^p(w),
L^\infty(w)\big)_{\theta,q}},
$$
and the last \lq\lq norm\rq\rq is (see
\cite{BL}),
$$
\|f^\ast\|_{L^{\bar{p},q}(w)}=\|f\|_{\Lambda_X^{\bar{p},q}(w)}.\qquad\qed
$$
\edem 

A consequence of this result is the following interpolation theorem with 
$L^\infty(X)\/$. 

\bteo
\label{te: II.6.3}
If $0<p,\bar{p}<\infty\/$ and  $T\/$ is a quasi-additive operator in 
$\lox{p}{X}{w}+L^\infty(X)\/$ such that,
\begin{eqnarray*}
T\,&:&\,L^\infty(X)\longrightarrow L^\infty(X),\\
T\,&:&\,\lox{p}{X}{w}\longrightarrow\lox{\bar{p},\infty}{\bar{X}}{\bar{w}},
\end{eqnarray*}
then for 
$q,\bar{q}\in(0,\infty)\/$ satisfying  $q/p=\bar{q}/\bar{p}>1,\/$ we have
$$ 
T\,:\,\lox{q,r}{X}{w}\longrightarrow\lox{\bar{q},r}{\bar{X}}{\bar{w}},\qquad0<r\le\infty.
$$
\eteo

\bdem
The argument used in the first part of the previous theorem to prove the inequality 
$$
K\big(f^\ast,t,L^p(w),L^\infty(w)\big)\le K\big(f,t,\Lambda^p_X(w),L^\infty(X)\big),
$$ 
still works if we substitute the spaces 
$L^p(w)$ and $\Lambda^p_X(w)\/$ by $L^{\bar{p},\infty}(\bar{w})\/$ and 
$\Lambda^{\bar{p},\infty}_{\bar{X}}(\bar{w})\/$ respectively. Therefore, with
$1/\bar{q}=(1-\theta)/\bar{p}\/$, it follows   that (see \cite{BL}), 
\begin{eqnarray*}
\|Tf\|_{\lox{\bar{q},r}{\bar{X}}{\bar{w}}}&=&\big\|(Tf)^\ast\big\|_{L^{\bar{q},r}(\bar{w})}=
\big\|(Tf)^\ast\big\|_{ \big(L^{\bar{p},\infty}(\bar{w}),L^\infty(w)\big)_{\theta,r}}\\
&\le&\|Tf\|_{\big(\lox{\bar{p},\infty}{}{\bar{w}},L^\infty\big)_{\theta,r}}
\end{eqnarray*}
and the rest is a consequence of Theorem~\ref{te: II.6.1}
 and Theorem~\ref{te: II.6.2}.$\qquad\qed$
\edem

\bobs
\label{obs: 6.4}
Observe that in the last two results, it is not necessary that the 
Lorentz spaces involved are quasi-normed.
If we assume that these spaces are quasi-normed (that is $W\in\Delta_2\/$) we can  use the
reiteration theorem to obtain a more general result analogous to the  interpolation  theorem  of 
Marcinkiewicz.
\eobs

\bteo
\label{te: II.6.5}
 Let $0<p_i,q_i,\bar{p}_i,\bar{q}_i\le\infty,\;i=0,1,\/$ with $p_0\ne
p_1,\,\bar{p}_0\ne\bar{p}_1\/$ and let  $T\/$ be a quasi-additive operator defined in 
$\lox{p_0,q_0}{X}{w}+\lox{p_1,q_1}{X}{w}\/$ satisfying
\begin{eqnarray*}
 T\,&:&\,\lox{p_0,q_0}{X}{w}\longrightarrow\lox{\bar{p}_0,\bar{q}_0}{\bar{X}}
{\bar{w}},\\
T\,&:&\,\lox{p_1,q_1}{X}{w}\longrightarrow\lox{\bar{p}_1,\bar{q}_1}{\bar{X}}{\bar{w}}.
\end{eqnarray*}
Assume that $W\in\Delta_2(X)$ and $\bar{W}\in\Delta_2(\bar{X})\/$. Then, for
$0<\theta<1,\,0<r\le\infty\/$,
$$
T\,:\,\lox{p,r}{X}{w}\longrightarrow\lox{\bar{p},r}{\bar{X}}{\bar{w}},
$$ 
where
$$
{1\over p}={1-\theta\over p_0}+{\theta\over p_1},\quad{1\over
\bar{p}}={1-\theta\over\bar{p}_0}+{\theta\over\bar{p}_1}.
$$
\eteo

\bdem
Let $0<s<\min\{p_0,p_1\}\/$ and take  $\theta_0,\theta_1\in(0,1)\/$ such that
$1/p_i=(1-\theta_i)/s,\,i=0,1\/$. Then,  by Theorem~\ref{te: II.6.2} it follows,
$$
\big(\lox{p_0,q_0}{X}{w},\lox{p_1,q_1}{X}{w}\big)_{\theta,r}\approx\bigg(\big(\lox{s}{X}{w},
L^\infty(X)\big)_{\theta_0,q_0},\big(\lox{s}{X}{w},L^\infty(X)\big)_{\theta_1,q_1}
\bigg)_{\theta,r}.
$$
Since  $W\in\Delta_2(X)\/$ we have that  $\lox{s}{X}{w}\/$ is a quasi-Banach space
and we can apply the reiteration  theorem (Theorem~{3.11.5} in \cite{BL}) to identify the above
space:  
$$
\big(\lox{s}{X}{w},\lox{\infty}{X}{w}\big)_{\eta,r}
$$ 
with
$\eta=(1-\theta)\theta_0+\theta\theta_1\/$. Applying again Theorem~\ref{te: II.6.2}, we finally
obtain
$$
\big(\lox{p_0,q_0}{X}{w},\lox{p_1,q_1}{X}{w}\big)_{\theta,r}\approx\lox{p,r}{X}{w}
$$ 
since
$(1-\eta)/s=(1-\theta)/p_0+\theta/p_1\/$. The same argument works for 
$\lox{\bar{p_i},\bar{q_i}}{\bar{X}}{\bar{w}}\/$ and analogously,
$$
\big(\lox{\bar{p}_0,\bar{q}_0}{\bar{X}}{\bar{w}},\lox{\bar{p}_1,\bar{q}_1}{\bar{X}}
{\bar{w}}\big)_{\theta,r}\approx\lox{\bar{p},r}{\bar{X}}{\bar{w}}.
$$
The result   now follows from Theorem~\ref{te: II.6.1}. $\qquad\qed$
\edem 

\bobs
\label{obs: 6.6}
 J. Cerd\`a and J. Mart\'{\i}n (\cite{CM}) have proved, under some conditions on the weights  $w_0,w_1\/$
that, 
$$
K\big(f,t,\lox{p_0,r_0}{}{w_0},\lox{p_1,r_1}{}{w_1}\big)\approx
K\big(f^\ast,t,L^{p_0,r_0}(w_0),L^{p_1,r_1}(w_1)\big)
$$
 for $0<p_0,p_1,r_0,r_1\le\infty\/$. This allows them to obtain a more general interpolation
theorem including  the case
$$
T\,:\,\big(\lox{p_0,q_0}{X}{w_0},\lox{p_1,q_1}{X}{w_1}\big)\longrightarrow
\big(\lox{\bar{p}_0,\bar{q}_0}{\bar{X}}{\bar{w}_0},\lox{\bar{p}_1,\bar{q}_1}
{\bar{X}}{\bar{w}_1}\big),
$$
with $w_0\ne w_1,\bar{w}_0\ne\bar{w}_1\/$. 
\eobs

The following Marcinkiewicz interpolation type result is a  direct consequence of  Theorem~\ref{te:
II.6.5}.

\bcor
\label{cor: II.6.7}
 Let $0<p_0<p_1\le\infty,\,W\in\Delta_2(X)$, and $\bar{W}\in\Delta_2(\bar{X})\/$.
If $T\/$ is a quasi-additive operator in 
$\lox{p_0}{X}{w}+\lox{p_1}{X}{w}\/$ such that,
\begin{eqnarray*} 
T\,&:&\,\lox{p_0}{X}{w}\longrightarrow\lox{p_0,\infty}{\bar{X}}{\bar{w}},\\
T\,&:&\,\lox{p_1}{X}{w}\longrightarrow\lox{p_1,\infty}{\bar{X}}{\bar{w}},
\end{eqnarray*} 
then, for $p_0<p<p_1\/$,
$$
T\,:\,\lox{p}{X}{w}\longrightarrow\lox{p}{\bar{X}}{\bar{w}}.
$$ 
\ecor

Before ending this chapter, we shall extend some of the results of  
section~\ref{se I: characteristic functions}
about boundedness of order  continuous operators. 
The first result is a generalization of  Corollary~\ref{cor: I.2.14}.

\bteo
\label{te: II.6.8}
 Let $L\subset\M(X)\/$ be a regular class and let $T:L\to\M(\bar{X})\/$ be a 
sublinear order continuous operator. If $0<p_0\le1,\;p_0\le p_1<\infty\/$ and 
$w_1\in B_{p_1/p_0,\infty}\/$, we have 
\bequ
\label{II.6.9}
\|Tf\|_{\Lambda^{p_1}_{\bar X}(w_1)}\le C\|f\|_{\Lambda^{p_0}_X(w_0)},\qquad f\in L, 
\eequ
if and only if there exists $C_0<\infty\/$ such that
$$
\big\|T\chi_B\big\|_{\Lambda^{p_1}_{\bar X}(w_1)}\le
C_0\big\|\chi_B\big\|_{\Lambda^{p_0}_X(w_0)},\qquad \chi_B\in L.
$$
\eteo

\bdem
We shall assume the last inequality and we shall prove (\ref{II.6.9}). By monotonicity, it is
sufficient to prove it for a 
 simple function $f\ge0\/$. Proceeding as in the proof of  Theorem~\ref{te: I.2.13}, 
$$
\|Tf\|_{\Lambda^{p_1}_{\bar X}(w_1)}^{p_0}=\big\||Tf|^{p_0}\big\|_{\Lambda^{p_1/p_0}_{\bar
X}(w_1)}\le\bigg\|\int_0^\infty p_0
t^{p_0-1}\big|T\chi_{\{f>t\}}(\cdot)\big|^{p_0}\,dt\bigg\|_{\Lambda^{p_1/p_0}_{\bar X}(w_1)}. 
$$
But the condition $w_1\in B_{p_1/p_0,\infty}\/$ implies that
$\|\cdot\|_{\Lambda^{p_1/p_0}_{\bar X}(w_1)}\/$ is equivalent to a  Banach function norm
(Theorem~\ref{te: II.5.2}) and it follows,
\begin{eqnarray*}
\|Tf\|_{\Lambda^{p_1}_{\bar X}(w_1)}^{p_0}&\le& C_1\int_0^\infty  p_0
t^{p_0-1}\big\|\big|T\chi_{\{f>t\}}\big|^{p_0}\big\|_{\Lambda^{p_1/p_0}_{\bar X}(w_1)}\,dt\\
&=&C_1\int_0^\infty  p_0 t^{p_0-1}\big\|T\chi_{\{f>t\}}\big\|_{\Lambda^{p_1}_{\bar
X}(w_1)}^{p_0}\,dt\\ 
&\le& C_1C_0^{p_0}\int_0^\infty p_0
t^{p_0-1}\big\|\chi_{\{f>t\}}\big\|_{\Lambda^{p_0}_X(w_0)}^{p_0}\,dt\\ 
&=&C^{p_0}\int_0^\infty p_0
t^{p_0-1}W_0\big(\lambda_f(t)\big)\,dt\\ 
&=&C^{p_0}\|f\|_{\Lambda^{p_0}_X(w_0)}^{p_0}.\qquad\qed
\end{eqnarray*}

\edem

Using the same idea, one can easily prove an analogous result for the weak-type case:

\bteo
\label{te: II.6.10}
 Let $L\subset\M(X)\/$ be a  regular class and let  $T:L\to\M(\bar{X})\/$ be a 
sublinear order continuous operator. If $0<p_0\le1,\;0<p_1<\infty\/$ and
$w_1\in B_{p_1/p_0}\/$, we have 
$$
\|Tf\|_{\Lambda^{p_1,\infty}_{\bar X}(w_1)}\le C\|f\|_{\Lambda^{p_0}_X(w_0)},\qquad f\in L,
$$ 
if and only if, there exists $C_0<\infty\/$ such that
$$
\big\|T\chi_B\big\|_{\Lambda^{p_1,\infty}_{\bar X}(w_1)}\le
C_0\big\|\chi_B\big\|_{\Lambda^{p_0}_X(w_0)},\qquad \chi_B\in L.
$$
\eteo

Combining the previous results with the general interpolation theorem (Theorem~\ref{te: II.6.5})
we obtain a   generalization  of Stein and Weiss theorem on restricted weak-type operators (\cite{SW}). 

\bteo
\label{te: II.6.11}
 Let $0<p_0,p_1,q_0,q_1\le\infty,\;p_0\neq p_1,\;q_0\neq q_1\/$ and let us assume that 
$T\,:\,\big(\Lambda^{p_0}_X(w)+\Lambda^{p_1}_X(w)\big)\to\M(\bar{X})\/$ is a 
sublinear order continuous operator satisfying
\begin{eqnarray*}
\big\|T\chi_B\big\|_{\Lambda^{q_0,\infty}(\bar{w})}&\le& C_0\|\chi_B\|_{\Lambda^{p_0}(w)},\quad
B\subset X,\\
\big\|T\chi_B\big\|_{\Lambda^{q_1,\infty}(\bar{w})}&\le& C_1\|\chi_B\|_{\Lambda^{p_1}(w)},\quad
B\subset X.
\end{eqnarray*} 
Then, if $W,\bar{W}\in\Delta_2\/$, we have
$$
T\,:\,\Lambda^{p,r}_X(w)\longrightarrow \Lambda^{q,r}_{\bar{X}}(\bar{w}),\qquad 0<r\le\infty,
$$ 
if
$${1\over p}={1-\theta\over p_0}+{\theta\over p_1},\qquad {1\over q}={1-\theta\over
q_0}+{\theta\over q_1},\qquad0<\theta<1.
$$
\eteo

\bdem
 If $\bar{W}\in\Delta_2\/$ there exists $t>0\/$ such that $\bar{w}\in B_t\/$ (see
\cite{CGS}). Since the classes $B_p\/$ are increasing in $p\/$, there exists an index $r\in(0,1)\/$
with
$r<p_i,\;i=0,1\/$ and such that $\bar{w}\in B_{q_i/r},\;i=0,1\/$. Since
$\big\|\chi_B\big\|_{\Lambda^{p_i}}=C_{p_i,r}\big\|\chi_B\big\|_{\Lambda^{p_i,r}},\;B\subset X,\/$
we have 
$$
\big\|T\chi_B\big\|_{\Lambda^{q_i,\infty}(\bar{w})}\le
C_i^\prime\|\chi_B\|_{\Lambda^{p_i,r}(w)},\quad B\subset X,\quad i=0,1. 
$$
But
$\Lambda^{p_i,r}(w)=\Lambda^r(\tilde{w}_i)\/$ with $\tilde{w}_i=W^{r/p_i-1}w\/$ (Remark~\ref{obs:
II.2.6}) and by  Theorem~\ref{te: II.6.10} it follows that 
$$
T\,:\,\Lambda^{p_i,r}_X(w)\longrightarrow \Lambda^{q_i,\infty}_{\bar{X}}(\bar{w}),\qquad i=0,1.
$$
Applying then  Theorem~\ref{te: II.6.5} we obtain the result.$\qquad\qed$
\edem

\bobs
\label{obs: 6.12}

(i) Since
$\big\|\chi_B\big\|_{\Lambda^{p_i}}=C_{p_i,r}\big\|\chi_B\big\|_{\Lambda^{p_i,r}},\;B\subset X,\/$
the previous  theorem is still true if  we substitute the spaces
$\Lambda^{p_i}(w)\/$ by  $\Lambda^{p_i,r_i}(w)\/$, with $0<r_i\le\infty,\;i=0,1\/$.

(ii) An analogous result holds, without the hypothesis $W\in\Delta_2\/$, if we change the space
$\Lambda^{p_1}(w)\/$ by $L^\infty\/$, and with  $1/p=(1-\theta)/p_0\/$. This is true since, in this
case, one can use Theorem~\ref{te: II.6.2} (where this hypothesis is not needed) to identify the
interpolated space. 
\eobs

\chapter{The  Hardy-Littlewood maximal operator in weighted Lorentz spaces} 
\label{ch:maximal}
\pagestyle{myheadings}
\markboth{CHAPTER 3. MAXIMAL OPERATOR}{M.J.\ Carro, J.A.\ Raposo, and J.\ Soria}

\section{Introduction} 

In the previous chapter, we have introduced and studied the Lorentz spaces $\Lambda^{p,q}_X(w)\/$.
Our purpose in this new chapter is to study the boundedness of the Hardy-Littlewood maximal
operator of the type
\bequ
\label{bpp}
M\,:\,\lo\longrightarrow\lo,
\eequ
and its weak version, $M:\lo\to\loin\/$ (see below for the definition of these spaces). Our goal is to find
necessary and/or sufficient conditions on the weights $u\/$ and $w\/$ (in $\R^n\/$ and $\R^m\/$
respectively) to have (\ref{bpp}). We shall also obtain some positive result
for the nondiagonal case
\bequ
\label{bpq}
M\,:\,\loz\longrightarrow\loinu.
\eequ

If $w=1\/$,    (\ref{bpp}) is equivalent to $M:L^p(u)\longrightarrow L^p(u)\/$ and
this problem was completely solved by Muckenhoupt (\cite{Muck1}), who obtained the condition 
$u\in
A_p,\;1<p<\infty\/$ (see (\ref{Ap})). On the other hand, if $u=1\/$, the characterization of (\ref{bpp}) is
equivalent to the boundedness of the Hardy operator $A:\lpwdec\to L^p(w)\/$ and was
obtained ($p>1\/$) by Ari\~no and Muckenhoupt (\cite{AM1}). The condition, in this case, is 
$w\in B_p\/$ (see Definition~\ref{def: I.6.3}). 

The problem (\ref {bpq}) and the corresponding strong-type boundedness when $w_0,w_1\/$ are power
weights ($w_i(t)=t^{\alpha_i}\/$)  reduces to $M:L^{p_0,q_0}(u)\to
L^{p_1,q_1}(u)\/$ and was solved, for some cases, by Chung, Hunt, and Kurtz in \cite{CHK, HK} 
(see also \cite {L}). The complete solution for the general diagonal case (\ref{bpp}) has been
an open question so far, although there are some partial results due to Carro-Soria (\cite{CS3}) and 
Neugebauer (\cite{Ne2}). In this chapter, we prove (see Theorem~\ref{te:principal})  a complete
characterization of (\ref{bpp}) and give necessary and sufficient conditions to have the weak-type
boundedness. 

The chapter is organized as follows: in section~\ref{se: gr} we present some general results which
will be used in what follows, in section~\ref{se:dfcd}, we show that
(\ref{bpp}) holds if and only if there exists $q\in(0,p)\/$  such that
$$
{W\big(u\big(\bigcup_j Q_j\big)\big)\over W\big(u\big(\bigcup_j E_j\big)\big)}\le
C\max_j\bigg({|Q_j|\over|E_j|}\bigg)^q,
$$ 
for every finite family of cubes  and sets $(Q_j,E_j)_j\/$, with $E_j\subset Q_j\/$.

In many cases this condition can be simplified. Among other results, we shall see that although the
weight  $w\/$ in (\ref{bpp}) has to be in some class $B_p\/$, the weight  $u\/$  need not be even in 
$A_\infty\/$ (see (\ref{Ain})). In fact, we shall show that (\ref{bpp}) can be true with weights $u\/$  which
are not doubling.\index{weight!doubling}

In section~\ref{se:wtin}, we study the weak boundedness. For example, we prove that if $u_0=u_1=u\in A_1\/$, 
  (\ref{bpq}) is equivalent to the case $u_0=u_1=1\/$. Finally, in section~\ref{se:applic}, we
completely solve the problem
$$
M\,:\,L^{p,q}(u)\longrightarrow L^{r,s}(u),
$$
and also (\ref{bpq}), whenever $u_0=u_1=u\/$ is a power weight.
\medskip

\noindent
{\bf Notation:} Letters $u,u_0,u_1,\dots\/$, will be used to denote weights in $\R^n\/$. They will
be nonnegative measurable functions, not identically zero and integrable on sets of finite
measure. If $A\subset\R^n\/$ is measurable, we shall denote by $u(A)\/$ the measure of  $A\/$ in
the space
$X=\big(\R^n, u(x)dx\big)\/$; that is, 
$$
u(A)=\int_A u(x)\,dx.
$$

The distribution function of $f\in\M(\R^n)\/$ in this space,  will be denoted by  $\lambda_f^u\/$ 
and the decreasing  rearrangement by  $f^\ast_u\/$. $\Lambda^{p,q}_u(w)\/$ will be the  Lorentz
space 
\break 
$\Lambda^{p,q}_X(w)\/$. If  $u=1\/$ we simply write $\Lambda^{p,q}(w)\/$. Hence, the \lq\lq
norm\rq\rq\ of $f\/$ in the weighted Lorentz space  $\Lambda^p_u(w)\/$ will be given by\index{$\lo$} 
$$
\|f\|_{\lo}=\bigg(\int_0^\infty\big(f^\ast_u(t)\big)^p w(t)\,dt\bigg)^{1/p},
$$
and, in the weak space, \index{$\loin$} 
$$
\|f\|_{\loin}=\sup_{t>0}W^{1/p}(t)f^\ast_u(t)=\sup_{t>0}t\,W^{1/p}\big(\lambda_f^u(t)\big).
$$
Let us also recall that  $w\/$ is a weight in  $\R^+\/$  with support in  $[0,u(\R^n)]\/$ and that 
$$
0\le W(t)=\int_0^t w(s)\,ds<\infty,\;t>0.
$$
Finally, the    Hardy-Littlewood maximal  operator $M\/$\index{Hardy-Littlewood maximal operator} is
defined by 
$$
Mf(x)=\sup_{x\in Q}{1\over|Q|}\int_Q|f(y)|\,dy,\qquad x\in\R^n,\;f\in\M(\R^n),
$$
where the supremum extends over all cubes $Q\/$, with sides parallel to the axes.

\section{Some general results}
\label{se: gr}

It is known that the decreasing rearrangement of $Mf\/$,  with respect to the Lebesgue measure, is
equivalent (\cite{BS}) to the function  $f^{\ast\ast}\/$:
\bequ
\label{rmax}
(Mf)^\ast(t)\approx f^{\ast\ast}(t)=Af^\ast(t),\qquad t>0.
\eequ
Since every decreasing and positive function in  $\R^+\/$ is equal  a.e. to the decreasing
rearrangement  of a measurable function in  $\R^n\/$, we deduce that the boundedness of 
$M:\llo\to\lloin\/$ (resp.
$M:\llo\to\llo\/$) is equivalent to the boundedness of  $A:L^p_{\hbox{\decp dec}}(w)\to
L^{p,\infty}(w)\/$ (resp. $A:L^p_{\hbox{\decp dec}}(w)\to L^p(w)\/$). Therefore 
(see section~\ref{se I: Hardy
operator}) a necessary and sufficient condition is $w\in B_{p,\infty}\/$ (resp.
$w\in B_p\/$). On the other hand, when  $w=1\/$, the boundedness 
$M:\Lambda^p_u(w)\to\Lambda^{p,\infty}_u(w)\/$  is equivalent to  $M:L^p(u)\to L^{p,\infty}(u)\/$
which is known to be   $u\in A_p\/$ (for  $p\ge1\/$), the Muckenhoupt class of weights  (\cite {Muck1,
MW, CF}). This motivates the following definition:

\bdefi
If $0<p<\infty\/$,  we write  
$w\in B_p(u)\/$\index{$B_p(u)$} (respectively  $w\in B_{p,\infty}(u)\/$)\index{$B_{p,\infty}(u)$} if the
boundedness 
$M:\lo\to\lo\/$ (resp. $M:\lo\to\loin\/$) holds. We shall also write $u\in A_p(w)\/$\index{$A_p(w)$} if 
$M:\lo\to\loin\/$ holds.
\edefi

That is, $u\in A_p(w)\Leftrightarrow w\in B_{p,\infty}(u)\/$. 
It is clear then that  $A_p(1)=A_p,\;1\le p<\infty,\/$ and that  
$B_p(1)=B_p,\;B_{p,\infty}(1)=B_{p,\infty},\;0<p<\infty\/$. On the other hand, since
$\Lambda^p\subset\Lambda^{p,\infty}\/$, we always have  that $B_p(u)\subset
B_{p,\infty}(u),\;0<p<\infty\/$. With this terminology, our purpose will be to identify the classes 
$B_p(u)\/$ and 
$B_{p,\infty}(u)\/$, or equivalently to identify  the classes  $A_p(w)\/$. 

Using (\ref{rmax}), the nondiagonal case (with $u=1\/$) is also very simple: the boundedness of 
$M:\lox{p_0,q_0}{}{w_0}\to\lox{p_1,q_1}{}{w_1}\/$ is equivalent to the boundedness of the Hardy
operator 
$A:L^{p_0,q_0}_{\hbox{\decp dec}} (w_0)\to L^{p_1,q_1}(w_1)\/$, which is known in most of the
cases.  

The natural question now is to see if there exists an analogue to (\ref{rmax}) 
when the rearrangement  $f^\ast\/$ is taken with respect to a weight  $u\ne1\/$; that is, if 
$(Mf)^\ast_u(t)\approx Af^\ast_u(t)\/$ or $(Mf)^\ast_u(t)\le C Af^\ast_u(t)\/$, when  $u\/$ is a
weight in 
$\R^n\/$ satisfying certain conditions. However, it was shown in \cite{CS3} (see also 
\cite{ LN1, LN2}) that  $(Mf)^\ast_u\approx Af^\ast_u\/$  if and only if   $u\approx 1\/$, and 
$(Mf)^\ast_u\le C\big(Af^{\ast p}_u\big)^{1/p}\/$  if and only if $M:L^p(u)\to
L^{p,\infty}(u),\;p>1\/$ (equivalently if  $u\in A_p\/$). The following theorem generalizes this
result, characterizing 
$M:\lo\to\loin\/$ in terms of an expression like  (\ref{rmax}).
\bteo
\label{adpp}
If $0<p<\infty\/$,  $M:\lo\to\loin\/$ is bounded if and only if
$$
(Mf)^\ast_u(t)\le C\bigg({1\over W(t)}\int_0^t\big(f^\ast_u\big)^p(s) w(s)\,ds\bigg)^{1/p},
\qquad t>0,\;f\in\M(\Rn).
$$
\eteo

\bdem
The inequality of the statement implies
$$
W^{1/p}(t)(Mf)^\ast_u(t)\le C\bigg(\int_0^t\big(f^\ast_u\big)^p (s) w(s)\,ds\bigg)^{1/p}
\le C\|f\|_{\lo},\qquad t>0,
$$
or, equivalently, $\|Mf\|_{\loin}\le C\|f\|_{\lo}\/$. This proves that this condition is sufficient. 
On the other hand, let us assume that 
$M:\lo\to\loin\/$ is bounded and let  $f\in\lo\/$. If  $f=f_0+f_1\/$,  with
$f_1\in L^\infty\/$, we have, for every  $t>0\/$,
\begin{eqnarray*}
(Mf)^\ast_u(t)&\le&(Mf_0)^\ast_u(t)+(Mf_1)^\ast_u(0)
\\
&\le& W^{-1/p}(t)\|Mf_0\|_{\loin}+\|Mf_1\|_{L^\infty(u)}
\\
&\le& CW^{-1/p}(t)\big(\|f_0\|_{\lo}+W^{1/p}(t)\|f_1\|_{L^\infty(u)}\big),
\end{eqnarray*}
and taking the infimum over all decompositions $f=f_1+f_0\/$ we obtain
\bequ
\label{rmaxu}
(Mf)^\ast_u(t)\le CW^{-1/p}(t)K\big(f,W^{1/p}(t),\lo,L^\infty(u)\big),\quad t>0. 
\eequ
By Theorem~\ref{te: II.6.2}, $K\big(f^\ast_u,W^{1/p}(t),L^p(w),L^\infty(w)\big)\/$  
is equivalent to the   $K$-functional  
  and therefore (see \cite {BL}) 
equivalent to 
$$
\bigg(\int_0^{W(t)}\big((f^\ast_u)^\ast_w\big)^p(s)\,ds\bigg)^{1/p}=
\bigg(\int_0^t(f^\ast_u)^p (s) w(s)\,ds\bigg)^{1/p}.
$$
This and (\ref{rmaxu})  give us the inequality we are looking for. $\qquad\qed$ 
\edem

\bobs
The same argument works for the nondiagonal case and hence (\ref{bpq})   is
equivalent, if $0<p_0,p_1<\infty\/$, to the inequality
$$
(Mf)^\ast_{u_1}(t)\le C\bigg({1\over
W_1^{p_0/p_1}(t)}\int_0^{W_0(t)}\big((f^\ast_{u_0})^\ast_{w_0}\big)^{p_0}(s)\,ds
\bigg)^{1/p_0},\quad
t>0.
$$ 
\eobs

Next result gives a necessary and sufficient condition to have the weak-type inequality,
$$
\|Mf\|_{\loinu}\le C\|f\|_{\loz}
$$
on characteristic functions $f=\chi_E,\;E\subset\Rn\/$.

\bteo
\label{adfc} 
Let $0<p_0,p_1<\infty\/$. Then, 
\bequ
\label{Mcar}
\big\|M\chi_E\big\|_{\loinu}\le C\big\|\chi_E\big\|_{\loz},\qquad E\subset\Rn\hbox{\ measurable,}
\eequ
if and only if, for every finite family of cubes  $(Q_j)_{j=1}^J\/$ and every family of measurable
sets $(E_j)_{j=1}^J\/$, with
$E_j\subset Q_j,$ for every $ j,\/$ we have that 
\bequ
\label{dQjEj}
\frac{W_1^{1/p_1}\big(u_1\big(\bigcup_j Q_j\big)\big)}{W_0^{1/p_0}\big(u_0\big(\bigcup_j
E_j\big)\big)}\le C\max_j\frac{|Q_j|}{|E_j|}. 
\eequ
\eteo

\bdem
To prove the necessary condition, let us consider $f=\chi_E\/$ with $E=\bigcup_j E_j\/$. If
$t>0\/$ is such that ${1/ t}>\max_j{|Q_j|\over|E_j|}\/$, then
$$
{1\over|Q_j|}\int_{Q_j}f(x)\,dx\ge{|E_j|\over|Q_j|}>t,
$$ 
and we have that  $Q_j\subset\{Mf>t\}\/$. 
This holds for every  $j=1,\dots,J\/$ and hence, $\bigcup_j Q_j\subset\{Mf>t\}\/$. Therefore,
$$
tW_1^{1/p_1}\big(u_1\big(\bigcup_j Q_j\big)\big)\le tW_1^{1/p_1}\big(\lambda^{u_1}_{Mf}(t)\big)
\le C\|f\|_{\loz}=CW_0^{1/p_0}\big(u_0\big(\bigcup_j E_j\big)\big),
$$
and consequently
$$
{W_1^{1/p_1}\big(u_1\big(\bigcup_j Q_j\big)\big)\over W_0^{1/p_0}\big(u_0\big(\bigcup_j E_j\big)
\big)}\le {C\over t}.
$$
Since  ${1/ t}>\max_j{|Q_j|\over|E_j|}\/$ is arbitrary, we have shown (\ref{dQjEj}).

Conversely, let us assume that  (\ref{dQjEj}) holds. Then, this inequality also holds
if the families  $(Q_j)_j,\;(E_j)_j\/$ are countable nonfinite. If $f=\chi_E\in\M(\Rn)\/$ and
$t>0\/$, for every $x\in\{Mf>t\}\/$, there exists a cube  $Q\/$, with $x\in Q\/$ and such that
$ \displaystyle\int_Q f(x)\,dx={|E\cap Q|}>t\/$. By definition we have that 
$Mf(y)>t$, for every $ y\in Q\/$ and hence, $Q\subset\{Mf>t\}\/$. Using again this argument, we
obtain a countable  family of cubes
$(Q_j)_j\/$ such that  $\{Mf>t\}=\bigcup_j Q_j\/$ and satisfying (with $E_j=E\cap Q_j\/$)
${|E_j|/|Q_j|}>t,\;j=1,2,\dots\/$ Then, 
$t\le\inf_j{|E_j|\over|Q_j|}=\big(\sup_j{|Q_j|\over|E_j|}\big)^{-1}\/$ and we have, applying 
 (\ref{dQjEj}),
\begin{eqnarray*}
tW_1^{1/p_1}\big(\lambda^{u_1}_{Mf}(t)\big)&\le&\displaystyle{W_1^{1/p_1}\big(u_1\big(\bigcup_j Q_j\big)\big)
\over\sup_j{|Q_j|\over|E_j|}}\\
&\le& C W_0^{1/p_0}\big(u_0\big(\bigcup_j E_j\big)\big)\leq C
W_0^{1/p_0}\big(u_0(E)\big).
\end{eqnarray*}
Taking now the supremum in  $t\/$ we obtain that 
$$
\|Mf\|_{\loinu}\le C W_0^{1/p_0}\big(u_0(E)\big)=C\|f\|_{\loz}.\qquad\qed
$$ 
\edem 

\bobs
\label{ob:cd}
In the condition  (\ref{dQjEj}),  we can assume that the sets $(E_j)_j\/$ 
are disjoint. In fact, for every finite family of cubes $(Q_j)_j\/$ there exists a subfamily 
$(Q_{j_k})_k\/$ of disjoint cubes such that $\bigcup_j Q_j\subset\bigcup_k Q^\ast_{j_k}\/$,
where every $Q^\ast_{j_k}\/$ is a dilation of  $Q_{j_k}\/$ with side three times bigger (see for
example \cite{St}). The sets  $(E_{j_k})_k\/$ are then disjoint and, if the condition holds, we
have

\begin{eqnarray*}
W_1^{1/p_1}\big(u_1\big(\bigcup_j Q_j\big)\big)\over W_0^{1/p_0}\big(u_0\big(\bigcup_j E_j\big)
\big) &\le \displaystyle{W_1^{1/p_1}\big(u_1\big(\bigcup_k Q_{j_k}^\ast\big)\big)\over
W_0^{1/p_0}\big(u_0\big(\bigcup_k E_{j_k}\big)\big)}
\le C\max_k{\big|Q_{j_k}^\ast\big|\over\big|E_{j_k}\big|} \\
&\le  C C_n\max_k{\big|Q_{j_k}\big|\over\big|E_{j_k}\big|}\le C C_n\max_j{|Q_j|\over|E_j|}.
\end{eqnarray*}
\eobs

\bcor
\label{cor: cnad}
If $M:\loz\to\loinu,\;0<p_0,p_1<\infty\/$, then
$$
{W_1^{1/p_1}\big(u_1(Q)\big)\over|Q|}\le C{W_0^{1/p_0}\big(u_0(E)\big)\over|E|},\quad E\subset Q,
$$
for every cube  $Q\subset\Rn\/$. In particular $W_0(t)>0,\;t>0\/$, and $u_0(x)>0\/$ 
a.e. $x\in\Rn.\/$
\ecor

The  following two propositions are a consequence of these results.

\bpro
\label{pro: cnu_0}
Let $0<p_0,p_1<\infty\/$ and let us assume that  $M:\loz\to\loinu\/$, then $u_0\notin L^1(\R^n)$.
\epro

\bdem
By Corollary~\ref{cor: cnad} we have that, for every cube  $Q\/$,
\bequ
\label{dEQ}
\frac{|E|}{|Q|}\le C\frac{W_0^{1/p_0} \big(u_0(E)\big)}{W_1^{1/p_1} \big(u_1(Q)\big)},
\quad E\subset Q.
\eequ
Let us assume that $u_0(\Rn)<\infty\/$. Then, 
if  $a\in\big(0,u_1(\Rn)\big)\/$ and $b\in(0,1)\/$ are such that 
$$
C{W_0^{1/p_0}\big(b u_0(\Rn)\big)\over W_1^{1/p_1}(a)}<5^{-n}, 
$$
with $C$ as in (\ref{dEQ}), we have that, if $u_1(Q)\ge a\/$,  the inequality ${u_0(E)/ u_0(Q)}\le b\/$
implies
$$
{|E|\over|Q|}\le C{W_0^{1/p_0}\big(u_0(E)\big)\over W_1^{1/p_1}\big(u_1(Q)\big)}\le 
C{W_0^{1/p_0}\big(b u_0(\Rn)\big)\over W_1^{1/p_1}(a)}<5^{-n};
$$
that is, if $Q\/$ is an arbitrary cube with  $u_1(Q)\ge a\/$,
\bequ
\label{dEQu_0}
E\subset Q,\;|E|\ge5^{-n}|Q|\Rightarrow u_0(E)>b u_0(Q).
\eequ

Let now  $Q_0\/$ be a cube with  $u_1(3Q_0)\ge a\/$ (for each cube  $Q\/$,  $kQ$ 
denotes another cube with the same center $Q\/$ and side $k\/$ times bigger than  $Q\/$). Let
$\tilde{Q}\subset3Q_0\/$ be a cube whose interior is disjoint with  $Q_0\/$ and such that
$|\tilde{Q}|=|Q_0|\/$. Then,  $5\tilde{Q}\supset 3Q_0\/$ and hence, $u_1(5\tilde{Q})\ge a\/$, and we
have by (\ref{dEQu_0})
$$
u_0(\tilde{Q})>b u_0(5\tilde{Q})\ge b u_0(Q_0).
$$
Therefore, $u_0(3Q_0)\ge u_0(Q_0)+u_0(\tilde{Q})\ge(1+b)u_0(Q_0)=\alpha u_0(Q_0)\/$, 
where $\alpha=1+b>1\/$. By the same argument, $u_0(9Q_0)\ge\alpha u_0(3Q_0)\ge\alpha^2 u_0(Q_0)\/$
and, in general, $u_0(3^n Q_0)\ge\alpha^n u_0(Q_0)\/$. Since $\lim_n u_0(3^n Q_0)=u_0(\Rn)\/$
and 
$\alpha>1\/$ we have that $u_0(\Rn)=\infty\/$ (observe that $u_0(Q_0)>0\/$ by  (\ref{dEQ}))
contradicting the initial assumption. $\qed$
\edem

\bpro
Let $0<p_0,p_1<\infty\/$. If $M:\lox{p_0}{u}{w}\to\loinx{p_1}{u}{w}\/$, 
then $p_1\le p_0\/$. If,  in addition, $w\notin L^1(\R^+)$, then $p_1=p_0\/$.
\epro

\bdem 
By Corollary~\ref{cor: cnad}, we have  that 
$$
W^{1/p_1-1/p_0}\big(u(Q)\big)\le C,
$$
for every cube  $Q\subset\Rn\/$, and by  Proposition~\ref{pro: cnu_0} 
$$
W^{1/p_1-1/p_0}(r)\le C,\quad r>0.
$$
Since $\lim_{t\to0}W(t)=0\/$, we obtain that $p_1\le p_0\/$. If $w\notin L^1(\R^+)$, 
$\lim_{t\to\infty}W(t)=\infty\/$ and the previous inequality holds  only if $p_0=p_1.\qquad\qed$
\edem

\bpro
\label{cnalpq}
 Let $0<p,q,r<\infty\/$. If $M:\lox{p,q}{u}{w}\to\lox{p,r}{u}{w}\/$, then
 $r\ge q\/$.
\epro

\bdem
 $|f|\le Mf\/$ for every $f\in\M(\Rn)\/$ and the hypothesis implies that
$\lox{p,q}{u}{w}\subset\lox{p,r}{u}{w}\/$. Since $\big(\Rn,u(x)dx\big)\/$ and  
$\big(\R^m,w(t)dt\big)\/$ are nonatomic spaces, and 
$\|f\|_{\lox{p,q}{u}{w}}=\|f^\ast_u\|_{L^{p,q}(w)}\/$, the previous embedding implies that 
$$
\bigg(\int_0^b g^r(t) t^{r/p-1}\,dt\bigg)^{1/r}\le C\bigg(\int_0^b g^q(t)
t^{q/p-1}\,dt\bigg)^{1/q},\qquad g\d,
$$ 
with $b=W\big(u(\R^n)\big)\/$. Equivalently,
$$
\sup_{g\d}{\displaystyle\int_0^b g(t)t^{r/p-1}\,dt\over\bigg(\displaystyle\int_0^b
g(t)^{q/r}t^{q/p-1}\,dt\bigg)^{r/q}}<\infty.
$$ 
Now, by  Theorem~\ref{te: I.5.7} in \cite{Saw2}, this supremum is finite if  $r<q.\quad\qed$
\edem

The following result is a consequence of the interpolation theorems developed in section~\ref{se I:
characteristic functions} of the previous chapter.

\bteo
\label{te:cnadfc}
Let $0<p_0,p_1<\infty\/$ and let us assume that  $W_1\in\Delta_2\/$. Then, the boundedness
$M:\Lambda^{p_0}_{u_0}(w_0)\to\Lambda^{p_1,\infty}_{u_1}(w_1)\/$ on characteristic  functions:
$$
\big\|M\chi_E\big\|_{\Lambda^{p_1,\infty}_{u_1}(w_1)}\le
C\|\chi_E\|_{\Lambda^{p_0}_{u_0}(w_0)},\qquad E\subset\Rn,
$$ 
implies
$$
M:\Lambda^{q_0,r}_{u_0}(w_0)\longrightarrow\Lambda^{q_1,r}_{u_1}(w_1),\qquad0<r\le\infty,
$$
for  $p_i<q_i<\infty,\;i=0,1,\;p_1/p_0=q_1/q_0\/$. In particular, we have  that
$M:\lox{q_0}{u_0}{w_0}\to\lox{q_1}{u_1}{w_1}\/$ if $p_1\ge p_0\/$.
\eteo

\bdem 
Since $M:L^\infty\to L^\infty\/$, the statement is an immediate consequence of Theorem~\ref{te: II.6.11}
(see also Remark~\ref{obs: 6.12}). $\qquad\qed$
\edem

Combining Theorem~\ref{adfc}  with the results of section~\ref{se:oper} of the previous chapter, on 
boundedness of continuous order operators on  Lorentz spaces, we obtain the following characterization for
the weak-type boundedness  in the case
$0<p_0\le1$ and $w_1\in B_{p_1/p_0}\/$.

\bteo
\label{cnsadnd}
 If $0<p_0\le1,\;0<p_1<\infty\/$, and  $w_1\in B_{p_1/p_0}\/$ we have the boundedness
$M:\loz\to\loinu\/$ if and only if,  for every finite family of cubes $(Q_j)_{j=1}^J\/$ and every
family of measurable sets 
$(E_j)_{j=1}^J\/$ with $E_j\subset Q_j,$ for every $ j,\/$ we have that 
$$
{W_1^{1/p_1}\big(u_1\big(\bigcup_j Q_j\big)\big)\over W_0^{1/p_0}\big(u_0\big(\bigcup_j E_j\big)
\big)}\le C\max_j{|Q_j|\over|E_j|}.
$$
\eteo

\bdem
 Since $M\/$ is an order continuous operator (Definition~\ref{def: I.2.6}) on the regular class $L=\loz\/$,
we have, by  Theorem~\ref{te: II.6.10},  that the boundedness   is equivalent 
to (\ref{Mcar}) and thus,  condition (\ref{dQjEj}) is necessary
and sufficient.$\qquad\qed$
\edem

\section{Strong-type boundedness in the diagonal case}
\label{se:dfcd}

In this section we shall give a characterization of the  $B_p(u)\/$ class in the more general case;
that is, we shall obtain a necessary and sufficient condition to have the boundedness
$$
M\,:\,\lo\longrightarrow\lo.
$$
To this end, we first need the two following technical lemmae:

\blem
\label{le:Bqp}
 Let $0<p<\infty\/$ and let us assume that, for every cube $Q\subset\Rn\/$, 
$$
{W\big(u(Q)\big)\over|Q|^p}\le C{W\big(u(E)\big)\over|E|^p},\qquad E\subset Q,
$$
with $C\/$ independent of $Q\/$. Then,  $w\in B_q\/$ for every $q>p\/$. In particular 
$W\in\Delta_2\/$.
\elem

\bdem 
The measure $u(x)dx\/$ is $\sigma$-finite and nonatomic and, therefore, 
\break $\big(\Rn,u(x)dx\big)\/$ is a resonant measure space. In these spaces, it holds that 
$$
\int_0^\infty f^\ast(s) g^\ast(s)\,ds=\sup_{h^\ast=g^\ast}\int f(x)h(x)u(x)\,dx
$$
(see \cite{BS}) for every measurable functions $f,g\/$. In our case, we have that, 
for $0<t<u(Q)\/$,

\begin{eqnarray*}
\big(u^{-1}\chi_Q\big)^{\ast\ast}_u(t)&=&{1\over t}\int_0^t\big(u^{-1}\chi_Q\big)^\ast_u(s)\,ds\\
&=&{1\over t}\sup\bigg\{\int_Qu^{-1}(x)\chi_E(x)\, u(x)\,dx\,:\,u(E)=t,\;E\subset Q\bigg\}\\
&=&{1\over t}\sup\big\{|E|\,:\,u(E)=t,\;E\subset Q\big\}\\
&\le& C{|Q|\over W^{1/p}\big(u(Q)\big)}
{W^{1/p}(t)\over t}. 
\end{eqnarray*}

Since the function $\big(u^{-1}\chi_Q\big)^{\ast\ast}_u\/$ 
is decreasing and,
$$
\big(u^{-1}\chi_Q\big)^{\ast\ast}_u(t)\ge\big(u^{-1}\chi_Q\big)^{\ast\ast}_u\big(u(Q)\big)
={|Q|\over
u(Q)}, 
$$ 
we obtain that
$$
{W^{1/p}\big(u(Q)\big)\over u(Q)}\le C{W^{1/p}(t)\over t},\qquad 0<t<u(Q),
$$
which is equivalent to  ${W^{1/p}(r)\over r}\le C{W^{1/p}(t)\over t},\;0<t<r<\infty\/$ and  
hence (see, for example \cite{Sor}) $w\in\bigcup_{q>p}B_q.\qquad\qed$
\edem

We shall also need the following result due to  Hunt-Kurtz (\cite{HK}).

\blem
\label{HK}
 If $t>0,\;E\subset\Rn,\/$ and we denote by 
$E_t=\big\{M\chi_E>t\big\}\/$, there exists a constant $\alpha>1\/$, depending only on the dimension,
such that
$$
\big(E_t\big)_s\supset E_{\alpha ts},\qquad s,t\in(0,1).
$$
\elem

It is known that the boundedness of  $M:L^p(u)\to L^p(u),\;p>1,\/$ implies $M:L^{p-\epsilon}(u)
\to L^{p-\epsilon}(u)\/$ for some $\epsilon>0\/$. Analogously,  $M:\Lambda^p(w)\to\Lambda^p(w)\/$
is equivalent to  $w\in B_p\/$, that implies  $w\in B_{p-\epsilon}\/$. 
As we are going to see next,
this also holds for the general case $M:\lo\to\lo\/$; that is, if $w\in
B_p(u),\;0<p<\infty,\/$ there exists $\epsilon>0\/$ such that $w\in B_{p-\epsilon}(u)\/$. 
In fact, in our next theorem, we prove a stronger result: if $M:\lo\to\lo\/$ on characteristics
functions, then 
$w\in B_{p-\epsilon}\/$. Part of its proof follows the same patterns developed in Theorem~2 of
\cite{HK}.

\bteo
\label{apfciaq}
 Let $0<p,r<\infty\/$ and let us assume that
$$
\big\|M\chi_E\big\|_{\Lambda^{p,r}_u(w)}\le C\|\chi_E\|_{\Lambda^p_u(w)},\qquad E\subset\Rn.
$$
Then, there exists $q\in(0,p)\/$ such that  $M:\lox{q}{u}{w}\to\lox{q}{u}{w}\/$.
\eteo

\bdem
First, we shall prove that the inequality of the statement is also true if we substitute
the indices $p,r$ by some others $\tilde p$ and $\tilde r$ with $\tilde p<p$ and $\tilde r<r$. 
 To this end, let us observe that the hypothesis of the theorem is equivalent (c.f.
Proposition~\ref{pro: II.2.5}) to the inequality
\bequ
\label{u(E_t)}
\int_0^1 t^{r-1}W^{r/p}\big(u(E_t)\big)\,dt\le BW^{r/p}\big(u(E)\big),\qquad E\subset\Rn,
\eequ
with $E_t,\;t>0,\/$ defined as in the previous lemma, and   $B\/$ is a constant independent of 
$E\/$. Let
$\alpha>1\/$, the constant, depending only on the dimension, of Lemma~\ref{HK}. We shall prove that,
for every 
$n=0,1,2,\dots\/$, the  inequality 
\bequ
\label{u(E_t)2}
\int_0^1 t^{r-1}W^{r/p}\big(u(E_t)\big){1\over n!}\log^n {1\over t}\,dt\le B(B\alpha^r)^n
 W^{r/p}\big(u(E)\big),\qquad E\subset\Rn, 
\eequ
holds.  If $n=0\/$, (\ref {u(E_t)2}) is  (\ref{u(E_t)}). By induction, we only have to show that 
(\ref {u(E_t)2})  implies an analogous inequality with $n+1\/$ instead of  $n\/$. To see this, and 
since
$s\in(0,1)\/$, we apply (\ref {u(E_t)2}) to the set $E=E_s\/$, obtaining
\bequ
\label{u((E_s)_t)}
\int_0^1 t^{r-1}W^{r/p}\big(u\big((E_s)_t\big)\big){1\over n!}\log^n {1\over t}\,dt
\le B(B\alpha^r)^n W^{r/p}\big(u(E_s)\big).
\eequ
By Lemma~\ref{HK}, the left hand side of (\ref{u((E_s)_t)}) is greater than or equal to
\begin{eqnarray*}
&&\int_0^{1/\alpha}t^{r-1}W^{r/p}\big(u(E_{\alpha st})\big){1\over n!}\log^n{1\over t}\,dt\\
=(\alpha s)^{-r}&&\int_0^s x^{r-1}W^{r/p}\big(u(E_x)\big){1\over n!}\log^n{\alpha s\over x}\,dx\\
\ge(\alpha s)^{-r}&&\int_0^s x^{r-1}W^{r/p}\big(u(E_x)\big){1\over n!}\log^n{s\over x}\,dx.
\end{eqnarray*}
Thus, from (\ref{u((E_s)_t)}) it follows that
$$
{1\over s}\int_0^s x^{r-1}W^{r/p}\big(u(E_x)\big){1\over n!}\log^n{s\over x}\,dx
\le  (B\alpha^r)^{n+1}s^{r-1}W^{r/p}\big(u(E_s)\big).
$$
Integrating in $s\in(0,1)\/$ both members of this inequality and by  (\ref{u(E_t)}), 
we deduce that
$$
\int_0^1{1\over s}\int_0^s x^{r-1}W^{r/p}\big(u(E_x)\big){1\over n!}\log^n{s\over x}
\,dx\,ds\le B(B\alpha^r)^{n+1}W^{r/p}\big(u(E)\big),
$$
and changing the order of integration in the left hand side, we obtain 
$$
\int_0^1 x^{r-1}W^{r/p}\big(u(E_x)\big){1\over(n+1)!}\log^{n+1}{1\over x}\,dx\le 
B(B\alpha^r)^{n+1}W^{r/p}\big(u(E)\big),
$$
as we wanted to prove.

Let now $R\in(0,1)\/$. The inequality (\ref{u(E_t)2}) can be written in the following form:
$$
\int_0^1 t^{r-1}W^{r/p}\big(u(E_t)\big){\big({R\over B\alpha^r}\log{1\over t}\big)^n\over n!}
\,dt\le BR^n W^{r/p}\big(u(E)\big),
$$
and summing in  $n\/$ we obtain, with $\delta=R/(B\alpha^r)>0\/$,
$$
\int_0^1 t^{r-\delta-1}W^{r/p}\big(u(E_t)\big)\,dt\le{B\over 1-R}W^{r/p}\big(u(E)\big),
$$
equivalently
$$
\big\|M\chi_E\big\|_{\Lambda^{\tilde{p},\tilde{s}}_u(w)}\le\widetilde{C}
\|\chi_E\|_{\Lambda^{\tilde{p}}_u(w)},\qquad E\subset\Rn,
$$
where $\tilde{p}=p(r-\delta)/r<p,\;\tilde{s}=r-\delta\/$ (see Proposition~\ref{pro: II.2.5}). 
In particular, 
$$
\big\|M\chi_E\big\|_{\Lambda^{\tilde{p},\infty}_u(w)}\le
\widetilde{C}\|\chi_E\|_{\Lambda^{\tilde{p}}_u(w)},\qquad E\subset\Rn.
$$
Also, by Theorem~\ref{adfc} and  Lemma~\ref{le:Bqp} we have that $W\in\Delta_2\/$. Then, we can
apply   Theorem~\ref{te:cnadfc} (with $u_0=u_1=u,\;w_0=w_1=w,\;p_0=p_1=\tilde{p}\/$) to conclude
the result.
$\qquad\qed$
\edem

An immediate consequence of Theorem~\ref{apfciaq} and Lemma~\ref{le:Bqp}
is the following: 

\bcor
\label{B(p-epsilon)}
Let $0<p<\infty\/$ and let $u\/$ be an arbitrary weight in $\Rn\/$. Then,
\begin{enumerate}
\item If $w\in B_p(u)\/$ there exists $q<p\/$ such that $w\in B_q(u)\/$.

\item  $B_p(u)\subset B_p\/$, that is, the boundedness $M:\lo\to\lo\/$ implies $M:\llo\to\llo\/$.
\end{enumerate}
\ecor

We can now obtain the characterization of the boundedness  
$M:\lo\to\lo\/$  or equivalently of the classes  $B_p(u)\/$ for every   $u\/$.

\bteo
\label{te:principal}
Let $u,w\/$ be weights in  $\Rn\/$ and $\R^+\/$ respectively. If $0<p<\infty\/$ 
the following results are equivalent:
\begin{enumerate}
\item[(i)] $M:\lo\to\lo\/$.

\item[(ii)]  $\big\|M\chi_E\big\|_{\Lambda^p_u(w)}\le C\|\chi_E\|_{\Lambda^p_u(w)},
\qquad E\subset\Rn\/$.

\item[(iii)]  $M:\lox{q}{u}{w}\to\lox{q}{u}{w}\/$, with $q\in(0,p)\/$.

\item[(iv)]  There exists $q\in(0,p)\/$ such that $\big\|M\chi_E\big\|_{\Lambda^{q,\infty}_u(w)}
\le C\|\chi_E\|_{\Lambda^q_u(w)},\quad E\subset\Rn\/$.

\item[(v)]  There exists $q\in(0,p)\/$ such that, for every finite family of cubes  $(Q_j)_{j=1}^J\/$
and every family of measurable sets 
$(E_j)_{j=1}^J\/$ with $E_j\subset Q_j,$ for every $ j,\/$ we have that
\bequ
\label{cprincipal}
{W\big(u\big(\bigcup_j Q_j\big)\big)\over W\big(u\big(\bigcup_j E_j\big)\big)}\le 
C\max_j\bigg({|Q_j|\over|E_j|}\bigg)^q. 
\eequ

\item[(vi)]  $\big((M\chi_E)^\ast_u(t)\big)^q\le C\displaystyle{W\big(u(E)\big)\over W(t)},\;t>0,\;E
\subset\Rn,\/$ with $q\in(0,p)\/$ independent of $t$ and $E\/$.
\end{enumerate}
\eteo

\bdem

(i) $\Rightarrow$ (ii) is immediate and  (ii)$\Rightarrow$(iii) is given in Theorem~\ref{apfciaq}.
(iii)$\Rightarrow$(iv) is also immediate and   (iv)$\Rightarrow$(i) is a consequence of 
Theorem~\ref{te:cnadfc} (since $W\in\Delta_2\/$ by Theorem~\ref{adfc} and  Lemma~\ref{le:Bqp}).  The
equivalence (iv)$\Leftrightarrow$(v) is Theorem~\ref{adfc}. On the other hand,  (vi) implies, for
$E\subset\Rn\/$,
$$
\big\|M\chi_E\big\|_{\Lambda^{q,\infty}_u(w)}=\sup_{t>0}W^{1/q}
(t)(M\chi_E)^\ast_u(t)\le CW^{1/q}\big(u(E)\big)=C\|\chi_E\|_{\Lambda^q_u(w)},
$$
which is condition  (iv). Finally, (iii) implies $M:\lox{q}{u}{w}
\to\loinx{q}{u}{w}\/$ and by Theorem~\ref{adpp}, we get  (vi).$\qquad\qed$
\edem

\bobs
\label{ob:disj}
\item As we mentioned in Remark~\ref{ob:cd}, we can assume that the sets $(E_j)_j\/$ 
in  (\ref{cprincipal}) are  disjoint. If the weight $u\/$ is doubling\index{weight!doubling} (i.e.,
$u(2Q)\le Cu(Q)$), then  also the cubes  $(Q_j)_j\/$ can be taken disjoint.

\item  Let us assume that the weight $w\/$ satisfies the following property: for every $\alpha>1\/$ 
there exists a constant $C_\alpha\/$ such that
\bequ
\label{Sum}
{W\big(\sum_j r_j\big)\over W\big(\sum_j t_j\big)}\le C_\alpha\max_j\bigg({W(r_j)
\over W(t_j)}\bigg)^\alpha,
\eequ
for every finite family of   positive numbers $\big\{(r_j,t_j)\big\}_{j=1}^m\/$, 
with $0<t_j<r_j,\;j=1,\dots,m\/$. Then, we only need to check condition  (\ref{cprincipal}) for a
unique cube 
$Q_j\/$ and a unique set  $E_j\/$, that is, it is equivalent to the inequality
\bequ
\label{cps}
\frac{W\big(u(Q)\big)}{|Q|^q }\le C\frac{W\big(u(E)\big)}{ |E|^q },\qquad E\subset Q, 
\eequ
for  $q<p$, and for every cube  $Q\subset\Rn\/$. To see this, observe that this condition is a 
consequence of (\ref{cprincipal}) and, if  (\ref{cps}) holds and  $(E_j)_j\/$ 
is a disjoint family with
$E_j\subset Q_j\/$, then
\begin{eqnarray*}
{W\big(u\big(\bigcup_j Q_j\big)\big)\over W\big(u\big(\bigcup_j E_j\big)\big)}&
\le& {W\big(\sum_j u(Q_j)\big)\over W\big(\sum_j u(E_j)\big)}\le C_\alpha\max_j\bigg({W(u(Q_j))
\over
W(u(E_j))}\bigg)^\alpha\\  &\le&  C C_\alpha\max_j\bigg({|Q_j|\over|E_j|}\bigg)^{\alpha q},\\ 
\end{eqnarray*}
and it is enough to take $\alpha>1\/$ with $\alpha q<p\/$. We observe that
every power weight  $w(t)=t^\alpha,\;\alpha>-1,\/$ satisfies the above condition.

\item  If $w=1,\;p>1,\/$ condition (\ref{cprincipal})  says that $u\in A_p\/$, since this
last condition is equivalent  (see, for example \cite{St}) to the existence of  $q\in(1,p)\/$ 
such that
$$
{u(Q)\over|Q|^q}\le C{u(E)\over|E|^q},\qquad E\subset Q,
$$ 
for every cube  $Q\/$,  which is equivalent to (\ref{cprincipal}) (by (ii)).

\item If $u=1\/$, then it is immediate to see that  (\ref{cprincipal}) is equivalent to  $w\in
B_p\/$.
\eobs

This last observation can be generalized in the following way.

\bteo
\label{sol-uA_1}
If $u\in A_1\/$ then, 
$$
M:\lo\to\lo\;\Leftrightarrow\; M:\llo\to\llo,\qquad0<p<\infty.
$$
With more generality, if $0<p<\infty,$ then $B_p(u)=B_p\/$ if and only if $u\in\bigcap_{q>1}A_q\/$. 
\eteo

\bdem
We already know, by  Corollary~\ref{B(p-epsilon)},  that $B_p(u)\subset B_p\/$. 
On the other hand, if $w\in B_p\/$, there exists $l<p\/$ such that $w\in B_l\/$ and we
have that ${W(r)/ W(t)}\le	C\big({r/ t}\big)^l,\;0<t<r<\infty\/$. Let $s>1\/$ be such
that  $sl<p\/$. If
$u\in\bigcap_{q>1}A_q\/$, in  particular $u\in A_s\/$, and it holds that  ${u(Q)/
u(E)}\le C\big({|Q|/ |E|}\big)^s,\;E\subset Q\/$. Therefore, for every family
$(Q_j)_j\/$ of cubes and $E_j\subset Q_j\/$ (pairwise disjoint),
\begin{eqnarray*}
{W\big(u\big(\bigcup_j Q_j\big)\big)\over W\big(u\big(\bigcup_j E_j\big)\big)}&
\le& {W\big(\sum_j u(Q_j)\big)\over W\big(\sum_j u(E_j)\big)}\le C\bigg({\sum_j u(Q_j)\over\sum_j
u(E_j)}\bigg)^l\\
&\le &C\max_j\bigg({u(Q_j)\over u(E_j)}\bigg)^l \le   C\max_j\bigg({|Q_j|\over
|E_j|}\bigg)^{sl},\\ 
\end{eqnarray*} 
which implies that $w\in B_p(u)\/$ (Theorem~\ref{te:principal} (v)).

Conversely, let us assume now that $B_p(u)=B_p\/$. For every $q<p\/$ the weight
$w(t)=t^{q-1}\/$ is in  $B_p=B_p(u)\/$ and by  (\ref{cprincipal}) we have that, with $q_1\in(q,p)\/$,
$$
\bigg({u(Q)\over u(E)}\bigg)^q={W\big(u(Q)\big)\over W\big(u(E)\big)}\lesssim
\bigg({|Q|\over|E|}\bigg)^{q_1},
$$
and thus, 
$$
{u(Q)\over|Q|^{q_1/q}}\le C{u(E)\over|E|^{q_1/q}},\qquad E\subset Q,
$$
for every cube $Q\/$ and hence, if $r=q_1/q\in(1,p/q)\/$, we obtain $u\in\bigcap_{s>r}A_s\/$. 
Since this argument works for  $q<p\/$, we deduce that   $u\in\bigcap_{s>1}A_s.\qquad\qed$
\edem

\bobs

(i) We know that the boundedness of  $M:\llo\to\llo\/$ holds if and only if $w\in B_p\/$ and, on
the other hand, 
$M:\Lambda^p_u\to\Lambda^p_u,\;p>1\/$ (that is, 
$M:L^p(u)\to L^p(u)\/$)
if and only if  $u\in A_p\/$. Then, one could think that the boundedness $M:\lo\to\lo,\;p>1,\/$ 
holds if and only if  $u\in A_p,\;w\in B_p\/$, that is,   $B_p(u)=B_p\/$ if $u\in A_p\/$. The
previous result shows that this is not true. In fact, as a consequence of it, 
$B_p(u)\neq
B_p,\;1<p<\infty,\/$ if $u\in A_p\setminus A_q\/$ with $1<q<p\/$.  

(ii)  The fact that $B_p\subset B_p(u)\/$ if $u\in\bigcap_{q>1}A_q\/$ was already proved in
\cite{CS3} and \cite{Ne2}, for the case $p\ge1\/$.
\eobs

Using the same idea one can analogously prove the following result which together with  Theorem~
\ref{sol-uA_1} improves Corollary~3.3 in  \cite{CS3} and  Theorem~4.1 in  \cite{Ne2}, to consider the whole
range $0<q<\infty$.

\bteo
If $1<p<\infty\/$ and $u\in A_p\/$, then $B_{q/p,
\infty}\subset B_q(u),\;0<q<\infty\/$.
\eteo

By Proposition~\ref{pro: cnu_0} we know that the boundedness  of  $M:\lo\to\lo\/$
implies  
$u(\Rn)=\infty\/$. This is, essentially, the best we can say about  $u\/$, since there are examples in
which  
$u\/$ is not in any $A_p\/$ class. In fact, the following result proves something stronger: the
weight 
$u\/$ could be not doubling. See also Proposition~\ref{W(r)/r^epsilon}.

\bteo
\label{te: pnd}
 If $u(x)=e^{|x|},\;x\in\R,\/$ and $w=\chi_{(0,1)}\/$, 
we have that $M:\lox{q}{u}{w}\to\lox{q}{u}{w}\/$ for every $q>1\/$. Therefore, it is not necessary,
in general, that the weight
$u\/$ is doubling to have the boundedness $M:\lo\to\lo\/$.
\eteo

\bdem
Since the weight $w\/$ satisfies condition (\ref{Sum}) of Remark~\ref{ob:disj},
it is enough to show that, for every cube  $Q\subset\R\/$ we have,
\bequ
\label{WuQE}
\frac{W(u(Q))}{|Q|}\le C\frac{W(u(E))}{|E|},\qquad E\subset Q. 
\eequ
If $u(E)\ge1\/$ then $u(Q)\ge1\/$, $W(u(E))=W(u(Q))=1\/$, and (\ref{WuQE}) holds
(with $C=1\/$) trivially. Thus, we can assume that $u(E)<1\/$. Then,  $W(u(E))=u(E)\/$ and 
(\ref{WuQE})
is equivalent, by Lebesgue differentiation theorem, to
\bequ
\label{W(u(Q))}
{W(u(Q))\over|Q|}\le Cu(x),\qquad\hbox{a.e.\ }x\in Q.
\eequ
To prove  (\ref{W(u(Q))}), we can assume $Q=(a,b)\/$ with $0\le a<b\/$ since   $u\/$ 
is an even function:
for cubes of the form $Q=(-a,b)\/$ we have
$$
{W(u(Q))\over|Q|}\le{W\big(2u(0,b)\big)\over|(0,b)|}\le 2{W\big(u(0,b)\big)\over|(0,b)|},
$$
and we use that the essential infimum of $u\/$ in $(0,b)\/$ and in  $(-a,b)\/$ coincide.
Let then  $Q=(a,b)\/$ with $0\le a<b\/$ and let us see  (\ref {W(u(Q))}). The infimum of  $u\/$ in 
$Q\/$
 is $e^a\/$ and, hence, we have to show that  ${W(u(Q))/|Q|}\le Ce^a\/$. If  
$e^b-e^a=u(Q)\ge1\/$
then $b\ge\log(1+e^a)\/$ and hence
$$
{W(u(Q))\over|Q|}={1\over b-a}\le{1\over\log(1+e^a)-a}={1\over\log(1+e^{-a})}\le Ce^a.
$$
If, on the contrary, $e^b-e^a=u(Q)<1\/$, we have that  $b-a<e^b-e^a<1\/$ and therefore,
$$
{W(u(Q))\over|Q|}={u(Q)\over|Q|}={e^b-e^a\over b-a}\le e^b=e^a e^{b-a}\le e e^a.\qquad\qed
$$
\edem

\section{Weak-type inequality}
\label{se:wtin}

In this section, we shall study necessary  and/or sufficient conditions 
for the weak-type inequality to hold,
$$
M\,:\,\lo\longrightarrow\loin.
$$
The following result (see \cite{CS3})  gives  a necessary condition.

\bteo
\label{te:c7p}
Let  $0<p_0,p_1<\infty\/$ and 
let us assume that $M:\loz\longrightarrow\loinu\/$. Then, there exists a constant $C>0\/$ such that
\bequ
\label{c7p}
\big\|u_0^{-1}\chi_Q\big\|_{(\loz)^\prime}\big\|\chi_Q\big\|_\lou\le C|Q| 
\eequ
for every cube  $Q\subset\Rn\/$. Here,  $\|\cdot\|_{(\loz)^\prime}\/$ denotes the norm in the
associate space (c.f. Definition~\ref{def:II.4.1},  Theorem~\ref{te: II.4.7}).
\eteo

In what follows, we shall study condition (\ref{c7p}) and we shall obtain some important consequences. For
example, combining the previous statement with  some of the results  in chapter~\ref{ch:lor},
we obtain a condition that reduces (depending on  $w_0\/$) the range of indices
$p_0\/$ for which the boundedness of  $M:\loz\to\loinu\/$ can be true. To this end, we
 define the index
$p_w\in[0,\infty)\/$ as follows:
\bequ
\label{defpw}
p_w=\inf\bigg\{p>0\,:\,{t^p\over W(t)}\in L^{p^\prime-1}\bigg((0,1),{dt\over t}\bigg)\bigg\},
\eequ
(where $p^\prime=\infty\/$ if $0<p\le1\/$).

\bteo
Let $0<p_1<\infty\/$ and let us assume that we have the boundedness of 
$M:\loz\longrightarrow\loinu\/$. Then, $p_0\ge p_{w_0}\/$. If $p_{w_0}>1\/$ the previous inequality
is strict.
\eteo

\bdem
If $M:\loz\longrightarrow\loinu\/$, we have, by  Theorem~\ref{te:c7p}, 
$(\Lambda^{p_0}_{u_0}(w_0))^\prime\neq\{0\}\/$. Then,  Theorem~\ref{te: II.4.11}
implies
 ${t^{p_0}/W(t)}\in L^{p_0^\prime-1}\big((0,1),{dt/ t}\big)\/$. That is,
$$
p_0\in\bigg\{p>0\,:\,{t^p\over W(t)}\in L^{p^\prime-1}\bigg((0,1),{dt\over t}\bigg)\bigg\}=I.
$$
The result we are looking for,   follows from the fact that  $I\/$ is an  interval in $[0,\infty)\/$,
unbounded in the right hand side, and open in the left hand side by 
$p_{w_0}\/$, if 
$p_{w_0}>1.\qquad\qed$
\edem

\bteo
If $p_0<1\/$ there are no weights $u_0,u_1\/$ such that  $M:L^{p_0}(u_0)
\to L^{p_1,\infty}(u_1),\;0<p_1<\infty\/$ is bounded.
\eteo

\bdem
It is enough to observe that  $L^{p_0}(u_0)=\Lambda^{p_0}_{u_0}(1)\/$ and that, if $w=1,\;p_w=1\/$
(c.f. (\ref{defpw})). $\qquad\qed$
\edem

To find equivalent integral expression to (\ref{c7p}), it will be useful to associate to each weight 
$u\/$ in 
$\Rn\/$ the family of functions
$\{\phi_Q\}_Q\/$ defined in the following way. For every cube  $Q\subset\Rn\/$,
\bequ
\label{defphiQ}
\phi_Q(t)=\phi_{Q,u}(t)={u(Q)\over|Q|}\int_0^t\big(u^{-1}\chi_Q\big)^\ast_u(s)\,ds\,,\quad t\ge 0.
\eequ
It will be also very useful the right derivative of the function $\phi_Q$
\bequ
\label{phiQprime}
\phi_Q^\prime(t)={u(Q)\over|Q|}\big(u^{-1}\chi_Q\big)^\ast_u(t),\quad t\ge 0.
\eequ
 Observe that $\phi_Q^\prime(t)=0\/$ if $t\ge u(Q)\/$, and that  $\phi_Q^\prime(t)
\le\phi_Q(t)/t,\;t\ge0\/$. Now, we can give an equivalent expression of  (\ref{c7p}) in
terms of these functions.   These integral expressions are quite similar to those for the classes
$B_{p_0,p_1,\infty}\/$ (Theorem~\ref{te: I.6.5}).

\bpro
\label{ce7p}
Let $\phi_Q=\phi_{Q,u_0}\/$.
\begin{enumerate}
\item[(a)] If $1<p_0<\infty\/$ each of the following expressions are equivalent  to (\ref{c7p}):
\begin{enumerate}
\item[(i)]
$$ \bigg(\int_0^{u_0(Q)}W_0^{1-p_0^\prime}(s)\,d\fiq^{p_0^\prime}(s)\bigg)^{1/p_0^\prime}
W_1^{1/p_1}\big(u_1(Q)\big)\le Cu_0(Q).
$$
\item[(ii)]
\begin{eqnarray*}
\bigg(\int_0^{u_0(Q)}\bigg({W_0\over\fiq}\bigg)^{-p_0^\prime}(t)
w_0(t)\,dt\bigg)^{1/p_0^\prime}W_1^{1/p_1}\big(u_1(Q)\big)&\le& Cu_0(Q)\\ 
 W_1^{1/p_1}\big(u_1(Q)\big)&\le& CW_0^{1/p_0}\big(u_0(Q)\big).
\end{eqnarray*}
\end{enumerate}
\item[(b)] If $0<p_0\le1\/$, (\ref{c7p}) is equivalent to
$$
{W_1^{1/p_1}\big(u_1(Q)\big)\over u_0(Q)}\le C{W_0^{1/p_0}(t)\over\fiq(t)},
\qquad 0<t\le u_0(Q).
$$
\end{enumerate}
We are assuming both in  (a) and   (b) that the inequalities are true for every cube 
$Q\subset\Rn\/$, and that the constant $C\/$ is independent of  $Q\/$.
\epro

\bdem
Observe that, by definition, $\displaystyle\int_0^t\big(u_0^{-1}
\chi_Q\big)^\ast_{u_0}(s)\,ds=\displaystyle{|Q|\over u_0(Q)}\fiq(t)\/$. Hence,  by Theorem~\ref{te: II.4.7}
(with
$r=u_0(Q)\/$),
\begin{eqnarray*}
\big\|u_0^{-1}\chi_Q\big\|_{(\loz)^\prime}&\approx&\displaystyle{|Q|\over u_0(Q)}
\bigg(\int_0^{u_0(Q)}W_0^{1-p_0^\prime}(s)\,d\fiq^{p_0^\prime}(s)\bigg)^{1/p_0^\prime}\\ 
&\approx&\displaystyle{|Q|\over u_0(Q)}\bigg(\int_0^{u_0(Q)}
\bigg({W_0\over\fiq}\bigg)^{-p_0^\prime}(s)w_0(s)\,ds\bigg)^{1/p_0^\prime}\\
&\qquad+&
{\displaystyle\int_0^\infty(u_0^{-1}\chi_Q)^\ast_{u_0}(s)\,ds\over
W_0^{1/p_0}\big(u_0(Q)\big)},\\  
\end{eqnarray*} 
in the case $p_0>1\/$, and 
$$
\big\|u_0^{-1}\chi_Q\big\|_{(\loz)^\prime}={|Q|\over u_0(Q)}\sup_{t>0}
{\phi_Q(t)\over W_0^{1/p_0}(t)},
$$
in the case $0<p_0\le1\/$. Since  $\big\|\chi_Q\big\|_\lou=W_1^{1/p_1}
\big(u_1(Q)\big)\/$ and \break   
$\displaystyle\int_0^\infty\big(u_0^{-1}\chi_Q\big)^\ast_{u_0}(s)\,ds=\displaystyle\int_{\Rn}u_0^{-1}(x)\chi_Q
(x)u_0(x)\,dx=|Q|\/$, the given expression can be immediately deduced  from (\ref{c7p}).
$\qquad\qed$
\edem

\bobs
In many cases, condition (\ref{c7p}) is also sufficient.
For example:

(i) If $w_0=w_1=1,\;p_0=p_1=p\ge1\/$ inequality (\ref{c7p}) is
$$
\big\|u_0^{-1}\chi_Q\big\|_{L^{p^\prime}(u_0)}\big\|\chi_Q\big\|_{L^p(u_1)}\le C|Q|,
$$
or equivalently,
\bequ\label{Ap}
\bigg({1\over|Q|}\int_Q u_1(x)\,dx\bigg)\bigg({1\over|Q|}\int_Q u_0^{-1/(p-1)}(x)\,dx\bigg)^{p-1}\le C,
\eequ
in the case $p>1\/$, and
$$
{u_1(Q)\over|Q|}\le C{u_0(E)\over|E|},\qquad E\subset Q,
$$
in the case $p=1\/$.
This is the so-called  $A_p\/$\index{$A_p$} condition for the pair  $(u_1,u_0)\/$ and it is known to be
sufficient for the boundedness
$M:L^p(u_0)\to L^{p,\infty}(u_1)\/$ or, in other terms,
$M:\lox{p}{u_0}{1}\to\loinx{p}{u_1}{1}\/$ (see \cite{GR}).

(ii) If $u_0=u_1=1,$ then $ \fiq(t)=t,\,0\le t\le1,\/$ and condition (\ref{c7p}) 
(or  any of the expressions of  Proposition~\ref{ce7p}) is equivalent to  $(w_0,w_1)\in
B_{p_0,p_1,\infty}\/$ (Theorem~\ref{te: I.6.5}). That is, the condition is also sufficient.

(iii) If $u_0=u_1=u,\,p_0=p_1=q\ge1\/$ and $w_0(t)=w_1(t)=t^{q/p-1},\,1<p<\infty,\/$ 
we are in the case $M:L^{p,q}(u)\longrightarrow L^{p,\infty}(u)\/$. Then condition  (\ref{c7p})
is equivalent to the inequality
\bequ
\label{Lpq7p}
\big\|u^{-1}\chi_q\big\|_{L^{p^\prime,q^\prime}(u)}\big\|\chi_Q\big\|_{L^{p,q}(u)}
\le C|Q|. 
\eequ
This case was studied by  Chung, Hunt, and Kurtz in  \cite{CHK} (see also \cite{HK}) 
where they proved that (\ref{Lpq7p}) is a necessary and sufficient condition. 
\eobs

The study of the properties of the functions $\{\phi_Q\}\/$ 
will allow us to better understand condition  (\ref{c7p}), and obtain some consequences. In the
following proposition, we summarize some of these properties.

\bpro
\label{6pphiQ}
Let $u\/$ be a weight in  $\Rn\/$ and, 
for every cube  $Q\subset\Rn\/$, let $\phi_Q=\phi_{Q,u}\/$. Then, 
\begin{enumerate}
\item[(i)] $\fiq(t)={u(Q)\over|Q|}\max\big\{|E|\,:\,E\subset Q,\,u(E)=t\big\},\quad 
0\le t\le u(Q). $
\item[(ii)] ${\fiq}_{|\,[0,u(Q)]}:\big[0,u(Q)\big]\longrightarrow\big[0,u(Q)\big]$ is strictly increasing,
onto, concave and absolutely continuous.

\item[(iii)] $\fiq(t)\ge t\in[0,u(Q)]\hbox{\ and\ }u\in A_1\hbox{\ if and only if\ }\fiq(t)
\le Ct,\quad 0\le t\le u(Q).  $
\item[(iv)] ${If\ }1<p<\infty,\,0<q\le p^\prime,\quad u\in A_p
\Leftrightarrow\big\|\fiq^\prime\big(u(Q)\,\cdot\big)\big\|_{L^{p^\prime,q}}\le C.  $
\item[(v)] $\hbox{If\ }u\in A_p,\;1\le p<\infty,\hbox{\ then\ } \fiq^\prime
\big(u(Q)t\big)<Ct^{-1/p^\prime},\;0<t\le1. $
\item[(vi)] $\hbox{If\ }1<p<\infty\hbox{\ and\ }\fiq^\prime\big(u(Q)t\big)
<Ct^{-1/p^\prime},\;0<t\le1,\hbox{\ then\ }u\in\bigcap_{q>p}A_q. $
\end{enumerate}
In (iii), (iv), (v), and (vi), $C\/$ represents a constant not depending on  $Q\/$.
\epro

\bdem

(i) Let $0\le t\le u(Q)\/$ and $E\subset Q\/$ be a measurable set such that $u(E)=t\/$.
 Then 
\begin{eqnarray*}
|E|&=&\int_{\Rn}u^{-1}(x)\chi_Q(x)\chi_E(x)\,
u(x)\,dx\le\int_0^\infty\big(u^{-1}\chi_Q\big)^\ast_u(s)\big(\chi_E\big)^\ast_u(s)\,ds\\
&=&\int_0^t\big(u^{-1}
\chi_Q\big)^\ast_u(s)\,ds=\displaystyle{|Q|\over
u(Q)}\fiq(t).\\
\end{eqnarray*}
On the other hand, since  $\big(Q,u(x)\chi_Q(x)dx\big)\/$ 
is a strongly resonant measure space (see \cite{BS}), there exists $f\ge0\/$  measurable in 
$Q\/$ with
$f^\ast_u=\big(\chi_E\big)^\ast_u=\chi_{[0,t)}\/$  and such that
\begin{eqnarray}
\label{eiphiQ}
{|Q|\over u(Q)}\fiq(t)&=&\int_0^\infty\big(u^{-1}\chi_Q\big)^\ast_u(s)
\big(\chi_E\big)^\ast_u(s)\,ds\nonumber\\
&=&\int_{Q}f(x)u^{-1}(x)\chi_Q(x)u(x)\,dx. 
\end{eqnarray}
It follows that  $f=\chi_{R}\/$, with $R\subset Q,\;u(R)=t\/$ and, by
(\ref{eiphiQ}), ${|Q|\over u(Q)}\fiq(t)=|R|\/$.

(ii) Let us see that ${\fiq}_{|\,[0,u(Q)]}\/$ is strictly increasing (the rest of the properties
are obvious). It is enough to prove that  $\big(u^{-1}\chi_Q\big)^\ast_u(t)>0\/$ for $0<t<u(Q)\/$.
On the contrary, there would exist  $t_0\in\big(0,u(Q)\big)\/$ such that
$\big(u^{-1}\chi_Q\big)^\ast_u(t_0)=0\/$ and hence $u\big(\big\{u^{-1}\chi_Q>0\big\}\big)\le
t_0\/$. Since $t_0<u(Q)\/$, we would have that $u\big(\big\{x\in
Q\,:\,u(x)=+\infty\big\}\big)=u\big(\big\{u^{-1}\chi_Q=0\big\}\big)>0\/$, 
which contradicts the fact that  $u\/$ is locally integrable. 

(iii) Since  $\fiq\/$ is concave in  $[0,u(Q)],\;t^{-1}\fiq(t)\/$ is decreasing and, thus,
$t^{-1}\fiq(t)\ge u(Q)^{-1}\fiq\big(u(Q)\big)=1\/$ for $0\le t\le u(Q)\/$. The second part of 
(iii) follows from  (i) and the fact that $u\in A_1\/$ if and only if
${u(Q)\over|Q|}\le C{u(E)\over|E|}\/$. 

(iv) By definition (see \cite{St}), $u\in A_p\/$
\begin{eqnarray*}
&\Leftrightarrow&{u(Q)\over|Q|}\bigg({1\over|Q|}\int_Q u^{-p^\prime/p}(x)\,dx\bigg)^{p/p^\prime}
\le
C\\
&\Leftrightarrow&{u(Q)^{1/p}\over|Q|}\bigg(\int_{\Rn}\big(u^{-1}(x)
\chi_Q(x)\big)^{p^\prime}u(x)\,dx\bigg)^{1/{p^\prime}}\le
C_1\\ 
&\Leftrightarrow&{u(Q)^{1/p}\over|Q|}\bigg(\int_0^\infty{\big(u^{-1}\chi_Q\big)^\ast_u}^{p^\prime}(s)\,ds
\bigg)^{1/{p^\prime}}\le
C_1
\Leftrightarrow\big\|\fiq^\prime\big(u(Q)\,\cdot\big)\big\|_{p^\prime}\le C_1. 
\end{eqnarray*}
Since
$L^{p^\prime,q}\subset L^{p^\prime}\/$ if $0<q\le p^\prime\/$, we have proved the sufficiency in 
(iv). Conversely, if $u\in A_p\/$ it is known that  $u\in A_{p-\epsilon}\/$ for some
$\epsilon>0\/$.  Therefore,
by the previous equivalence, $\big\|\fiq^\prime\big(u(Q)\,\cdot\big)\big\|_{(p-\epsilon)^\prime}
\le
C\/$. Observe that
$\fiq^\prime\big(u(Q)t\big)={u(Q)\over|Q|}\big(u^{-1}\chi_Q\big)^\ast_u\big(u(Q)t\big)=0\/$ for
$t\ge1\/$. Since $L^{p_1,r}(X)\subset L^{p_0,q}(X)\/$ if $p_0<p_1\/$ and $X\/$ is of finite measure 
(see \cite{BS}), we have that
$$
\big\|\fiq^\prime\big(u(Q)\,\cdot\big)\big\|_{p^\prime,q}\le\big\|\fiq^\prime\big(u(Q)\,
.\big)\big\|_{(p-\epsilon)^\prime}\le
C.
$$

(v) The case $p=1\/$ is a consequence of  (iii) and the inequality $\fiq^\prime(t)
\le\fiq(t)/t\/$. If $u\in A_p,\;1<p<\infty\/$ then, by (iv), we have that 
$\big\|\fiq^\prime\big(u(Q)\,\cdot\big)\big\|_{p^\prime,\infty}\le\big\|\fiq^\prime\big(u(Q)\,\cdot\big)\big\|_{p^\prime}\le
C\/$. Since $\fiq^\prime\big(u(Q)\,\cdot\big)\/$ is decreasing, right continuous and equal to zero in 
$[1,\infty)\/$, 
$\big\|\fiq^\prime\big(u(Q)\,\cdot\big)\big\|_{p^\prime,\infty}=\sup_{0<t<1}t^{1/p^\prime}
\fiq^\prime\big(u(Q)t\big)\/$
and we obtain (v).

(vi) The hypothesis implies  $\big\|\fiq^\prime\big(u(Q)\,\cdot\big)\big\|_{L^{p^\prime,\infty}}
\le C\/$. But $\fiq^\prime\big(u(Q)\,\cdot\big)\/$ is supported in $[0,1]\/$ and therefore,
$$
\big\|\fiq^\prime\big(u(Q)\,\cdot\big)\big\|_{L^{q^\prime}}\le\big\|\fiq^\prime\big(u(Q)
\,\cdot\big)\big\|_{L^{p^\prime,\infty}}\le C,
$$
for $q>p\/$ and (vi) follows from (iv).$\qquad\qed\/$ 
\edem

Let us see a useful consequence of the two above propositions.

\bpro
\label{pro: cnadnd}
 If $0<p_0,p_1<\infty\/$ and $M:\loz\to\loinu\/$, then  
\begin{eqnarray*}
(i)\quad&\displaystyle{W_1^{1/p_1}\big(u_1(Q)\big)\over u_0(Q)}\le
C\displaystyle{W_0^{1/p_0}(t)\over\fiq(t)},
\qquad 0<t\le u_0(Q),\\ 
(ii)\quad&\displaystyle{W_1^{1/p_1}\big(u_1(Q)\big)\over u_0(Q)}\le
C\displaystyle{W_0^{1/p_0}\big(u_0(S)\big)
\over u_0(S)},\quad S\subset Q.
\end{eqnarray*}
Here $\phi_Q=\phi_{Q,u_0}\/$ and we are assuming that the previous inequalities are satisfied for
every cube  $Q\subset\Rn\/$, with
$C\/$ independent of  $Q\/$.
\epro

\bdem
 (ii) is a consequence of  (i) and Proposition~\ref{6pphiQ} (iii). 
Hence, we only have to show  (i). If $p_0\le1\/$, (i) is Proposition~\ref{ce7p} (b).
For $p_0>1\/$ we have, by Proposition~\ref{ce7p} (a.i),
$$
\bigg(\int_0^t W_0^{1-p_0^\prime}(s)\,d\fiq^{p_0^\prime}(s)\bigg)^{1/p_0^\prime}W_1^{1/p_1}
\big(u_1(Q)\big)\le Cu_0(Q),
$$
for $0<t\le u_0(Q)\/$. Since $W_0\/$ is increasing, it follows that
$$
W_0^{-1/p_0}(t)\bigg(\int_0^t d\fiq^{p_0^\prime}(s)\bigg)^{1/p_0^\prime}W_1^{1/p_1}
\big(u_1(Q)\big)\le Cu_0(Q),
$$
and we obtain (i).$\qquad\qed$
\edem

The following result establishes that in the boundedness
$M:\lox{p_0}{u}{w_0}\to\loinx{p_1}{u}{w_1}\/$ we can always assume that  $u\equiv1$. 

\bteo
\label{te: MubMb}
 If $0<p_0,p_1<\infty\/$ and $M:\lox{p_0}{u}{w_0}\to\loinx{p_1}{u}{w_1}\/$, 
we have that
$$
M\,:\,\lloz\longrightarrow\lloinu,
$$
that is, $(w_0,w_1)\in B_{p_0,p_1,\infty}\/$.
\eteo

\bdem
If $p_0>1\/$ then, by Proposition~\ref{ce7p} (a.ii) and using that  $\fiq(t)\ge t,\;0<t<u(Q)\/$
(Proposition~\ref{6pphiQ}  (iii)) we have the inequalities
\begin{eqnarray*}
\bigg(\int_0^{u(Q)}\bigg({W_0(t)\over t}\bigg)^{-p_0^\prime}w_0(t)\,dt\bigg)^{1/p_0^\prime}
W_1^{1/p_1}\big(u(Q)\big)&\le& Cu(Q),\\ 
W_1^{1/p_1}\big(u(Q)\big)&\le& CW_0^{1/p_0}\big(u(Q)\big),
\end{eqnarray*}
for every cube $Q\subset\Rn\/$. But we have proved  (Proposition~\ref{pro: cnu_0}) that
$u(\Rn)=\infty\/$.
Then, for every  $r>0\/$ there exists a cube  $Q\/$ with $u(Q)=r\/$. Therefore, the two previous
inequalities are equivalent to the condition (a.iii) of Theorem~\ref{te: I.6.5}  and hence, 
$(w_0,w_1)\in\BBinnux{p_0,p_1}\/$.

If $p_0\le1\/$  we apply Proposition~\ref{ce7p}(b) 
to obtain the expression of  Theorem~\ref{te: I.6.5} to conclude the same.$\qquad\qed\/$.
\edem

\bcor
\label{Bpinfi(u)}
 $B_{p,\infty}(u)\subset B_{p,\infty},\;0<p<\infty\/$.
\ecor

\bobs
\label{obs7p}
(i) If $M:\Lambda^1_u(w)\to\Lambda^{1,\infty}_u(w)\/$ then, 
by  Corollary~\ref{Bpinfi(u)}, $w\in B_{1,\infty}\/$. The bigger class $B_p\/$ contained in 
$B_{1,\infty}\/$ is $B_1\/$ and if  $w\in B_1\/$, Theorem~\ref{cnsadnd} characterizes the previous
boundedness. That is, the problem $M:\Lambda^1_u(w)\to\Lambda^{1,\infty}_u(w)\/$ only remains open
in the case
$w\in B_{1,\infty}\setminus B_1\/$.

(ii) Condition (\ref{c7p}) is sufficient in the case
$M:\Lambda^1_u(w)\to\Lambda^{1,\infty}_u(w)\/$ if the weight  $w\/$ satisfies 
$$
\sum_{j=1}^n{W(t_j)\over W(r_j)}r_j\le C{W\big(\sum_j t_j\big)
\over W\big(\sum_j r_j\big)}\sum_j r_j,
$$
for every finite families of numbers  $\big\{(t_j,r_j)\big\}_j\/$ with $0<t_j<r_j,\;j=1,
\dots,n\/$. To see this, we observe that  if the previous condition and  (\ref{c7p}) hold (which
implies  (\ref{dEQ})) we have, for every finite family of disjoint cubes  $(Q_j)_j\/$ and every
family of sets 
$(E_j)_j\/$ with
$E_j\subset Q_j\/$,
$$
\sum_j{|E_j|\over|Q_j|}u(Q_j)\le C\sum_j{W(u(E_j))\over W(u(Q_j))}u(Q_j)
\le C{W\big(u\big(\bigcup E_j\big)\big)\over W\big(u\big(\bigcup Q_j\big)\big)}\sum_j u(Q_j)
.
$$
This also holds (since $W\in\Delta_2\/$) if the cubes  $Q_j\/$ 
are \lq\lq almost\rq\rq\ disjoint (if $\sum_j\chi_{Q_j}\le k=k_n\/$). If $0\le f\in\M(\Rn)\/$  is
bounded and with compact support and   $t>0\/$, the level set  $E_t=\{Mf>t\}\/$ is contained
(see \cite{LN1}) in a finite union of cubes $(Q_j)_j\/$ such that $\sum_j\chi_{Q_j}\le k=k_n\/$
(only depending on the dimension) and  ${1\over|Q_j|}\int_{Q_j}f(x)\,dx>t,$ for every $ j\/$. Then, 
\begin{eqnarray*}
t\sum_j u(Q_j)&\le&\int_{\R^n} f(x)\bigg(\sum_j{u(Q_j)\over|Q_j|}u^{-1}(x)\chi_{Q_j}(x)\bigg)u(x)\,dx\\
&\le&\|f\|_{\Lambda^1_u(w)}\bigg\|\sum_j{u(Q_j)\over|Q_j|}u^{-1}\chi_{Q_j}
\bigg\|_{(\Lambda^1_u(w))^\prime}\\
\end{eqnarray*}
and, to prove that  $\|Mf\|_{\Lambda^{1,\infty}_u(w)}\le C\|f\|_{\Lambda^1_u(w)}\/$ is bounded, 
we only have to see that
$$
\bigg\|\sum_j{u(Q_j)\over|Q_j|}u^{-1}\chi_{Q_j}\bigg\|_{(\Lambda^1_u(w))^\prime}
\lesssim{\sum_j u(Q_j)\over W\big(u\big(\bigcup_j Q_j\big)\big)},
$$
(since then we have that $tW(u(E_t))\le tW\big(u\big(\bigcup_j Q_j\big)\big)
\le C\|f\|_{\Lambda^1_u(w)}\/$).
But, since $\sum_j\chi_{Q_j}\le k\/$,
\begin{eqnarray*}
\int_0^t\bigg(\sum_j{u(Q_j)\over|Q_j|}u^{-1}\chi_{Q_j}\bigg)_u^\ast(s)\,ds&=\sup_{u(E)=t}
\displaystyle\int_E\sum_j{u(Q_j)\over|Q_j|}u^{-1}(x)\chi_{Q_j}(x)\,u(x)\,dx\\ 
&=\sup_{u(E)=t}\displaystyle\sum_j{|E\cap Q_j|\over|Q_j|}u(Q_j)\\ 
&\lesssim\displaystyle{W(t)\over W\big(u\big(\bigcup_j Q_j\big)\big)}\sum_j u(Q_j),
\end{eqnarray*}
and  (see Theorem~\ref{te: II.4.7}(i))    the inequality we are looking for follows.

(iii) If $w\in L^1(\R^+)$, the condition on  $w\/$ of the previous observation can be weakened up:
it is sufficient that it holds \lq\lq near  0\rq\rq, that is, with $\sum r_j<\epsilon\/$ (for some
fixed $\epsilon>0\/$).

(iv) It is easy to see (using the two previous observations) that, for the weights
$w=\chi_{(0,1)},\;w(t)=\big(\log^+(1/t)\big)^\alpha,\;w(t)=\big(\log^+\log^+(1/t)\big)^\alpha,
\;w(t)=t^\alpha\/$,
condition (\ref{c7p}) is also equivalent (in the case
$p_0=p_1=1,\;u_0=u_1,\;w_0=w_1=w\/$.)
\eobs

\bteo
\label{Bq/p}
 Let $1\le p<\infty,\,0<q<\infty\/$. Then, 

\medskip
(i)  $B_{q/p}\subset B_{q,\infty}(u)\/$ implies $u\in\bigcap_{r>p}A_r\/$.

(ii) If $q\le1,\;B_{q,\infty}\subset B_{q,\infty}(u)\/$ implies $u\in A_1\/$.
\eteo

\bdem
(i) It is immediate to check (Theorem~\ref{te: I.6.6}) that the weight $w(t)=t^{r-1}\/$ is in 
$B_{q/p}\subset B_{q,\infty}(u)\/$, for  $0<r<q/p\/$. Hence, by  Proposition~\ref{pro:
cnadnd} (i),   we have that,
for every cube  $Q\subset\Rn\/$ and $0<t<u(Q)\/$,
$$
\fiq(t)\le Cu(Q){W^{1/q}(t)\over W^{1/q}\big(u(Q)\big)}=C u(Q)^{1-r/q}\,t^{r/q}.
$$
Since $\fiq^\prime(t)\le \fiq(t)/t,\;0<t<u(Q),\/$ we obtain that  $\fiq^\prime\big(u(Q)t\big)
\le Ct^{r/q-1},\;0<t<1,\/$  and by  Proposition~\ref{6pphiQ}(vi), $u\in A_s\/$ for every
$s>q/r\/$. Since this holds for  $0<r<q/p\/$, we conclude that $u\in\bigcap_{r>p}A_r\/$.

(ii) Observe that  (using, for example,  Theorem~\ref{te: I.6.5})  $w(t)=t^{q-1}\in B_{q,\infty}
\subset B_{q,\infty}(u)\/$. From  Proposition~\ref{pro: cnadnd}(i) (with
$w_0=w_1=t^{q-1}\;,u_0=u_1=u\;,p_0=p_1=q\/$)  we obtain now that  $\fiq(t)\le Ct,\;0<t<u(Q),\/$ and
the result  is a consequence of Proposition~\ref{6pphiQ} (iii). $\qquad\qed$
\edem

The following result  is related to the  $A_\infty\/$ class. 
Let us recall that (see, for example, \cite{St})   $u\in A_\infty=\bigcup_{p\ge1}A_p\/$ if and only
if there exist constants 
$\epsilon,C>0\/$ so that
\bequ\label{Ain}
{u(Q)\over|Q|^\epsilon}\ge C{u(E)\over|E|^\epsilon}
\eequ
for every cube  $Q\/$ and every measurable set  $E\subset Q\/$. In general, 
the boundedness of $M:\lo\rightarrow\loin\/$
does not imply $u\in A_\infty\/$ (Theorem~\ref{te: pnd}). But we can still  give a sufficient
condition on the weight $w\/$.

\bpro
\label{W(r)/r^epsilon}
Let us assume that there exist constants $\epsilon,C>0\/$ so that
$$
{W(r)\over r^\epsilon}\ge C{W(t)\over t^\epsilon},\quad0<t<r<\infty.
$$
Then $A_p(w)\subset A_\infty,\;0<p<\infty\/$. 
\epro

\bdem 
We can assume  $\epsilon<p\/$. If $u\in A_p(w)\/$, by Proposition~\ref{pro: cnadnd},
$$
\fiq(t)\le C_1{W^{1/p}(t)\over W^{1/p}\big(u(Q)\big)}u(Q)\le C_2
\bigg({t\over u(Q)}\bigg)^{\epsilon/p}u(Q),
$$
for every $Q\/$ and $0<t<u(Q)\/$. Then,
$$
\fiq^\prime\big(u(Q)t\big)\le{1\over u(Q)t}\fiq\big(u(Q)t\big)\le C_2 t^{\epsilon/p-1},\quad0<t<1,
$$
and, from Proposition~\ref{6pphiQ} (vi), it follows that  $u\in A_{q^\prime}\/$, for the range 
$q<p/(p-\epsilon).\qquad\qed$ 
\edem

\bobs The condition on $w\/$ of  Proposition~\ref{W(r)/r^epsilon} is not equivalent to 
  $w\in\bigcup_p B_p\/$ (as it happens with the classes
$A_p\/$). For example, the weight $w=\chi_{(0,1)}\in B_{1,\infty}\/$ does not satisfy that
condition. Even if we impose  $w\notin L^1(\R^+)$ as the example 
$w(t)={1/( 1+t)}\/$ shows.

Power weights $w(t)=t^\alpha\/$  trivially satisfy the condition of 
Proposition~\ref{W(r)/r^epsilon}. Another nontrivial example, of a function satisfying such
condition (with
$\epsilon=1/2\/$) is $w(t)=1+\log^+{(1/ t)}\/$.
\eobs

Up to now, we have only seen necessary conditions for the weak-type boundedness. 
Let us see now some sufficient conditions. For example, from  Theorem~\ref{adpp}, we deduce the two
following results.

\bteo
 Let $0<p,q<\infty\/$ and let us assume that $M:\lo\to\loin\/$ is bounded. 
Let $\tilde{w}\/$ another weight. Then we have that 
$M:\lox{q}{u}{\tilde{w}}\to\Lambda^{q,\infty}_u(\tilde{w})\/$ is bounded in each of the following
cases:  
\begin{enumerate}
\item[(i)] If $0<q\le p\/$ and the condition 
$$
{\widetilde{W}^{1/q}(r)\over W^{1/p}(r)}\le C{\widetilde{W}^{1/q}(t)\over W^{1/p}(t)},
\qquad 0<t<r<\infty
$$
holds.

\item[(ii)] If $p<q<\infty\/$ and it satisfies
$$
\bigg({1\over\widetilde{W}(t)}\int_0^t\bigg({W(s)\over\widetilde{W}(s)}\bigg)^{q\over q-p}
\tilde{w}(s)\,ds\bigg)^{q-p\over q}\le C{W(t)\over\widetilde{W}(t)},\qquad t>0.
$$
\end{enumerate}
\eteo

\bdem
Let us fix  $t>0\/$ and let us consider the weights
$$
w_1(s)={w(s)\over W(t)}\chi_{(0,t)}(s),\quad w_0(s)={\tilde{w}(s)\over
\widetilde{W}(t)}\chi_{(0,t)}(s),\quad s>0.
$$
Under the hypothesis of  (i) we have, by Theorem~\ref{twoop} (d)
$$
\sup_{g\d}{\bigg(\displaystyle\int_0^\infty g^p(t)
w_1(t)\,dt\bigg)^{1/p}\over\displaystyle\bigg(\int_0^\infty g^q(t) w_0(t)\,dt\bigg)^{1/q}}
=\sup_{s>0}{W_1^{1/p}(s)\over W_0^{1/q}(s)}\le C.
$$
Therefore,
$$
\bigg({1\over W(t)}\int_0^t (f^\ast_u)^p(s) w(s)\,ds\bigg)^{1/p}\le C\bigg({1\over \widetilde{W}(t)}
\int_0^t (f^\ast_u)^q (s)\tilde{w}(s)\,ds\bigg)^{1/q}
$$ 
and we obtain  (i).

To see  (ii), observe that the hypothesis implies, with $r=q/p>1\/$,
$$
B=B(t)=\bigg(\int_0^\infty\bigg({W_1(s)\over W_0(s)}\bigg)^{r^\prime}w_0(s)\,ds\bigg)^{1/r^\prime}+
{W_1(\infty)\over W^{1/r}_0(\infty)}\le C+1<\infty,
$$
and by Theorem~{I.5.7} in \cite{Saw2}, we have
\begin{eqnarray*}
\sup_{g\d}{\bigg(\displaystyle\int_0^\infty g^p
(s)w_1(s)\,ds\bigg)^{1/p}\over\bigg(\displaystyle\int_0^\infty g^q(s) w_0(s)\,ds\bigg)^{1/q}}
&=&\Bigg(\sup_{g\d}{\displaystyle\int_0^\infty g(s) w_1(s)\,ds\over\bigg(\displaystyle\int_0^\infty g^r(s)
w_0(s)\,ds\bigg)^{1/r}}\Bigg)^{1/p}\\  &\le& (C^\prime B)^{1/p}\le
\big(C^\prime(C+1)\big)^{1/p}=C^{\prime\prime}<\infty.
\end{eqnarray*}
In particular,
$$
\bigg({1\over W(t)}\int_0^t\big(f^\ast_u\big)^p(s) w(s)\,ds\bigg)^{1/p}\le C^{\prime\prime}\bigg({1\over 
\widetilde{W}(t)}\int_0^t\big(f^\ast_u\big)^q(s) \tilde{w}(s)\,ds\bigg)^{1/q},
$$
for every $t>0\/$ and $f\in\M(\Rn)\/$ and by Theorem~\ref{adpp} we get the result.$\qquad\qed$
\edem

\bcor
Let $0<p<\infty\/$ and let us assume that  $M:\lo\to\loin\/$ is bounded. If $0<q\le p\/$ and 
$\tilde{w}\/$ is a weight such that  $\widetilde{W}^{1/q}/W^{1/p}\/$  is decreasing, then 
$M\,:\,\lox{q}{u}{\tilde{w}}\longrightarrow\loinx{q}{u}{\tilde{w}}\/$ is also bounded.
\ecor
 
The following result completes Theorem~\ref{te: MubMb}.

\bteo
\label{aua_1}
If $u\in A_1\/$, we have the boundedness $M:\Lambda^{p_0}_u(w_0)
\to\Lambda^{p_1,\infty}_u(w_1)\/$, if and only if 
$M:\Lambda^{p_0}(w_0)\to\Lambda^{p_1,\infty}(w_1)\/$. In particular
$B_{p,\infty}(u)=B_{p,\infty},\;0<p<\infty,\/$ if $u\in A_1\/$.
\eteo

\bdem
The \lq\lq only if\rq\rq\  condition is Theorem~\ref{te: MubMb}.
To see the other implication, let us observe that if   $u\in A_1\/$  we have, by
Theorem~\ref{adpp} (applied to 
$w\equiv1,\;p=1\/$) that $(Mf)_u^\ast(t)\lesssim Af_u^\ast(t),\;t>0,\;f\in\M(\Rn),\/$ where 
 $A\/$ the Hardy operator. Since the boundedness of 
$M:\Lambda^{p_0}(w_0)\to\Lambda^{p_1,\infty}(w_1)\/$ is equivalent to  $A:L_{\hbox{\decp dec}}^{p_0}(w_0)\to
L^{p_1,\infty}(w_1)\/$, we get
$$
W^{1/p_1}_1(t)(Mf)_u^\ast(t)\lesssim\|Af_u^\ast\|_{L^{p_1,\infty}(w_1)}\lesssim 
\|f_u^\ast\|_{L^{p_0}(w_0)}=\|f\|_{\Lambda_u^{p_0}(w_0)},\qquad t>0,
$$
which is equivalent to  $\|Mf\|_{\Lambda^{p_1,\infty}_u(w_1)}\lesssim\|f\|_{\Lambda^{p_0}_u(w_0)}.
\qquad\qed$
\edem

Using the same idea, we can prove the following result analogous to Theorem~3.2 in \cite{CS3}.

\bteo
\label{ B_{p_0,p_1}}
 Let $0<p_0,p_1<\infty,\;1\le p<\infty\/$ and let us assume $u\in A_p\/$.
 Then,
\begin{enumerate}
\item[(i)] If  $(w_0,w_1)\in B_{p_0,p_1}\/$, we have that  $M:\Lambda^{pp_0}_u(w_0)
\to\Lambda^{pp_1}_u(w_1)\/$ is bounded.

\item[(ii)] If $(w_0,w_1)\in B_{p_0,p_1,\infty}\/$,  then $M:\Lambda^{pp_0}_u(w_0)\to
\Lambda^{pp_1,\infty}_u(w_1)\/$ is bounded.
\end{enumerate}
\eteo

\bdem
 If $u\in A_p\/$, by Theorem~\ref{adpp}, $(Mf)_u^\ast(t)^p\lesssim 
A\big((f_u^\ast)^p\big)(t)\/$. The hypothesis in (i) implies that 
$A:L_{\hbox{\decp dec}}^{p_0}(w_0)\to
L^{p_1}(w_1)\/$ and it follows that
\begin{eqnarray*}
\|Mf\|^p_{\Lambda^{pp_1}_u(w_1)}&=&\bigg(\int_0^\infty\big((Mf)_u^\ast\big)^{pp_1}(t)w_1(t)\,dt
\bigg)^{1/p_1}\lesssim\big\|A\big((f_u^\ast)^p\big)\big\|_{L^{p_1}(w_1)}\\ 
&\lesssim&\|(f_u^\ast)^p\|_{L^{p_0}(w_0)}=\|f\|^p_{\Lambda^{pp_0}_u(w_0)}.
\end{eqnarray*}
Analogously, one can easily prove  (ii). $\qquad\qed$
\edem

If $p>1,\;p_0\le p_1\/$,  (ii) can be improved in the following way:

\bcor
Let $0<p_0\le p_1<\infty\/$. If $u\in A_p,\;p>1,\/$ 
and  $(w_0,w_1)\in B_{p_0,p_1,\infty}\/$, then
$M:\Lambda^{pp_0}_u(w_0)\to\Lambda^{pp_1}_u(w_1)\/$ is bounded.
\ecor

\bdem
If  $u\in A_p\/$ then $u\in A_q\/$ for some $q\in(1,p)\/$, 
and by the previous theorem, $M:\Lambda^{qp_0}_u(w_0)\to\Lambda^{qp_1,\infty}_u(w_1)\/$ is bounded.
Applying now the interpolation theorem (Theorem~\ref{te: II.6.3}) we obtain the boundedness of 
$M:\Lambda^{pp_0}_u(w_0)\to\Lambda^{pp_1,pp_0}_u(w_1)\subset\Lambda^{pp_1}_u(w_1).\qquad\qed$
\edem

\bobs
In \cite{Ne2},  Neugebauer introduces the classes $A_p^\ast\/$\index{$A_p^\ast\/$} 
which are defined by 
$$
A_p^\ast=\big\{\,u\in A_\infty\,:\,p=\inf\{q\ge1\,:\,u\in A_q\}\,\big\},\qquad1\le p<\infty.
$$
These classes are pairwise disjoint and the union of all of them is $A_\infty\/$. 
With this terminology, Theorem~\ref{sol-uA_1} states that $B_p(u)=B_p,\;0<p<\infty,\/$ if
$u\in A_1^\ast\/$ (and this condition is also necessary). From Theorems~\ref{Bq/p} and \ref{
B_{p_0,p_1}} we can deduce   (using the property   $w\in B_p\Rightarrow w\in
B_{p-\epsilon}\/$) that, for
$1<p<\infty,\;0<q<\infty\/$,
\medskip

(a)  $u\in A^\ast_p\,\Rightarrow\,B_{q/p}\subset B_q(u)\/$,
 
(b) $B_{q/p}=B_q(u)\,\Rightarrow u\in A_p^\ast\/$.

\medskip
This cannot be improved (as it happens in the case $p=1\/$), that is, it is not true that the
characterization of the classes
$B_q(u)\/$ when 
$u\in A_\infty\/$, is 
$B_q(u)=B_{p/q}\/$ for  $p\/$ such that $u\in A^\ast_p\/$. In fact, it is easy to check that the
weight 
$u(x)=1+|x|,\;x\in\R,\/$ is in  $A_2^\ast\/$ and, by Theorem~\ref{te:principal} (see also Remark~
\ref{ob:disj}(ii)),
$w=\chi_{(0,1)}\in B_{3/2}(u)\/$ (condition  (\ref{cprincipal}) holds with $q=1\/$). However, 
$w\notin B_{3/4}\/$ and thus, $B_{3/2}(u)\neq B_{3/4}\/$.
\eobs

We shall study now a sufficient condition to have the weak-type boundedness. 
To this end, we shall need the following notation. We shall associate to each weight  $u \/$ in
$\Rn\/$ a function
$\Fiu\/$ which is connected with the family of functions $\big\{\fiq\big\}_Q\/$ introduced in 
(\ref{defphiQ}). This new function is defined by 
\bequ
\label{def-Phi_u}
\Fiu(t)=\sup_{Q}\fiq^\prime\big(u(Q)\,t\big),\quad0<t<\infty, 
\eequ
where the supremum is taken over all cubes $Q\subset\Rn\/$. Let us observe that 
$$
\Fiu(t)=\sup_{Q}{u(Q)\over|Q|}\big(u^{-1}\chi_Q\big)^\ast_u\big(u(Q)\,t\big),\quad0<t<\infty.
$$ 
Therefore, $\Fiu(t)=0\/$ if $t\ge1\/$. Moreover, by Proposition~\ref{6pphiQ}(iii),
$\displaystyle\int_0^t\Fiu(s)\,ds\ge t,\;0<t<1\/$. We know that if $u\equiv1,\;(Mf)_u^\ast(t)\approx
A(f_u^\ast)(t),\;t>0,\/$ where $A\/$ is the  Hardy operator. And although this does not happen when 
$u\ne1\/$ (see \cite{CS3}) there exists  the following  positive partial result  due
to  Leckband and Neugebauer (\cite{LN1}).

\bteo
\label{te:LN1}
Let $u\/$ be a weight in $\Rn\/$. Then, 
for each $f\in\M(\Rn)\/$ we have that 
$$
\big(Mf\big)^\ast_u(t)\le C\int_0^\infty\Fiu(s)f_u^\ast(st)\,ds,\quad0<s<\infty,
$$
where $C\/$ is a constant depending only on the dimension.
\eteo

Using this theorem,  we can find a sufficient condition to have the boundedness
$M:\Lambda^{p_0}_u(w_0)\to\Lambda^{p_1,\infty}_u(w_1)\/$. Observe the analogy with the
corresponding expressions of Proposition~\ref{ce7p}.

\bteo
\label{csamdnd}
 Let $0<p_1<\infty\/$.
\begin{enumerate}
\item[(a)] If  $1<p_0<\infty\/$ and there exists a constant $C<\infty\/$ such that 

(i)
$\displaystyle\bigg(\int_0^r\bigg(\int_0^{t/r}\Fiu(s)\,ds\bigg)^{p_0^\prime}W_0^{-p_0^\prime}(t)w_0(t)
\,dt\bigg)^{1/p_0^\prime}W_1^{1/p_1}(r)\le C,\quad r>0,$

(ii) $W_1^{1/p_1}(r)\le CW_0^{1/p_0}(r),\quad r>0,$

then $M:\Lambda^{p_0}_u(w_0)\to\Lambda^{p_1,\infty}_u(w_1)\/$ is bounded. 

\item[(b)] If $0<p_0\le1\/$, we have the same if the following condition holds
$$
\int_0^{t/r}\Fiu(s)\,ds\le C{W_0^{1/p_0}(t)\over W_1^{1/p_1}(r)},\quad0<t<r<\infty.
$$
\end{enumerate}
\eteo

\bdem
Let us consider the operator
$$
Tg(t)=\int_0^\infty\Fiu(s)g(st)\,ds=\int_0^\infty k(t,s)g(s)\,ds,
$$
(with $k(t,s)=(1/t)\Fiu(s/t)\/$) acting on functions $g\d\/$. By Theorem~\ref{te:LN1}, 
the boundedness of 
\bequ
\label{adT}
T\,:\,L^{p_0}_{\hbox{\decp dec}}(w_0)\longrightarrow L^{p_1,\infty}(w_1)
\eequ
implies the result. The characterization of  (\ref{adT})  can be obtained as a
direct application of  Theorem~4.3 in \cite{CS2}.
$\qquad\qed$
\edem

\bobs
\label{pobs}

(i) The sufficient condition of the previous theorem is not necessary in general. 
For example, if $w_0=w_1=w=\chi_{(0,1)},\;p_0=p_1=1\/$, such condition is
$\displaystyle\int_0^t\Fiu(s)\,ds\lesssim t,\;0<t<1,\/$ and by  Proposition~\ref{6pphiQ}, it is equivalent
to 
$u\in A_1\/$. That is,
$A_1\subset A_1(\chi_{(0,1)})\/$. However, for this weight,   condition (\ref{c7p}) is
necessary and sufficient in this case (see Remark~\ref{obs7p}). This condition is  
(Proposition~\ref{ce7p}b)
$\fiq(t)\lesssim\displaystyle{u(Q)\over W(u(Q))}W^{-1}(t),\;t>0,\/$ and, by Proposition~\ref{6pphiQ}(i), it 
is equivalent to 
$$
{W\big(u(Q)\big)\over|Q|}\le C{W\big(u(E)\big)\over|E|},\qquad E\subset Q.
$$
It is a simple exercise to check that (with $w=\chi_{(0,1)}\/$) the weight $u(x)=1+|x|,\;x\in\R,\/$ 
satisfies this condition. However,  $u\notin A_p\/$ if $p\le2\/$, and it follows that 
$A_1(\chi_{(0,1)})\neq A_1\/$, that is, the condition of Theorem~\ref{csamdnd} is not necessary.
 
(ii) Let us assume that the weight $u\/$ satisfies the following property: for every $r>0\/$  there
exists a cube 
$Q_r\/$ so that

\medskip
\ \ \ (a)\ \ $u(Q_r)\approx r\/$,

\medskip
\ \ \ (b)\ \ $\displaystyle\int_0^t\Fiu(s)\,ds\approx\displaystyle\int_0^t\phi_{Q_r}(u(Q_r)s)\,ds,\qquad
0<t<1\/$.

Then the condition of the previous theorem is equivalent to  (\ref{c7p}) and,
hence, both are equivalent to the boundedness of 
$M:\Lambda^{p_0}_u(w_0)\to\Lambda^{p_1,\infty}_u(w_1)\/$. This follows immediately from the
expressions of Theorem~\ref{csamdnd} and by  Proposition~\ref{ce7p} (with
$u_0=u_1=u\/$) since, fixed $r>0\/$  one can substitute, in the expressions of 
Theorem~\ref{csamdnd}, $r\/$ by
$u(Q_r)\/$ and the integral
$\displaystyle\int_0^{t/r}\Fiu(s)\,ds\/$ by $\displaystyle\int_0^{t/r}\fiq(u(Q_r)s)\,ds\/$ to obtain the
condition of  Proposition~\ref{ce7p} and viceversa.

(iii) Every power weight $u(x)=|x|^\alpha,\;x\in\Rn,\;\alpha>0\/$ satisfies the previous condition.
In fact, if $Q\/$ is a cube centered at the origin, one can easily see that
$$
\fiq(u(Q)t)\approx\Fiu(t)\approx t^{-\alpha\over n+\alpha},\qquad 0<t<1.
$$
Therefore, the characterization of the boundedness of 
$M:\lox{p_0}{u}{w_0}\to\lox{p_1,\infty}{u}{w_1}\/$
is given in this case by the expressions of Theorem~\ref{csamdnd} (which is now equivalent to 
(\ref{c7p}). See also Theorem~5.7.
\eobs

\section{Applications} 
\label{se:applic}

We can use now the result of the previous sections,  
 to characterize in its total generality, the boundedness of 
\bequ
\label{aLpq}
M:L^{p,q}(u)\to L^{r,s}(u),
\eequ
completing results of \cite{Muck1,  CHK, HK, L}.

\bteo Let $p,r\in(0,\infty),\;q,s\in(0,\infty]\/$.
\begin{enumerate}
\item[(a)] If either $p<1\/$,   $p\neq r\/$ or $s<q\/$, there are no weights $u\/$
 satisfying the boundedness   $M:L^{p,q}(u)\to L^{r,s}(u)\/$.

\item[(b)] The boundedness of 
$$
M\,:\,L^{1,q}(u)\longrightarrow L^{1,s}(u)
$$
only holds if  $q\le1,\;s=\infty\/$ and, in this case, a necessary and sufficient condition is 
$u\in A_1\/$.

\item[(c)] If $p>1\/$ and $0<q\le s\le\infty\/$, a necessary and sufficient condition 
to have the boundedness of 
$$
M\,:\,L^{p,q}(u)\longrightarrow L^{p,s}(u)
$$
is:
\begin{enumerate}
\item[(1)]  If $q\le1,\;s=\infty\,:\quad \displaystyle{u(Q)\over|Q|^p}\le C{u(E)\over|E|^p},
\quad E\subset Q\/$.

\item[(2)] If $q>1\/$ or   $s<\infty\,:\quad u\in A_p\/$.
\end{enumerate}
\end{enumerate}
\eteo

\bdem
(a) By Proposition~\ref{cnalpq}, the boundedness of 
\bequ
\label{alpqrs}
M\,:\,L^{p,q}(u)\longrightarrow L^{r,s}(u) 
\eequ
implies $s\ge q\/$. On the other hand, (\ref{alpqrs}) implies the boundedness of  $M:L^{p,q}(u)\to
L^{r,\infty}(u)\/$
which is  equivalent to $M:\lox{q}{u}{t^{q/p-1}}\to\lox{r,\infty}{u}{1}\/$ and from
Corollary~\ref{cor: cnad} it follows that 
\bequ
\label{u(Q)|Q|}
{\big(u(Q)\big)^{1/r}\over|Q|}\le C{\big(u(E)\big)^{1/p}\over|E|},\qquad E\subset Q.
\eequ
By the Lebesgue differentiation theorem, we have that necessarily  $p\ge1\/$. Moreover, applying
(\ref{u(Q)|Q|})  with $E=Q\/$ we obtain that $\big(u(Q)\big)^{1/r-1/p}\le C\/$ and since  $u(Q)\/$
can take any value in $(0,\infty)\/$ (Proposition~\ref{pro: cnu_0}), we get  $p=r\/$.

(b) Since the norm of a characteristic function in $L^{1,q}\/$ does not depend on  $q\/$, the
boundedness in (b) implies
$\|M\chi_E\|_{L^{1,s}(u)}\le C\|\chi_E\|_{L^1(u)},\;E\subset Q\/$. If $s<\infty\/$ we obtain, from
Theorem~\ref{apfciaq}, the boundedness of 
$M:L^p(u)\to L^p(u)\/$ with $p<1\/$ and, by (a), we have a contradiction.  Therefore, 
$s=\infty\/$. On the other hand, the boundedness of  $M:L^{1,q}(u)\to L^{1,\infty}(u)\/$ is the same
than 
$M:\lox{q}{u}{t^{q-1}}\to\lox{q,\infty}{u}{t^{q-1}}\/$ and by Corollary~\ref{cor: cnad} we obtain 
${u(Q)/|Q|}\le C{u(E)/|E|},\;E\subset Q,\/$ which is  $u\in A_1\/$. We know that this
condition is sufficient if 
$q\le1\/$. But from Corollary~\ref{Bpinfi(u)} we get $t^{q-1}\in B_{q,\infty}\/$,  and this is
only possible (see Theorem~\ref{te: I.6.5}) if $q\le1\/$.

(c) If $q\le s<\infty\/$, by Theorem~\ref{apfciaq} and using interpolation, the boundedness of 
 $M:L^{p,q}(u)\to L^{p,s}(u)\/$ is equivalent to the boundedness of $M:L^p(u)\to L^p(u)\/$ and a
necessary and sufficient condition is 
$u\in A_p\/$. In the case $M:L^{p,q}(u)\to L^{p,\infty}(u)\/$ (that is
$s=\infty\/$) we have two possibilities: (i) if $q>1\/$ a necessary and sufficient condition is 
(see \cite{CHK}) $u\in A_p\/$, and  (ii) if $q\le1\/$, from Corollary~\ref{cor: cnad} we obtain
the condition 
${u(Q)/|Q|^p}\le{u(E)/|E|^p}\/$. In \cite{CHK} it is shown that this condition is
sufficient in the case
$q=1\/$ and  (since
$L^{p,q}\subset L^{p,1}\/$ if $q\le1\/$) also in the case
$q\le1.\qquad\qed$
\edem

\bobs 

(i) The cases  $M:L^1(u)\to L^{1,\infty}(u)\/$ and $M:L^{p,q}(u)\to L^{p,s}(u)\/$ 
with  $1<p\le q\le s\le\infty\/$ were solved by  Muckenhoupt in \cite{Muck1}, giving rise to the 
$A_p\/$ classes. Chung, Hunt, and Kurtz (\cite{CHK, HK}) solved the case $M:L^{p,q}(u)\to L^{p,s}(u)\/$  with
$p>1,\;1\le q\le s\le\infty\/$ and  Lai (\cite{L}) proved the necessity of the condition  $p=r\/$
to have (\ref{aLpq}). 

(ii) The conditions obtained  in the cases of weak-type inequalities
($q,s<\infty\/$) coincide with  (\ref{c7p}).
\eobs

The following result characterizes the boundedness of 
\bequ
\label{bLpqu}
M:\lox{p_0}{u}{w_0}\to\lox{p_1,\infty}{u}{w_1} 
\eequ 
when $u(x)=|x|^\alpha,\;x\in\Rn\/$.

\bteo 
Let $u(x)=|x|^\alpha,\;x\in\Rn,\;\alpha>-n\/$.
\begin{enumerate}
\item[(a)] If $\alpha\le0\/$, {\rm(\ref{bLpqu})} is equivalent to $(w_0,w_1)\in B_{p_0,p_1,\infty}\/$.

\item[(b)] If $\alpha>0\/$    {\rm(\ref{bLpqu})} holds if and only if  $(\bar{w}_0,\bar{w}_1)\in
B_{p_0,p_1,\infty}\/$,
 where for each  $i=0,1\/$, 
$$
\bar{w}_i(t)=w_i\big(t^{n+\alpha\over n}\big)\,t^{\alpha\over n},\qquad t>0.
$$
\end{enumerate}
\eteo

\bdem
  (a) is consequence of Theorem~\ref{aua_1}, since  $u(x)=
|x|^\alpha\in A_1\/$ if $\alpha\le0\/$. To prove (b) we use the condition of 
Theorem~\ref{csamdnd}
 that  (by Remark~\ref{pobs}(iii)) is necessary and sufficient.  The final condition is obtained 
making the change of variables $t=\bar{t}^{(n+\alpha)/n},\;r=\bar{r}^{(n+\alpha)/n}\/$ 
and comparing the expression we get with the corresponding to the classes  
$B_{p_0,p_1,\infty}\/$ (Theorem~\ref{te: I.6.5}).$\qquad\qed$
\edem

\vfill\eject

\vfill\eject
\addcontentsline{toc}{chapter}{Index}
\markboth{INDEX}{M.J. Carro, J.A. Raposo, and J. Soria}

\input{CRS.ind}

\end{document}